\documentclass{amsart}
\usepackage{amsmath}
\usepackage{amsthm}
\usepackage{amssymb}
\usepackage{amsaddr}
\usepackage{mathrsfs}
\usepackage{enumerate}
\usepackage{graphicx}
\usepackage{fullpage}

\newcommand{\sQ}{\mathcal Q}
\newcommand{\sG}{\mathcal G}
\newcommand{\sL}{\mathcal L}
\newcommand{\sJ}{\mathscr J}
\newcommand{\p}[1]{\left( {#1} \right)}
\newcommand{\norm}[1]{\left\| #1 \right\|}
\newcommand{\sC}{\mathcal C}
\newcommand{\br}[1]{\left[ #1 \right]}
\newcommand{\Br}[1]{\left\{ #1 \right\}}
\newcommand{\stackvec}[2]{\br{\begin{array}{c} #1 \\ #2 \end{array}}}
\newcommand{\sX}{\mathscr X}
\newcommand{\sY}{\mathscr Y}
\newcommand{\sZ}{\mathscr Z}
\newcommand{\sV}{\mathscr V}
\newcommand{\dual}[1]{\left< #1 \right>}
\newcommand{\alphab}{\boldsymbol\alpha}
\newcommand{\xb}{\mathbf x}
\newcommand{\nb}{\mathbf n}
\newtheorem{theorem}{Theorem}
\newtheorem{lemma}{Lemma}
\newtheorem{corollary}{Corollary}
\newtheorem{proposition}{Proposition}

\newtheorem{remark}{Remark}

\title[Approximation of Operator-Valued Riccati Equations]{On the Approximation of Operator-Valued Riccati Equations in Hilbert Spaces}
\author{James Cheung}
\address{Independent Researcher in Los Angeles, CA, 90045}
\begin{document}

\begin{abstract}
In this work, we present an abstract theory for the approximation of operator-valued Riccati equations posed on Hilbert spaces. It is demonstrated here that the error of the approximate solution to the operator-valued Riccati equation is bounded above by the approximation error of the governing semigroup, under the assumption of boundedness on the semigroup and compactness on the coefficient operators. One significant outcome of this result is the correct prediction of optimal convergence for finite element approximations of the operator-valued Riccati equations for when the governing semigroup involves parabolic, as well as hyperbolic processes. We derive the abstract theory for the time-dependent and time-independent operator-valued Riccati equations in the first part of this work. In the second part, we derive optimal error estimates for the finite element approximation of the functional gain associated with model weakly damped wave and thermal LQR control systems. These theoretical claims are then corroborated with computational evidence. 
\end{abstract}

\maketitle


\tableofcontents

\section{Introduction}
In recent years, the reinforcement learning community has been rediscovering the usefulness of optimal control in the context of determining action policies. In general, reinforcement learning is involved with the direct numerical approximation of the dynamic programming principle, where the Hamilton-Jacobi-Bellman equation is solved to obtain the value of each possible state in the state space. A simplification of these equations can be made if the cost function to be minimized is quadratic and if the dynamical system governing the state is linear. This simplified problem is known in the optimal control community as the linear quadratic regulator (LQR) \cite{stengel1994optimal}. For these control systems, the dynamic programming principle gives an optimal feedback law which depends on the solution of Riccati equations. Aside from optimal controls, Riccati equations can be found in the K\'alm\'an filter \cite{kalman1960new}, in two-player linear quadratic games, and $H^\infty$-control \cite{MCMILLAN1993239}. The broad implications of Riccati equations make it a very important problem to study in science, engineering, and mathematics. 

Many application problems require solving a Riccati equation on infinite dimensional topologies. for example, in the control of systems whose state is governed by partial differential equations. We will call the Riccati equations studied in this work \emph{operator-valued Riccati equations} to emphasize its potentially infinite dimensional nature. The two forms of this equation studied in this work are the operator-valued Differential Algebraic Riccati Equation (DARE), given by 
$$
\left\{
\begin{aligned}
    \frac{d\Sigma(t)}{dt} &= A\Sigma(t) + \Sigma(t) A^* - \p{\Sigma BB^* \Sigma}(t) + C^*C \\
    \Sigma(0) &= \Sigma_0,
\end{aligned}
\right.
$$
and the operator-valued Algebraic Riccati Equation (ARE), given by
$$
    A\Sigma + \Sigma A^* - \Sigma BB^* \Sigma + C^*C = 0,
$$
where $A: \mathcal D(A) \rightarrow H$ is a generator of a $C_0$-semigroup $S(t)\in \sL\p{H}$, $B: U \rightarrow H$ and $C: H \rightarrow Y$ are bounded operators, and $H, U,$ and $Y$ are the Hilbert spaces associated with the application problem. A more complete discussion of this equation and the notation used here is given in \S\ref{section: problem setting}.

In practice, we can almost never obtain an analytical solution to an operator-valued Riccati equation. Instead, we rely heavily on numerical approximation techniques to approximate the solution to this equation with a finite dimensional approximation. An abstract representation of the approximate DARE is given by
$$
    \left\{
    \begin{aligned}
        \frac{d\Sigma_n(t)}{dt} &= A_n\Sigma_n(t) + \Sigma_n(t) A^*_n - \p{\Sigma_n B_nB^*_n \Sigma_n}(t) + C^*_nC_n \\
        \Sigma_n(0) &= \Sigma_{0,n},
    \end{aligned}
    \right.
$$
and an abstract approximate ARE is given by
$$
    A_n\Sigma_n + \Sigma_n A^*_n - \Sigma_n B_nB^*_n \Sigma_n + C^*_nC_n = 0,
$$
where $A_n$ denotes the generator of the approximate semigroup $\widehat S_n(t) \in \sL\p{H_n}$, $B_n: U\rightarrow H_n$ is the approximation to $B$, $C_n: H_n \rightarrow Y$ is the approximation to $C$, and $H_n$ is an approximating subspace of $H$. In the following paragraphs, we will provide a short literature review discussing the available approximation theory for operator-valued Riccati equations.

\subsection{Approximation Theory for Operator-Valued Riccati Equations}

The approximation theory for operator-valued Riccati equations starts with an investigation regarding whether its approximate solution converges as the approximation space $H_n$ approaches the Hilbert space $H$. In \cite{banks, banks1994approximation, ito1998approximation} a series of sufficient conditions were presented to guarantee the convergence of $\Sigma_n$ to $\Sigma$ in the operator norm. These classical convergence conditions are summarized as follows:
\begin{enumerate}
    \item (Semigroup Convergence) $S_n(t)\phi \rightarrow S(t)\phi$ for all $z \in H$ and uniformly for $t \in \mathbb R_+$. 
    \item (Adjoint Semigroup Convergence) $S_n^*(t) \phi \rightarrow S^*(t) \phi$ for all $z \in H$ and uniformly for $t \in \mathbb R_+$.
    \item $B_n u \rightarrow B u$ for all $u \in U$ and $B^*_n\phi \rightarrow B^* \phi$ for all $\phi \in H$.
    \item $C_n \phi \rightarrow C \phi$ for all $\phi \in H$ and $C^*_n y \rightarrow C^*y$ for all $y \in Y$.
\end{enumerate}
Additionally, for the approximation to the operator-valued ARE, we require the preservation of exponential stability (POES), i.e. that the approximate control system is exponentially stable, to ensure convergence of the approximate solution. These conditions are used extensively in the literature to demonstrate convergence of approximate solutions to operator-valued Algebraic Riccati Equations. For reference, some works using this approach are cited here \cite{borggaard2004strong, burns2013approximating, burns2001numerical, ito1998approximation}. While the theoretical framework developed in \cite{banks, banks1994approximation, ito1998approximation} is attractive for determining convergence of the approximate solutions, it does not provide an avenue for deriving error estimates. There are several works in the literature that provide error estimates for the approximate solution to operator-valued Riccati equations. In general, the applicability of these works are limited only to the case where $A: \mathcal D(A) \rightarrow H$ generates an analytic or compact semigroup, e.g. in the case of parabolic partial differential equations.

In \cite{ito1987strong}, a general inequality for deriving error estimates for the approximate solution of operator-valued ARE were derived, and from this formula, an error of $\mathcal O(h)$ was predicted for the finite element approximation for when $A: \mathcal D(A) \rightarrow H$ is the generator for a parabolic system. The error estimate implicitly requires that the classical convergence conditions presented in the previous paragraph are satisfied.

The error estimate for parabolic systems has been improved upon in  \cite{kroller1991convergence}, where an error bound of $\mathcal O\p{h^2\log(h^{-1})}$ was derived, where the error estimate applies to both the operator-valued ARE and DARE. The derivation of this estimate involves using the mild form of the DARE, i.e.
$$
    \Sigma(t)\phi := S(t) \Sigma_0 S^*(t)\phi + \int_0^t S(t-s)\p{C^*C - \p{\Sigma BB^*\Sigma}(s)} S^*(t-s)\phi ds \quad \forall \phi \in H,
$$
and in utilizing the mild form of the equation, the authors were able to utilize existing finite element error estimates for analytic semigroups in their analysis. However, the singularity that arises in parabolic problems with $L^2(\Omega)$ initial conditions adds a parasitic $t^{-1}$ term in the semigroup approximation error estimate. This in turn leads to the parasitic $\log\p{h^{-1}}$ term in the error estimate. 

Following this pattern of increasing the form strength the Riccati equation analyzed, \cite{burns2022optimal} derived an $\mathcal O\p{h^{k+1}}$ error estimate for parabolic systems using the Bochner integral form of the operator-valued ARE, i.e.,
$$
    \Sigma := \int_0^\infty S(s)\p{C^*C - \Sigma BB^*\Sigma} S^*(s) ds \quad \textrm{ on } H.
$$
The Brezzi-Rappaz-Raviart theorem was used to exploit fixed-point form of this equation in deriving error estimates. The derivation of error estimates in this context relies on specific compactness assumptions regarding the operators $C^*C$ and $BB^*$. The theory is also only applicable if $S(t)$ is a compact semigroup. 

An alternative theory for determining error estimates is presented in \cite[Chapter 4]{lasiecka2000control}, where the authors make the fundamental assumption of analyticity regarding the semigroup $S(t)$. This allows the semigroup approximation error to be bounded by the resolvent approximation error. This resolvent error estimate is then used to derive error estimates for the operator-valued ARE approximation. Higher-order error estimates, i.e. $\mathcal O\p{h^{k+1}}$, are generally not able to be derived from this theory due to the  fundamental limitations imposed by measuring the resolvent approximation error in the $\sL\p{H}$ norm. 

\subsection{Contributions Made in This Work}

In this work, we derive an abstract theory for the approximation of operator-valued Riccati equations. The theory is capable of predicting error estimates for when the process semigroup $S(t)$ is only bounded. This means, as motivating examples, that our results can be applied to processes governed by parabolic as well as hyperbolic partial differential systems without loss of generality. Additionally, our result only assumes pointwise convergence of the approximate semigroup $S_n(t)$ to $S(t)$ and compactness of $BB^*$ and $C^*C$. Hence our theory reduces the classical four sufficient conditions for convergence given in \cite{banks, banks1994approximation, ito1998approximation} into just two necessary conditions. 

In a summarized form, our theory determines that if
\begin{enumerate}
    \item $\norm{S(t)\phi - S_n(t)\phi}_H \rightarrow 0$ for all $\phi \in H$, and
    \item The coefficient operators $BB^*$ and $C^*C$ are compact operators belonging to $\sL\p{H,K}$, where $K$ is a compact subspace of $H$,
\end{enumerate}
then for the time-dependent problem
$$
    \norm{\Sigma(t) - \Sigma_n(t)}_{\sL\p{H}} \leq C \norm{S(t)-S_n(t)}_{\sL\p{K,H}}
$$
uniformly for all compact subsets of $t \in \mathbb R_+$, including $[0, \infty)$, and 
$$
    \norm{\Sigma - \Sigma_n}_{\sL\p{H}} \leq C \int_0^\infty e^{-\alpha t}\norm{S(t)-S_n(t)}_{\sL\p{K,H}} dt
$$
for the time-independent problem. Note that convergence of $\Sigma(t)$ to $\Sigma_n(t)$ is guaranteed because the pointwise convergence of linear operators implies uniform convergence on compact sets. This result then implies that the error associated with the solution approximation for the operator-valued Riccati equations is bounded only by the error of the semigroup approximation. This is significant because the \emph{existing approximation theory for time-dependent processes can be used to predict the convergence rate for the approximate solution to operator-valued Riccati equations}. As a corollary, optimal error estimates for various approximation methods can to be derived in our theory. This result then allows control system designers to choose which subspace $K$ that their weighting functions, sensors, and actuators map to and from to satisfy prespecified computational performance requirements for their system. Computational efficiency becomes especially important in sensor and actuator placement problems \cite{burns2009distributed, burns2010optimal, burns2015infinite, tang2017optimal, morris2016study} where many operator-valued Riccati equations needs to be solved in succession.

\subsection{Structure of the Paper}
In the remainder of this work, we will refer to the operator DARE and ARE as \emph{time-dependent} and \emph{time-independent} Riccati equations respectively. This is to avoid any confusion for readers not familiar with the implication of "differential" in the DARE, especially for when the "algebraic" part of the equation involves differential operators.

We will begin the main part of this work with a section describing the problem setting. We first present the necessary background material material in general functional analytic notation, linear operators, $C_0$-semigroups, and finally the abstract operator approximation. From there, we then define the Bochner integral form of the operator-valued Riccati equation in its time-dependent and time-independent forms and their associated abstract approximations. From there, we present several examples of real approximation methods that fall into our abstract approximation framework. We then motivate the practical importance of our work with several application problems. Finally, we conclude the section with a statement of the Brezz-Rappaz-Raviart theorem. 

Then in the next two sections, we move on to derive the main results in this work -- the abstract error estimates for the approximate solution to the time-dependent and time-dependent operator-valued Riccati equations. Here, we exploit the fixed-point form present in the Bochner integral representation to derive the aforementioned result. 

Finally, we apply our abstract approximation result to several model problems. These problems include: a simple scalar perturbation system, one and two-dimensional thermal control system, and a one-dimensional weakly damped wave system. With these prototypical examples, we demonstrate that our approximation theory can apply to simple systems as well as more complex systems governed by parabolic and hyperbolic partial differential equations. For each model problem, an error estimate is derived using our abstract results and then a computational example is presented to corroborate the derived error estimates.

\section{Problem Setting} \label{section: problem setting}
This section is concerned with establishing the theoretical framework necessary to analyze the approximation of operator-valued Riccati equations. We begin our discussion by defining the notation that we will use throughout this work. Then, we discuss the Bochner integral form of the time-dependent and time-independent operator Riccati equations and their approximations. Afterwards, we move on to motivate the practical implications our analysis by state several examples to show that the operator-valued Riccati equations are pervasive in optimal control, state estimation, and game theory. And finally, this section is concluded with a brief discussion of the Brezzi-Rappaz-Raviart theorem. 

\subsection{Notation}
In this subsection, we will define the notation that we will utilize throughout this work.

\subsubsection{Bounded Linear Operators on Hilbert Spaces}
In this work, we are interested in the Riccati equation posed on Hilbert spaces. To that end, let us denote $H$ as a separable Hilbert space and 
$$
\p{\cdot,\cdot}_H: H\times H \rightarrow \mathbb R_+
$$ 
as the inner product on $H$. The norm on $H$ is then defined as
$$
\norm{\phi}_H := \p{\phi,\phi}_H^\frac12
$$ 
for any $\phi \in H$. Further, we will denote $\sL\p{H}$ as the space of bounded linear operators on $H$ with the $\sL(H)$ norm is then defined as
$$
    \norm{T}_{\sL\p{H}} := \sup_{\phi \in H} \frac{\norm{T\phi}_H}{\norm{\phi}_H} \quad \forall \phi \in H,
$$
for any bounded linear operator $T \in \sL\p{H}$. 

The space of all compact operators on $H$ will be denoted as $\sJ_\infty$ to be consistent in notation with the trace-class spaces the solution of the operator-valued Riccati equations belong to. To be more specific, 
$$
    \sJ_\infty := \Br{T \in \sL(H): T\phi \in K, \ \ \forall \phi \in H \textrm{ and for any compact } K \subset H}.
$$
The $\sJ_\infty$ norm is known to coincide with the $\sL\p{H}$ norm for the operators that we discuss in this work. Although the $\sL(H)$ and $\sJ_\infty$ norms are equivalent, we will denote the norm of $\sJ_\infty$-valued operators as $\norm{\cdot}_{\sJ_\infty}$ to emphasize the compactness of the operator being measured. Additionally, it is known that $\sJ_\infty$ is a two-sided *-ideal in $\sL\p{H}$, which can be restated as the following.

\begin{proposition} \label{proposition: two-sided ideal}
Let $X \in \sJ_\infty$ and $Y \in \sL(H)$, then we have that $XY \in \sJ_\infty$, $YX \in \sJ_\infty$. 
\end{proposition}
\begin{proof}
See \cite[Theorem 10.69]{axler2020measure} for a proof. 
\end{proof}

Symmetric operators are prevalent in the study of operator-valued Riccati equations. We will denote $\sJ_\infty^s \subset \sJ_\infty$ as the space of all symmetric compact operators defined on $H$. In mathematical notation, 
$$
    \sJ_\infty^s := \Br{T \in \sJ_\infty : T\phi = T^*\phi \quad \forall \phi \in H}.
$$
It is then easy to see that $\sJ_\infty^s$ is a topological vector space. 

The space of bounded linear operators mapping a Hilbert space $H_1$ to a Hilbert space $H_2$ will be denoted as $\sL\p{H_1,H_2}$, whose norm is defined as
$$
    \norm{T}_{\sL\p{H_1,H_2}} := \sup_{\psi \in H_1} \frac{\norm{T\psi}_{H_2}}{\norm{\psi}_{H_1}},
$$
for any bounded linear operator $T \in \sL\p{H_1,H_2}$. For any bounded linear operator $T\in\sL\p{H_1,H_2}$, we have that $T^*\in\sL\p{H_2,H_1}$, where $T^*$ is defined as the adjoint of $T$. Additionally, it is known that 
$$
\p{T\psi,\eta}_{H_2} = \p{\psi,T^*\eta}_{H_1}
$$ 
for all $\psi\in H_1$ and $\eta \in H_2$. This then implies that $\norm{T}_{\sL\p{H_1,H_2}} = \norm{T^*}_{\sL\p{H_2,H_1}}$ \cite[Theorem 10.11]{axler2020measure}.

Let us now define $K \subset H$ as a compact subspace relative to the Hilbert space $H$. Throughout this work, we will utilize the following embedding property regarding the operator spaces $\sL\p{H,K}$ and $\sL\p{K,H}$.

\begin{proposition} \label{proposition: properties of compact operators}
Let $H$ be a Hilbert space and $K$ be a compact subspace of $H$. Then the following inclusions are satisfied
$$
\sL\p{H,K}\ \subset \sJ_\infty \subset \sL\p{K,H}.  
$$
\end{proposition}
\begin{proof}
We begin by proving the first assertion presented in the above statement. From basic definitions, we have that
\begin{equation*}
    \norm{T}_{\sJ_\infty} = \sup_{\phi \in H} \frac{\norm{T\phi}_H}{\norm{\phi}_H} 
    \leq \sup_{\phi \in H} \frac{\norm{T\phi}_K}{\norm{\phi}_H} = \norm{T}_{\sL(H,K)},
\end{equation*}
since $K \subset H$. This result implies that $\sL(H,K) \subset \sJ_\infty.$

We then see that 
\begin{equation*}
    \norm{T}_{\sL(K,H)} = \sup_{\phi \in K} \frac{\norm{T\phi}_H}{\norm{\phi}_K} \leq \sup_{\phi \in H} \frac{\norm{T\phi}_H}{\norm{\phi}_H} = \norm{T}_{\sL(H)} = \norm{T}_{\sJ_\infty},
\end{equation*}
again because $K \subset H$. This result then implies that $\sL(K,H) \subset \sJ_\infty$. This concludes the proof. 
\end{proof}
\begin{remark}
    In Lemma \ref{lemma: t1} of Appendix \ref{sec: technical appendix}, we provide a proof that $\sL(H,K)$ is also a two-sided *-ideal in $\sL(H)$. This is not a surprising result, since we already know that $\sL(H,K) \subset \sJ_\infty$.
\end{remark}

The final class of operator spaces of interest in this work are the continuous Bochner spaces $\sC(I; \sV)$, where $I \subset \mathbb R_+$ is a compact interval in the nonnegative real numbers and $\sV \subset \sL(H)$ is an arbitrary operator space. This space is characterized by
$$
    \sC(I; \sV) := \Br{ V(t) \in \sV \textrm{ for all } t \in I: \norm{V(t)}_{\sV} \textrm{ is a continuous function with respect to } t \in I }.
$$
The $\sC(I; \sV)$ norm is then defined as
$$
    \norm{V(\cdot)}_{\sC(I; \sV)} := \sup_{t \in I} \norm{V(t)}_{\sV}.
$$
With the necessary functional analytic notation defined, we are now ready to discuss $C_0$-semigroups and their approximations.

\subsubsection{$C_0$-Semigroups and Approximation on Hilbert Spaces} \label{subsubsec: semigroups}

Let $A: \mathcal D(A) \rightarrow H$, where $\mathcal D(A) \subset H$ is the domain of $A$, be the (infinitesimal) generator of a $C_0$-semigroup $S(t) \in \sL(H)$ (See \cite{goldstein2017semigroups} for a discussion on $C_0$-semigroups). For reference, $\mathcal D(A)$ is understood to be a dense subspace of $H$. By the definition of $C_0$-semigroups, we have that there exists constants $M, \omega \in \mathbb R_+$ so that 
\begin{equation} \label{eqn: semigroup bound}
    \norm{S(t)}_{\sL(H)} \leq M e^{\omega t} \quad \forall t \in \mathbb R_+. 
\end{equation}
It is then known that 
$$
    \frac{d}{dt}S(t)\phi = AS(t)\phi = S(t)A\phi, 
$$
for all $\phi \in \mathcal D(A)$ and $t \in \mathbb R_+$. This property is derived from the definition of the time-derivative, the semigroup property , i.e. $S(t+s) = S(t)S(s)$, and the fundamental definition of the generator of a semigroup. Given these properties, we see that $z = S(t) \phi \in H$ is the solution operator to the following initial value problem on Hilbert spaces
$$
\left\{
    \begin{aligned}
        \frac{dz(t)}{dt} &= Az \quad \textrm{ in } t\in I \\
        z(0) &= \phi.
    \end{aligned}
\right.
$$
This fact is useful in that it allows time-dependent infinite-dimensional systems to be cast into the familiar framework of initial value problems. We will utilize this observation many times in the analysis presented in this work. 

Denoting $A^*:\mathcal D(A^*) \rightarrow H$ as the adjoint of $A$, then \cite[Theorem 4.3]{goldstein2017semigroups} implies that $A^*$ is the generator of a semigroup $S^*(t) \in \sL(H)$ and that $S^*(t)$ is the adjoint of $S(t) \in \sL(H)$ generated by $A$. By extension, we have also that 
$$
    \frac{d}{dt} S^*(t) \phi = A^*S(t) \phi = S^*(t)A^* \phi
$$
for all $\phi \in \mathcal D(A^*)$ and that $q = S^*(t) \psi \in H$ is also the solution operator to an initial value problem
$$
\left\{
    \begin{aligned}
        \frac{dq(t)}{dt} &= A^*q \quad \textrm{ in } t\in I \\
        q(0) &= \psi.
    \end{aligned}
\right.
$$
for all $\psi \in H$,

We now define $H_n \subset H$ as a subset of $H$ and $\pi_n: H\rightarrow H_n$ as the orthogonal projection operator defined by
\begin{equation} \label{def: projection}
    \pi_n v := \sum_{j=1}^\infty \p{\phi_j, v}_H \phi_j,
\end{equation}
where $\Br{\phi_j}_{j=1}^\infty$ are a set of basis functions spanning $H_n$. We have assumed an infinite sum in \eqref{def: projection} because we generally do not need any assumptions on the dimensionality of $H_n$ in the abstract theory that we present in this work. If $H_n$ is a finite dimensional space, we can simply take $\phi_j := 0$ for all $j > n$. Let now $\iota_H: H_n \rightarrow H$ be the canonical injection operator. It we then have that
\begin{equation} \label{def: injection}
    \iota_H v_n:= \pi_n^* v_n \quad \forall v_n \in H_n.
\end{equation}
The projection and injection operators will play an important role in extending operators defined on $\sL(H_n)$ to $\sL(H)$, as we will see in the following paragraph and through the remainder of this section. 

Let now $A_n: \mathcal D(A_n) \rightarrow H_n$ be the generator of a $C_0$-semigroup $\widehat S_n(t) \in \sL(H_n)$. Again, by definition, there exists constants $M_n, \omega_n \in \mathbb R_+$ such that
\begin{equation} \label{eqn: Hn semigroup bound}
    \norm{\widehat S_n(t)}_{\sL(H_n)} \leq M_n e^{\omega_n t} \quad \forall t \in \mathbb R_+.
\end{equation}
We can now extend this semigroup to $\sL(H)$ by defining $S_n(t) := \iota_H \widehat S(t) \pi_n$. Since the projection and injection operators are bounded above by 1, it becomes clear that
\begin{equation} \label{eqn: semigroup approximation bound}
    \norm{S_n(t)}_{\sL(H)} \leq M_n e^{\omega_n t} \quad \forall t \in \mathbb R_+,
\end{equation}
where the constants $M_n, \omega_n \in \mathbb R_+$ used above are the same constants used in \eqref{eqn: Hn semigroup bound}. It should be made clear that $S_n(t) \in \sL(H)$ is not a $C_0$-semigroup on $H$. This is because $\lim_{t\rightarrow 0} S_n(t)\phi \neq \phi$ for all $\phi \in H$. It is however a $C_0$-semigroup on $H_n$. The definition of $S_n(t)\in\sL(H)$ becomes useful when measuring the distance between $S(t) \in \sL(H)$ and $\widehat S_n(t) \in \sL(H_n)$ in the $\sL(H)$ metric topology. 

Throughout this work, we will assume that $A_n: \mathcal D(A_n) \rightarrow H$ is defined in a way so that it generates a parameterized sequence of $C_0$-semigroup approximations $\Br{S_n(t)}_{n=1}^\infty \in \sL(H)$ satisfying
\begin{equation} \label{eqn: pointwise convergence}
    \lim_{n\rightarrow \infty} \norm{S(t)\phi - S_n(t)\phi}_{H} = 0 \quad \forall \phi \in H
\end{equation}
for any $t \in \mathbb R_+$. This then implies that
$$
\lim_{n \rightarrow \infty} \norm{S(t) - S_n(t)}_{\sL(K,H)} = 0
$$
for any compact subspace $K\subset H$ (as a consequence of \cite[Theorem 5.3-2]{ciarlet2013linear}). The sufficient and necessary condition to ensure \eqref{eqn: pointwise convergence} is the following.  
\begin{theorem}[Trotter-Kato \cite{ito1998trotter}] \label{theorem: Trotter-Kato}
Let $A\in \mathcal D(A) \rightarrow H$ be the generator of the $C_0$-semigroup $S(t) \in \sL\p{H}$ and $A_n\in \mathcal D(A_n) \rightarrow H$ be the generator of the approximation $\widehat S_n(t) \in \sL(H_n)$ and $\rho(A) \subset \mathbb C$ and $\rho(A_n) \subset \mathbb C$ denote the resolvent set of $A$ and $A_n$ respectively, then the following two statements are equivalent.
\begin{enumerate}
\item \label{item:a} There exists a $s\in \rho(A) \cap \bigcap_{n=1}^\infty \rho(A_n)$ such that 
$$
\lim_{n\rightarrow \infty} \norm{\p{A-s\mathcal I_H}^{-1}\phi - \iota_H \p{A_n - s\mathcal I_{H_n}}^{-1} \pi_n\phi}_H = 0
$$ 
for all $\phi \in H$, and
\item \label{item:b}$\lim_{n\rightarrow \infty} \norm{S(t)\phi-S_n(t)\phi}_H$ for all $\phi \in H$ and $t\in I$, where $S_n(t) := \iota_H \widehat S_n(t) \pi_n$. 
\end{enumerate}
If either \eqref{item:a} or \eqref{item:b} is true, then \eqref{item:a} holds true for any $\emph{Re}\p{s} > \omega$, for some $\omega \in \mathbb R$.  
\end{theorem}

We will primarily rely on Theorem \ref{theorem: Trotter-Kato} in this work to prove the pointwise convergence of the finite element approximation of the weakly damped wave semigroup in Lemma \ref{eqn: damped wave pointwise convergence}. 

\subsection{Bochner Integral Operator-Valued Riccati Equations}
The first problem that we wish to approximate the solution to is given by seeking a $\Sigma(\cdot) \in \sC(I; \sJ_\infty^s)$ that satisfies
\begin{equation} \label{eqn: Equation 1}
    \Sigma(t) = S(t) \Sigma_0 S^*(t) + \int_0^tS(s)\p{F - \Sigma G \Sigma}(s) S^*(s)ds, \quad\forall t \in I
\end{equation}
where $I:=[0,\tau]$ for some $\tau \in \mathbb R_+$, $\Sigma_0 \in \sJ_\infty^s$, $F(\cdot),G(\cdot) \in \sC\p{I, \sJ_\infty^s}$ are self-adjoint coefficient operators, and $S(t) \in \sL(H)$ is again a $C_0$-semigroup. This equation is shown to be equivalent to the weak differential operator Riccati equation found in, e.g. \cite[Equation 6.45]{curtain2012introduction}, as we demonstrate in Proposition \ref{proposition: equivalence}. 

From (\cite[Theorem 3.6]{burns2015solutions}), we have the following simplified statement of the well-posedness theorem for the Bochner integral representation of operator-valued Riccati equations.

\begin{theorem}[Burns-Rautenberg] \label{theorem: burns-rautenberg}
Let $H$ be a separable complex Hilbert space, $I=[0,\tau]$ be a bounded subset with $\tau \in \mathbb R_+$, and $S(t) \in \sL\p{H}$ be a $C_0$-semigroup on $H$. In addition, assume that 
\begin{enumerate}
\item $\Sigma_0 \in \sJ_\infty^s$;
\item $F\p{\cdot}, G\p{\cdot} \in \sC\p{I;\sJ_\infty^s}$, 
\end{enumerate}
where $F\p{t}$ and $G\p{t}$ are non-negative and self-adjoint on $H$ for all $t \in I$. Then the equation
$$
    \Sigma(t) = S(t) \Sigma_0 S^*(t) + \int_0^t S(t-s)\p{F - \Sigma G \Sigma}(s) S^*(t-s)ds,
$$
where the integral is understood to be a Bochner integral, has an unique solution $\Sigma(\cdot) \in \sC\p{I; \sJ_\infty^s}$. Furthermore $\Sigma(t)\in \sJ_\infty^s$ is non-negative and self-adjoint on $H$ for all $t \in I$.
\end{theorem}

The operator-valued Riccati equation, in most applications, is stated in its differential form. To justify the error analysis presented in this work, we must prove that the Bochner integral representation of the Riccati equation is equivalent to the differential form of the Riccati equation. We do so in the following.
\begin{proposition}[Equivalence of Bochner Integral and Differential Form] \label{proposition: equivalence}
Let $H$ be a separable Hilbert space and $S(t) \in \sL(H)$ be a $C_0$-semigroup on $H$. Then $\Sigma(\cdot) \in \sC\p{I, \sJ_\infty^s}$ satisfies 
\begin{equation}\label{eqn: a5}
    \Sigma(t) = S(t) \Sigma_0 S^*(t) + \int_0^t S(t-s) \p{F - \Sigma G\Sigma}(s) S^*(t-s)ds \quad \forall t \in I,
\end{equation}
if and only if $\Sigma(\cdot)$ satisfies
\begin{equation}\label{eqn: a6}
\left\{
\begin{aligned}
\frac{d}{dt} \p{ \psi, \Sigma(t) \phi}_H &= \p{\psi, A\Sigma(t)\phi}_H 
    + \p{\psi, \Sigma(t) A^*\phi}_H 
    - \p{\psi, (\Sigma G \Sigma)(t)\phi}_H
    + \p{\psi, F(t)\phi}_H\\
\p{\psi, \Sigma(0)\phi}_H &= \p{\psi, \Sigma_0\phi}_H,
\end{aligned}
\right.
\end{equation}
for all $\phi,\psi \in \mathcal{D}(A^*)$ and $t \in I$, where $I:= [0,\tau]$ is a compact time interval, $\Sigma_0 \in \sJ_\infty^s$, $F(\cdot),G(\cdot)\in\sC\p{I;\sJ_\infty^s}$ are non-negative and self-adjoint for all $t\in I$, and $A: \mathcal D(A) \subset H \rightarrow H$ is the generator of $S(t) \in \sL(H)$.
\end{proposition}
\begin{remark}
Equation \eqref{eqn: a6} is known as the \emph{weak form} of the time-dependent operator-valued Riccati equation.
\end{remark}
\begin{proof}
    Recall that \eqref{eqn: a5} is a well-posed problem by virtue of Theorem \ref{theorem: burns-rautenberg}. We begin this analysis by testing this integral equation with a test function $\phi \in \mathcal D(A^*)$. This yields the following mild form of the Riccati equation
\begin{equation} \label{eqn: z1}
    \Sigma(t)\phi = 
        S(t) \Sigma_0 S^*(t) \phi 
        + \int_0^t S(t-s)\p{F - \Sigma G \Sigma}(s) S^*(t-s)ds \quad \forall \phi \in \mathcal D(A^*),
\end{equation}
for every $t \in I$. It is clear that the solution of \eqref{eqn: a5} satisfies \eqref{eqn: z1}, since strong solutions satisfy weak formulations. We then differentiate with respect to $t\in I$ to see that
\begin{equation} \label{eqn: aaa1}
    \frac{d}{dt} \Sigma(t)\phi = 
        A\Sigma(t)\phi + 
        \Sigma A^*\phi - 
        \Sigma G \Sigma(t)\phi + 
        F(t)\phi \quad \forall \phi \in \mathcal D(A^*)
\end{equation}
for every $t \in I$. This is made clear by seeing that, for any $\phi \in \mathcal D(A^*)$, 
\begin{equation} \label{eqn: aaa2}
\begin{aligned}
    \frac{d}{dt}\br{S(t)\Sigma_0 S^*(t)\phi}&= \p{\frac{d}{dt}S(t)} \Sigma_0 S^*(t)\phi + S(t) \Sigma_0 \p{\frac{d}{dt}S^*(t) \phi} \\
      &= AS(t)\Sigma_0 S^*(t) \phi + S(t) \Sigma_0 S^*(t) A^* \phi \\
      &= A\Sigma(t)\phi + \Sigma A^*\phi \\
      & \quad - \int_0^t AS(t-s)\p{F - \Sigma G \Sigma}(s) S^*(t-s)\phi ds \\
      & \quad - \int_0^t S(t-s)\p{F - \Sigma G \Sigma}(s) S^*(t-s)A^*\phi ds \\
      &= A\Sigma(t)\phi + \Sigma A^*\phi - \int_0^t \frac{\partial}{\partial t}\br{S(t-s)\p{F - \Sigma G \Sigma}(s)S^*(t-s)}\phi ds,
\end{aligned}
\end{equation}
and that
\begin{equation} \label{eqn: aaa3}
\begin{aligned}
    &\frac{d}{dt}\int_0^t S(t-s) \p{F - \Sigma G \Sigma}(s) S^*(t-s)\phi ds = \\
    &\qquad F(t)\phi - \Sigma G \Sigma(t)\phi + \int_0^t \frac{\partial}{\partial t} \br{S(t-s)\p{F - \Sigma G \Sigma}(s) S^*(t-s)}\phi ds 
\end{aligned} 
\end{equation}
for all $t \in I$ after applying the Leibniz integral rule and seeing that $S(0)\phi = \phi$ for all $\phi \in   \mathcal D(A^*)$. Summing together \eqref{eqn: aaa2} and \eqref{eqn: aaa3}, and then taking the inner product with the test function $\psi \in \mathcal D(A^*)$ yields the following weak equation
\begin{equation} \label{eqn: aaa10}
    \frac{d}{dt} \p{\psi, \Sigma(t)\phi}_H = 
    \p{\psi, A\Sigma(t)\phi}_H + 
    \p{\psi, \Sigma A^*\phi}_H - 
    \p{\psi, \Sigma G \Sigma(t)\phi}_H + 
    \p{\psi, F(t)\phi}_H \quad \forall \phi, \psi \in \mathcal D(A^*).
\end{equation}
The weak form of the initial condition is then derived by evaluating \eqref{eqn: a5} at $t=0$ and testing the result against $\phi,\psi \in \mathcal D(A^*)$, i.e.
\begin{equation} \label{eqn: aaa4}
    \p{\psi, \Sigma(0)\phi}_H = \p{\psi, \Sigma_0 \phi}_H.
\end{equation}
This identity arises from the recalling that $S(0)\phi = S^*(0)\phi = \phi$ for all $\phi \in \mathcal D(A^*)$. The combination of \eqref{eqn: aaa10} and \eqref{eqn: aaa4} is equivalent to \eqref{eqn: a6}, and thus we have shown that the solution of \eqref{eqn: a5} also satisfies \eqref{eqn: a6}. 

We now prove the result in the other direction. Let us now suppose that $\Sigma(\cdot) \in \sC\p{I; \sJ_\infty^s}$ satisfies \eqref{eqn: a6}. Then integrating both sides of the equation gives us
\begin{equation*}
    \p{\psi, \Sigma(t)\phi}_H = \p{\psi, S(t) \Sigma_0 S^*(t) \phi}_H + \p{\psi, \int_0^t S(t-s)\p{F - \Sigma G \Sigma}(s) S^*(t-s)\phi ds}_H
\end{equation*}
for all $\phi, \psi \in \mathcal D(A^*)$. Theorem \ref{theorem: burns-rautenberg} demonstrates that the test functions $\phi, \psi \in \mathcal D(A^*)$ are not needed in order for the integral equation \eqref{eqn: a5} to be well-posed. In light of this, we see that the operator $\Sigma(\cdot) \in \sC(I;\sJ_\infty^s)$  satisfying \eqref{eqn: a6} also satisfies \eqref{eqn: a5}. This concludes the proof. 
\end{proof}

The well-posedness of the Bochner integral form of the time-dependent operator-valued Riccati equation and its equivalence to the differential form of the Riccati equation have thus been established. This result indicates to us that any result derived from the analysis of the Bochner integral form of the time-dependent equation also applies to the differential form of the equation. This is an important result because the analysis in the main parts of this work is dependent on the favorable fixed point structure present in the Bochner integral form of the equation that does not exist in its differential form.

It is known that, under the assumption of stabilizability or detectability of a system (depending on the context of the control system), there exists an unique non-negative steady-state solution to the time-dependent Riccati equation \cite[Part V, Theorem 3.1]{bensoussan2007representation}. Without going into the additional complexities associated with the concept of stabilizability and detectability, \emph{we simply assume in this work that $S(t)\in \sL(H)$ is an exponentially stable $C_0$-semigroup and that this semigroup is derived from an initial stabilization of the system}. Additionally, we will require that $F,G \in \sJ_\infty^s$ be constant in time. We formalize the definition of the steady-state Bochner integral Riccati equation and show that its solution is unique in the following.

\begin{proposition} \label{lemma: steady-state well-posedness}
Let $H$ be a separable Hilbert space and $S(t) \in \sL(H)$ be an exponentially stable semigroup such that there exists constants $M, \alpha \in \mathbb R_+$ satsifying
$$
    \norm{S(t)}_{\sL(H)} \leq Me^{-\alpha t} \quad \forall t \in \mathbb R_+.
$$
Additionally, assume that $F,G \in \sJ_\infty^s$ are non-negative self-adjoint operators. Then, the steady-state of \eqref{eqn: a5} is given by the solution of
\begin{equation} \label{eqn: a7}
    \Sigma_\infty = \int_0^\infty S(s)\p{F - \Sigma_\infty G \Sigma_\infty} S^*(s)ds.
\end{equation}
\end{proposition}
\begin{proof}
Under the assumptions given in the statement of this proposition, we know that there exists a unique steady-state for the weak form of the time-dependent operator-valued Riccati equation given by
\begin{equation}
\left\{
\begin{aligned}
\frac{d}{dt}\p{\psi,\Sigma(t)\phi}_H &= \p{\psi, A\Sigma(t)\phi}_H + \p{\psi, \Sigma(t) A^*\phi} 
- \p{\psi, (\Sigma(t) G \Sigma)(t)\phi}_H 
+ \p{\psi, F(t)\phi}_H \\
\p{\psi, \Sigma(0)\phi}_H &= \p{\psi,\Sigma_0\phi}_H,
\end{aligned}
\right.
\end{equation}
for all $\phi,\psi \in \mathcal D(A^*)$ and for any non-negative self-adjoint initial state $\Sigma_0 \in \sJ_\infty^s$ \cite[Part V, Proposition 2.2, 2.3]{bensoussan2007representation}. From the definition of the steady-state solution \cite[Part V, Proposition 2.1]{bensoussan2007representation}, we have that $\lim_{t\rightarrow \infty }\frac{d}{dt} \p{\psi, \Sigma(t)\phi}_H = 0$ for all $\phi,\psi \in \mathcal D(A^*)$. Denoting $\Sigma_\infty = \lim_{t\rightarrow \infty} \Sigma(t)$, we observe that the system
\begin{equation}
\left\{
\begin{aligned}
\frac{d}{dt}\p{\psi, \Sigma_\infty \phi}_H &= 
\p{\psi, A\Sigma_\infty \phi}_H 
+ \p{\psi, \Sigma_\infty A^* \phi}_H 
- \p{\psi, \Sigma_\infty F \Sigma_\infty\phi}_H 
+ \p{\psi,G\phi}_H \\
\p{\psi, \Sigma(0)\phi}_H &= \p{\psi, \Sigma_\infty\phi}_H,
\end{aligned}
\right.
\end{equation}
is stationary (i.e. $\frac{d}{dt}\Sigma_\infty = 0$ for almost all $t \in \mathbb R_+$) for all $\phi, \psi \in \mathcal D(A^*)$ if we choose $\Sigma_0 := \Sigma_\infty$. In other words the solution of the differential equation will remain at the steady-state if it starts at the steady state. 

We then utilize Proposition \ref{proposition: equivalence} to see that $\Sigma_\infty \in \sJ_\infty^s$ also satisfies
\begin{equation*}
\Sigma_\infty = S(t) \Sigma_\infty S^*(t) + \int_0^t S(t-s)\p{F - \Sigma_\infty G \Sigma_\infty} S^*(t-s)ds.
\end{equation*}
for any $t \in \mathbb R_+$. Taking the limit as $t \rightarrow \infty$ in the above equation results in 
\begin{equation*}
\begin{aligned}
\Sigma_\infty &= \lim_{t\rightarrow\infty} \int_0^t S(t-s)\p{F - \Sigma_\infty G \Sigma_\infty} S^*(t-s)ds \\
 &= \int_0^\infty S(s)\p{F - \Sigma_\infty G \Sigma_\infty} S^*(s)ds,
\end{aligned}
\end{equation*}
after a change of variables and seeing that $\lim_{t\rightarrow \infty}S(t) = 0$ since $\norm{S(t)}_{\sL(H)} \leq Me^{-\alpha t}$. 
\end{proof}


With this, we have demonstrated that the Bochner integral form of operator-valued Riccati equations are equivalent to its weak differential forms. Additionally, we have demonstrated that both the time-dependent and time-independent Bochner integral equations are well-posed. We now proceed to discuss the abstract approximation of these equations in the following subsection.

\subsection{An Abstract Approximation of Operator-Valued Riccati Equations}
The discussion presented in this subsection is aimed at demonstrating that the approach taken in the error analysis for the approximation of the Bochner integral form of the operator-valued Riccati equations in the following sections of this work has widespread implications. It is our goal here to demonstrate that many operator approximation methods can fit into the abstract theory that we present in this work. 

\subsubsection{Time-Dependent Abstract Approximation}
Recall that the time-dependent problem of interest is to seek a solution $\Sigma(\cdot) \in \sC(I; \sJ_\infty^s)$ that satisfies
$$
    \Sigma(t) = S(t)\Sigma_0 S^*(t) + \int_0^t S(t-s)\p{F - \Sigma G \Sigma}(s) S^*(t-s)ds,
$$
where $\Sigma_0 = \sJ_\infty^s$, $F(\cdot), G(\cdot) \in \sC\p{I; \sJ_\infty^s}$, and $S(t) \in \sL(H)$ is a $C_0$-semigroup. The approximate problem we are interested in analyzing involves seeking a solution $\Sigma_n(\cdot) \in \sC(I; \sJ_\infty^s)$ that satisfies
\begin{equation} \label{eqn: approximation 1}
    \Sigma_n(t) = S_n(t)\Sigma_0 S_n^*(t) + \int_0^t S_n(t-s)\p{F - \Sigma_n G \Sigma_n}(s) S_n^*(t-s)ds,
\end{equation}
where we have approximated $S(t) \in \sL(H)$ with $S_n(t) \in \sL(H)$. Since $S_n(t)$ is, in general, not a semigroup on all of $H$, we cannot directly apply Proposition \ref{proposition: equivalence} to demonstrate that \eqref{eqn: approximation 1} has a unique solution. However, we recall that $S_n(t) := \iota_H \widehat S_n(t) \pi_n$, where $\widehat S_n(t) \in \sL(H_n)$ is a semigroup on $H_n$. We prove that \eqref{eqn: approximation 1} is well-posed in the following by showing that its solution is the unique extension of the solution of the approximate operator-valued Riccai equation posed on $H_n$.

\begin{proposition} \label{proposition: approximation well-posedness}
    Let $H_n \subset H$ be a separable Hilbert space that represents the space of functions approximating the elements of $H$ and let $S_n(t) \in \sL(H)$ be the time-dependent approximation to the $C_0$-semigroup $S(t) \in \sL(H)$, where $S_n(t) = \iota_H \widehat S_n(t) \pi_n$ with $\widehat S_n(t) \in \sL(H_n)$ being the restriction of $S_n(t)$ to $\sL(H_n)$ (see \S\ref{subsubsec: semigroups}). Then there exists a unique solution $\Sigma_n(\cdot) \in \sC\p{I; \sJ_\infty^s}$ that satisfies \eqref{eqn: approximation 1} with $I:=[0,\tau]$ and $\tau \in \mathbb R_+$.
\end{proposition}
\begin{proof}
    We begin the analysis by first assuming that there exists an unique solution $\Sigma_n(\cdot) \in \sC\p{I, \sJ_\infty^s}$ that satisfies \eqref{eqn: approximation 1}. Given this assumption, there then exists an unique operator $\widehat \Sigma_n(t) \in \sJ_\infty^s\p{H_n}$ so that $\Sigma_n(t) = \iota_H \widehat \Sigma_n(t) \pi_n$. To see this, simply assume that there exists two operators $\widehat \Sigma_n^{(1)}, \widehat \Sigma_n^{(2)} \in \sL\p{H_n}$ such that $\Sigma_n(t) = \iota_H \widehat \Sigma_n^{(1)}(t) \pi_n$ and $\Sigma_n(t) = \iota_H \widehat \Sigma_n^{(2)}(t) \pi_n$. It then follows that $\iota_H \p{\widehat \Sigma_n^{(1)}(t) - \widehat \Sigma_n^{(2)}(t)} \pi_n = 0$, which implies that $\widehat \Sigma_n^{(1)}(t) = \widehat \Sigma_n^{(2)}(t)$. In the following, we demonstrate that $\widehat \Sigma_n(\cdot) \in \sJ_\infty^s(H_n)$ is the unique solution of the same operator-valued Bochner integral Riccati equation posed on $H_n$, which then validates the initial assumption of existence and uniqueness of the solution to \eqref{eqn: approximation 1},

    We take \eqref{eqn: approximation 1} and apply the projection operator $\pi_n$ and injection operator $\iota_H$ as follows
    \begin{equation}
    \pi_n \Sigma_n(t) \iota_H \pi_n = \pi_n S_n(t) \Sigma_0 S_n^*(t) \iota_H \pi_n + \int_0^t \pi_n S_n(t-s)\p{F - \Sigma_n G \Sigma_n}(s) S_n^*(t-s)\iota_H \pi_n ds
    \end{equation}
    on $H$. Placement of $\pi_n$ and $\iota_H$ in the integral is achieved by observing that it is not a time-dependent operator. From here, we first 
    observe that $\pi_n \Sigma_n(t) \iota_H \pi_n = \widehat \Sigma_n(t) \pi_n$ and that
    \begin{equation*}
    \begin{aligned}
        \pi_n S_n(t) \Sigma_0 S^*_n(t) \iota_H \pi_n &=\pi_n \iota_H \widehat S_n(t) \pi_n \Sigma_0 \iota_H \widehat S_n^*(t) \pi_n \iota_H \pi_n \\
        &= \pi_n \widehat S_n(t) \p{\pi_n \Sigma_0 \iota_H \pi_n}\widehat S_n^*(t) \pi_n
    \end{aligned}
    \end{equation*}
    for $t\in I$, after recalling that $\pi_n \iota_H = \pi_n$, $S_n(t) = \iota_H \widehat S_n(t) \pi_n$, and seeing that the range of $\widehat S_n^*(t)$ is $H_n$ implying that $\iota_H \widehat S_n^*(t) = \iota_H \pi_n \widehat S_n^*(t)$. With this, we have that 
    \begin{equation} \label{eqn: p1}
         \pi_n S_n(t) \Sigma_0 S^*_n(t) \iota_H \pi_n = \widehat S_n(t) \widehat \Sigma_{0,n} \widehat S^*_n(t)
    \end{equation}
    for $t \in I$ on $H_n$, where we have defined $\widehat \Sigma_{0,n} := \pi_n \Sigma_0 \iota_H \pi_n$, and observed that the restriction of the domain and range of the operator equation from $H$ to $H_n$ allows us to remove the $\pi_n$ operators from the equation. Following a similar strategy, we have that
    \begin{equation} \label{eqn: p2}
    \begin{aligned}
        \int_0^t \pi_n S_n(t-s)F(s) S^*_n(t-s) \pi_nds &= \int_0^t \pi_n \iota_H \widehat S_n(t-s)\p{\pi_n F(s) \iota_H \pi_n }\widehat S^*_n(t-s)\pi_nds\\
        &= \int_0^t \widehat S_n(t-s) \widehat F_n(s) \widehat S^*_n(t-s)ds,
    \end{aligned}
    \end{equation}
    on $H_n$ with $\widehat F_n(s) := \pi_n F(s)\iota_H\pi_n$. Finally, letting $\widehat G_n(t) := \pi_n G(t) \iota_H \pi_n$ and recalling that $\Sigma_n(t) := \iota_H \widehat \Sigma_n(t) \pi_n$, it then follows that 
    \begin{equation} \label{eqn: p3}
    \begin{aligned}
        &\int_0^t \pi_n S_n(t-s)\p{\Sigma_n G \Sigma_n}(s) S^*_n(t-s)\pi_n ds \\
        &\quad=  
        \int_0^t \pi_n \iota_H \widehat S_n(t-s) \pi_n \p{ \Sigma_n G \iota_H \Sigma_n}(s) \pi_n \iota_H \widehat S_n^*(t-s)\pi_n ds \\
        &\quad
        =\int_0^t \pi_n \iota_H \widehat S_n(t-s) \pi_n \br{\iota_H \widehat \Sigma_n \p{\pi_n G \iota_H \pi_n} \iota_H \widehat \Sigma_n \pi_n}(s) \pi_n \iota_H \widehat S^*_n(t-s) \pi_nds
        \\
        &\quad 
        =\int_0^t \widehat S_n(t-s) \p{\widehat \Sigma_n \widehat G_n \widehat \Sigma_n}(s) \widehat S^*_n(t-s) ds,
    \end{aligned}
    \end{equation}
    for $t \in I$ on $H_n$. With \eqref{eqn: p1}, \eqref{eqn: p2}, and \eqref{eqn: p3} we see that \eqref{eqn: approximation 1} becomes the following equation when restricting the domain and range of both sides of the operator equation to $H_n$
    \begin{equation} \label{eqn: restriction}
        \widehat \Sigma_n(t) = \widehat S_n(t) \widehat \Sigma_0 \widehat S^*_n(t) + \int_0^t \widehat S_n(t-s)\p{\widehat F_n - \widehat \Sigma_n \widehat G_n \widehat \Sigma_n}(s) S^*_n(t-s)ds \textrm{ on } H_n.
    \end{equation}
    
    Since $\widehat S_n(t)\in \sL(H_n)$, $\widehat \Sigma_0 \in \sJ_\infty^s(H_n)$, and $\widehat F_n(\cdot), \widehat G_n(\cdot) \in \sC\p{I, \sJ_\infty^s(H_n)}$, Proposition \ref{theorem: burns-rautenberg} implies that there exists a unique $\widehat \Sigma_n(\cdot) \in \sC\p{I; \sJ_\infty^s(H_n)}$ that satisfies \eqref{eqn: restriction}. Because $\Sigma_n(t) = \iota_H \widehat \Sigma_n(t) \pi_n$ for all $t \in I$, uniqueness of the solution $\widehat \Sigma_n(\cdot) \in \sC\p{I, \sJ^s_\infty(H_n)}$ to \eqref{eqn: restriction} then implies the uniqueness of the solution $\Sigma_n\p{\cdot} \in \sC\p{I, \sJ_\infty^s}$ to \eqref{eqn: approximation 1}.
\end{proof}

The above result then implies the solution to \eqref{eqn: approximation 1} also satisfies the following equation
\begin{equation} \label{eqn: approximation 2}
\left\{
\begin{aligned}
\frac{d}{dt} \p{ \psi_n, \Sigma_n(t) \phi_n}_H &= \p{\psi_n, A_n\Sigma_n(t)\phi_n}_H 
    + \p{\psi_n, \Sigma_n(t) A_n^*\phi_n}_H 
    - \p{\psi_n, (\Sigma_n G \Sigma_n)(t)\phi_n}_H
    + \p{\psi_n, F(t)\phi_n}_H\\
\p{\psi_n, \Sigma_n(0)\phi_n}_H &= \p{\psi_n, \Sigma_0\phi_n}_H,
\end{aligned}
\right.
\end{equation}
for all $t \in I$ and for all $\phi_n, \psi_n \in \mathcal D(A_n^*)$, where $A_n \in \mathcal D(A_n) \rightarrow H_n$ is again the generator of the approximate semigroup $\widehat S_n(t) \in \sL(H)$, with $S_n(t) := \iota_H \widehat S_n(t) \pi_n$. 

\subsubsection{Time-Independent Abstract Approximation}
Similarly, recall that the time-independent problem of interest is to seek a $\Sigma_\infty \in \sJ_\infty^s$ that satisfies
$$
    \Sigma_\infty = \int_0^\infty S(s)\p{F - \Sigma_\infty G \Sigma_\infty}S^*(s)ds,
$$
where $F, G \in \sJ_\infty^s$, and $S(t) \in \sL(H)$ is a $C_0$-semigroup. The correspondent approximation is then to seek a solution $\Sigma_{\infty, n} \in \sJ_\infty^s$ that satisfies
\begin{equation} \label{eqn: approximation 3}
    \Sigma_{\infty,n} = \int_0^\infty S_n(s)\p{F - \Sigma_{\infty, n} G \Sigma_{\infty,n}}S^*_n(s)ds,
\end{equation}
where $S_n(t) \in \sL(H)$ again approximates $S(t) \in \sL(H)$. The well-posedness of \eqref{eqn: approximation 3} follows from a similar argument used in Proposition \ref{proposition: approximation well-posedness}. Likewise, the differential form of the approximation is given by the following
\begin{equation} \label{eqn: approximation 5}
    \p{\psi_n, A_n \Sigma_{\infty, n}\phi_n}_H + \p{\psi_n, \Sigma_{\infty, n} A_n^*\phi_n}_H - \p{\psi_n, \Sigma_{\infty, n} G \Sigma_{\infty, n}\phi_n}_H + \p{\psi_n, F\phi_n}_H = 0 \quad \forall \phi_n, \psi_n \in \mathcal D(A_n^*),
\end{equation}
where $A_n: \mathcal D(A_n) \rightarrow H_n$ is again the generator to $\widehat S_n(t) \in \sL(H)$ with $S_n(t) := \iota_n S(t) \pi_n$.

\subsubsection{Galerkin Approximation} \label{sec: Galerkin Approximation}
The Galerkin approximation framework is likely the most utilized framework for approximating the solution of infinite dimensional equations with the solution of a finite dimensional problem. Its widespread popularity stems from its strong theoretical foundation and its flexibility in regards to choosing an approximating function space. Examples of popular Galerkin approximation methods are Finite Element methods \cite{ciarlet2002finite, strang1974analysis}, Discontinuous Galerkin methods \cite{cockburn2012discontinuous}, Isogeometric Analysis \cite{cottrell2009isogeometric}, Spectral Galerkin methods \cite{bernardi1997spectral}, and Fourier Transform methods \cite{vondvrejc2014fft}. Here, we provide a brief overview of how the Galerkin approximation can be applied to the problems studied in this work. 

The Galerkin approximation of a semigroup generator $A: \mathcal D(A) \rightarrow H$ is given by 
$$
    A_n^g := \pi_n^g A \iota_H\pi_n^g,
$$
where $\pi_n^g: H \rightarrow H_n^g$ is the orthogonal projection operator mapping $H$ into $H_n^g$. The subspace $H_n^g \subset H$ is chosen so that
\begin{equation} \label{eqn: projection galerkin}
    \p{\psi, \iota_H\pi_n^g A\iota_H \pi_n^g \phi}_H < \infty \quad \forall \phi, \psi \in H,
\end{equation}
where $\iota_H: H_n^g \rightarrow H$ is again the canonical injection operator. The Galerkin approximation of \eqref{eqn: Equation 1} is most commonly given in the weak form by
\begin{equation} \label{eqn: approximation 4}
\left\{
\begin{aligned}
\frac{d}{dt} \p{ \psi_n^g, \Sigma_n^g(t) \phi_n^g}_H 
    &= \p{\psi_n^g, A_n^g \Sigma_n^g(t)\phi_n^g}_H 
    + \p{\psi_n^g, \Sigma_n^g(t) (A_n^g)^*\phi_n^g}_H 
    - \p{\psi_n^g, (\Sigma_n^g G \Sigma_n^g)(t)\phi_n^g}_H
    + \p{\psi_n^g, F(t)\phi_n^g}_H\\
\p{\psi_n^g, \Sigma_n(0)\phi_n^g}_H &= \p{\psi_n^g, \Sigma_0\phi_n^g}_H,
\end{aligned}
\right.
\end{equation}
for all $\phi_n^g, \psi_n^g \in H_n^g$ and $t \in I$. It is then clear that \eqref{eqn: approximation 4} is a variant of the class of problems given by \eqref{eqn: approximation 2}, where $\mathcal D(A_n^g) = H_n^g$. 

Since \eqref{eqn: approximation 4} is posed on a finite-dimensional topology (of dimension $n$), we can represent this problem by seeking a $\widehat \Sigma_n^g(\cdot) \in  \sC\p{I; \mathcal \sJ_\infty^s\p{H_n}}$ that satisfies the following equation
\begin{equation} \label{eqn: approximation 6}
\left\{
\begin{aligned}
    \frac{d}{dt} \widehat\Sigma_n^g(t) &= A_n^g \widehat\Sigma_n^g(t) 
                                    + \widehat\Sigma_n^g(t) (A_n^g)^* 
        - (\widehat\Sigma_n^g \widehat G_n^g \widehat\Sigma_n^g)(t) + \widehat F_n^g(t) \\
    \widehat\Sigma_n^g(0) &= \Sigma_{0,n}^g,
\end{aligned}
\right.
\end{equation}
where $\Sigma_{0,n}^g := \pi_n^g \Sigma_0 \iota_H \pi_n^g$, $\widehat F_n^g(t) := \pi_n^g F(t) \iota_H  \pi_n^g$, and $\widehat G_n^g(t) := \pi_n^g G(t) \iota_H \pi_n^g$ are the Galerkin approximations to $\Sigma_0 \in \sJ_\infty^s$ and $F(\cdot), G(\cdot) \in \sC\p{I; \sJ_\infty^s}$ respectively. This becomes apparent after recalling \eqref{eqn: projection galerkin}, and seeing that the domain and range of $\Sigma_n^g(t)$ must be $H_n^g$ in order to satisfy \eqref{eqn: approximation 4}. Integrating \eqref{eqn: approximation 6} then gives us the following integral equation posed on $H_n^g$
\begin{equation} \label{eqn: qqq}
\widehat\Sigma_n^g(t) = \widehat S_n^g(t) \Sigma_{0,n}^g (\widehat S_n^g)^*(t) + \int_0^t \widehat S_n^g(t-s)\p{\widehat F_n^g - \widehat\Sigma_n^g \widehat G_n^g \widehat\Sigma_n^g}(s) (\widehat S_n^g)^*(t-s)ds.
\end{equation}
Let us now define $\Sigma_n^g(t):= \iota_H \widehat \Sigma_n^g(t) \pi_n^g$ and $S_n^g(t) := \iota_H \widehat S_n^g(t) \pi_n^g$. It then follows from inspection, by injecting the range of \eqref{eqn: qqq} into $H$, that $\Sigma_n^g(t) \in \sC(I; \sJ_\infty^s)$ satisfies 
\begin{equation} \label{eqn: www}
\Sigma_n^g(t) = S_n^g(t) \Sigma_0 (S_n^g)^*(t) + \int_0^t S_n^g(t-s)\p{F - \Sigma_n^g G \Sigma_n^g}(s) (S^g_n)^*(t-s) ds,
\end{equation}
from which it becomes clear that \eqref{eqn: www} satisfies the abstract form of the approximate problem given in \eqref{eqn: approximation 3}. Therefore, the Galerkin approximation of the time-dependent approximation of the operator-valued Riccati equation can be analyzed within the context of the abstract theory presented in this work. A similar argument applies to the Galerkin approximation variant of the time-independent problem given in \eqref{eqn: approximation 4}.  

\subsubsection{Asymptotic Approximation}
Oftentimes, the generator $A: \mathcal D(A) \rightarrow H$ is parameterized by a presumably small parameter $\epsilon \in \mathbb R_+$. The operator perturbation theory \cite{kato2013perturbation} tells us that this operator can be represented in the following asymptotic series
\begin{equation*}
    A = \sum_{j=0}^\infty \epsilon^j A_j,
\end{equation*}
where each $A_j$ must necessarily map $\mathcal D(A)$ into $H$. Truncating the sum to $n \in \mathbb N_0$ terms gives us the asymptotic approximation
\begin{equation*}
    A_n^a := \sum_{j=0}^n \epsilon^j A_j.
\end{equation*}
Let now $S(t)\in \sL(H)$ and $S_n^a(t) \in \sL(H)$ be the semigroups generated by $A$ and $A_n^a$ respectively. Then, the conditions laid out in the semigroup approximation theory presented in \cite[Chapter 1, \S7]{goldstein2017semigroups} must be satisfied in order for $S_n^a(t) \rightarrow S(t)$ as $n \rightarrow \infty$. If convergence can be established in the limit, then as an approximation, we can simply replace $S(t) \in \sL(H)$ with $S_n^a(t) \in \sL(H)$ in \eqref{eqn: Equation 1} without any issue because the domain and range of $S(t)$ and $S_n^a(t)$ coincide. This describes the process of the asymptotic approximation of operator-valued Riccati equation. As an example, The application of asymptotic operator approximation can be found in applications involving the Born first order approximation for wave scattering in RADAR and quantum mechanical systems \cite{cheney2009fundamentals, gubernatis1977born}.

\subsection{Application Problems}
This subsection will describe how the discussion presented in the preceding discussion may be applied to many important problems of practical interest. This framework may be applied to the K\'alm\'an filter, Linear Quadratic Regulator, and Two-Player Differential Games, among other potential applications not described here (e.g. $H^\infty$-optimal control). In the following discussion we will first pose the application problem, highlight its use of a Riccati equation, and finally demonstrate how the Riccati equation fits into the abstract framework described in this section. For each of the application problems, we will present the associated time-dependent Riccati equation in the weak form and will defer the derivation of the steady-state form of these equations to the reader. The reader is reminded that the differential and Bochner integral form of the operator-valued Riccati equation are shown to be equivalent in Proposition \ref{proposition: equivalence}.

\subsubsection{K\'alm\'an Filter}
The K\'alm\'an filter \cite{kalman1960new} has become the predominant method used in the field of state estimation. It has been so successful that state estimation practitioners have made many attempts to extend its application to estimating the state of nonlinear systems \cite{ribeiro2004kalman}. Here, we describe the K\'alm\'an filter for linear-quadratic systems. 

First, let us define $H, X,$ and $Y$ as separable Hilbert spaces associated with the state, noise, and sensor output respectively. Then the cost functional is defined as
$$
    \mathcal J(\hat z) := \mathbb E\br{\int_0^\tau \norm{z(t; \eta(t)) - \hat z(t; \eta(t), \nu(t))}^2_H dt},
$$
where $\mathbb E(\cdot)$ denotes the expectation operator with respect to the Gaussian measures associated with the random processes $\eta(t)\in X$ and $\nu(t) \in Y$ representing the process and measurement noise with covariance operators $R_\eta \in \sL(X)$ and $R_\nu \in \sL(Y)$ respectively, and $z(\cdot), \hat z(\cdot) \in L^2([0,\tau],H)$ is the solution of the following set of stochastic differential equations (stated here with an abuse of notation)
\begin{equation*}
\left\{
\begin{aligned}
    \frac{dz}{dt} &= Az(t) + B\eta(t), \quad &z(0) = z_0 \\
    y(t) &= Cz(t) + \nu(t) & \\
    \frac{d\hat z}{dt} &= A\hat z(t) + L_t\p{y(t) - C\hat z}, \quad & \hat z(0) = \hat z_0
\end{aligned}
\right.
\end{equation*}
for all $t \in [0, \tau]$, where $A: \mathcal D(A) \rightarrow H$ is a generator of a $C_0$-semigroup, $B: X\rightarrow H$ is the state noise operator, $C: H \rightarrow Y$ is the sensing operator. The functional gain $L_t: Y \rightarrow H$ updates the state estimate $\hat z(t)$ based on the difference between the sensor data $y(t)$ and predicted sensor output $C\hat z(t)$ for all $t \in [0, \tau]$. This operator is given by
$$
    L_t= \Sigma(t) R^{-1}_\nu C,
$$
where $\Sigma(t) \in \sL(H)$ is the covariance operator given by the solution of the following operator-valued Riccati equation
\begin{equation*}
\left \{
\begin{aligned}
    \frac{d}{dt}\p{\psi, \Sigma(t) \phi}_H &= 
    \p{\psi, A\Sigma(t)\phi}_H 
    + \p{\psi, \Sigma(t) A^* \phi}_H
    - \p{\psi, \p{\Sigma C^* R^{-1}_\nu C\Sigma}(t)\phi}_H 
    + \p{\psi, \p{BR_\eta B}(t) \phi}_H \\
    \p{\psi, \Sigma(0) \phi}_H &= \p{\psi, \Sigma_0 \phi}_H
\end{aligned}
\right. 
\end{equation*}
for all $t \in [0, \tau]$ and for all $\phi, \psi \in \mathcal D(A)$. It then becomes clear that the above equation is a variation of \eqref{eqn: a6} after setting $F:= BR_\eta B^*$ and $G:= C^*R_\nu^{-1}C$. 
 
\subsubsection{Linear Quadratic Regulator} \label{sec: LQR}
The linear quadratic regulator (LQR) \cite{stengel1994optimal} is a fundamental closed-loop control method that forces the state of a linear system to the desired terminal state. Like the K\'alm\'an filter, the LQR has been used extensively in the control of many linear systems. It has also been applied to nonlinear systems, where the state system has been linearized about a desired set point \cite{tedrake2009lqr, boby2014robust, prasad2014optimal}.

Let $H$, $Y$, and $U$ be separable Hilbert spaces associated with the state, weighting operator, and control respectively. The cost function to be minimized is a quadratic functional given by 
$$
    \mathcal J(u) := \p{z(\tau;u), \Pi_\tau z(\tau;u)}_H + \int_0^\tau \br{\norm{C(t)z(t;u)}_Y^2 + \p{u(t), R u(t)}_U} dt,
$$
where $\Pi_\tau\in \sJ_\infty^s, C(t) \in \sL(H, Y)$ and $R \in \sL(U)$ are weighting operators whose values can be tuned to adjust the importance of minimizing the terminal state, running cost, and the control effort respectively in the cost functional. The system to be controlled is given by 
$$
\left\{
\begin{aligned}
    \frac{dz}{dt} &= \mathcal Az(t) + B(t)u(t) \\
    z(0) &= z_0
\end{aligned}
\right.
$$
for all $t \in [0, \tau]$, where $B(t): U \rightarrow H$ is the actuation operator. The optimal control, given as a function of $z(t)$, is the following
$$
    u_{opt}(t) = -\mathcal K_tz(t),
$$
where $\mathcal K_t: H \rightarrow U$ is the functional gain given by 
$$
    \mathcal K_t = R^{-1}(t)B(t)^*\Pi(t),
$$
and $\Pi(t) \in \sJ_\infty^s$ is the solution of the following operator-valued Riccati equation
$$
\begin{aligned}
    \frac{d}{dt}\p{\psi,\Pi(t)\phi}_H &= 
    -\p{\psi, \mathcal A^*\Pi(t)\phi}_H 
    -\p{\psi, \Pi(t)\mathcal A\phi}_H 
    +\p{\psi, \p{\Pi BR^{-1}B^*\Pi}(t)\phi}_H 
    -\p{\psi, (C^*C)(t)\phi}_H \\
    \p{\psi, \Pi(\tau) \phi}_H &= \p{\psi, \Pi_\tau\phi}_H
\end{aligned}
$$
for all $t \in [0, \tau]$ and for every $\phi, \psi \in \mathcal D(A)$. A change of variables $t \mapsto s$, where $s=\tau-t$ transforms the terminal value problem into the initial value problem that we seek in the abstract description of the operator-valued Riccati equation. Taking $A := \mathcal A^*$, $F(t) := C^*C(t)$, and $G(t):= \p{BR^{-1}B^*}(t)$ then allows us to see that the Riccati equation associated with LQR control can be analyzed within the framework described in this work.  

\subsubsection{Linear Quadratic Two-Player Differential Games}
In a linear quadratic two-player differential game \cite{bernhard1979linear}, there exists two sets of strategies $u_1(t) \in U_1$ and $u_2(t) \in U_2$, associated with player 1 and player 2 respectively, that aims to find the saddle-point of the following cost functional
$$
    \mathcal J(u_1, u_2) := \p{z(\tau; u_1, u_2), P_\tau z(\tau; u_1, u_2)}_H + \int_0^\tau \br{\p{z(t; u_1, u_2), Q z(t; u_1, u_2)}_H + \norm{u_1(t)}_{U_1}^2 - \norm{u_2(t)}_{U_2}^2} dt.
$$
The dynamical system affected by these competing strategies is given by
$$
\left\{
\begin{aligned}
    \frac{dz}{dt} &= \mathcal Az(t) + B_1 u_1(t) + B_2 u_2(t) \\
            z(0) &= z_0.
\end{aligned}
\right.
$$
The optimal strategy for each player is then determined by
$$
    u_{1,opt}(t) = -B_1 P(t) z(t)  \qquad 
    u_{2,opt}(t) = B_2^* P(t) z(t),
$$
where $P(t) \in \sL(H)$ is the solution of the following operator-valued Riccati equation:
$$
\begin{aligned}
    \frac{d}{dt}\p{\psi,P(t)\phi}_H &= 
    -\p{\psi, \mathcal A^*P(t)\phi}_H 
    -\p{\psi, P(t)\mathcal A\phi}_H 
    +\p{\psi, \p{PRP(t)}\phi}_H 
    -\p{\psi, P\phi}_H \\
    \p{\psi, P(\tau) \phi}_H &= \p{\psi, P_\tau\phi}_H
\end{aligned}
$$
for every $t \in [0, \tau]$ and $\phi, \psi \in \mathcal D(A)$, where $R= B_1B_1^* - B_2B_2^*$. Like in the LQR Riccati equation presented in the previous example, our desired initial-value problem is obtained through a variable transformation. Of course, setting $A:= \mathcal A^*$, $G:= R$, and $F:= Q$ allows us to see that the Riccati equation associated with linear quadratic two-player games fits within the abstract framework of this work. 

\subsection{The Brezzi-Rappaz-Raviart Theorem} \label{subsection: BRR}

The Brezzi-Rappaz-Raviart (BRR) Theorem \cite[Part IV, Theorem 3.3]{girault2012finite} is an important result in numerical analysis that provides a framework for deriving error bounds for the approximation of certain classes of nonlinear equations. We demonstrate in the following sections that approximations to both the time-dependent and time-independent operator-valued Riccati equation can be analyzed within this framework. A brief statement of the problem setting and the BRR theorem itself is presented in the following paragraph. 

Consider the fixed point problem
\begin{equation} \label{eqn: FPP}
u = T\sG(\lambda, u),
\end{equation}
where $u \in \sX$, $\lambda \in \Lambda$ belongs to  a compact interval in $\mathbb R_+$, $T \in \sL(\sY; \sX)$ is a linear operator, and $\sG: \Lambda \times \sX \rightarrow \sY$ is a $C^2$ smooth mapping. Let us introduce $T_n\in \sL(\sY; \sX)$ as an operator approximation to $T \in \sL\p{\sY;\sX}$. We then have that \eqref{eqn: FPP} can be approximated by the following:
\begin{equation}\label{eqn: approximate FPP}
	u_n = T_n\sG(\lambda, u_n).
\end{equation}
Now suppose that there exists another space $\sZ \hookrightarrow \sY$ with continuous imbedding such that
\begin{equation} \label{eqn: perturbation condition}
	D_u \sG(\lambda, u) \in \sL(\sX; \sZ ) \quad  \forall \lambda \in \Lambda \textrm{ and } u \in \sX.
\end{equation}
Then, additionally under the convergence assumptions
\begin{equation} \label{eqn: pointwise convergence condition}
	\lim_{n \rightarrow \infty} \norm{\p{T - T_n}  g}_{\sX} = 0 \quad \forall g \in \sY
\end{equation}
and
\begin{equation} \label{eqn: uniform convergence condition}
	\lim_{n \rightarrow 0} \norm{T - T_n}_{\sL(\sZ; \sX)} = 0
\end{equation}
the following holds.
\begin{theorem}[Brezzi--Rappaz--Raviart] \label{theorem: BRR}
Assume that conditions \eqref{eqn: perturbation condition}, \eqref{eqn: pointwise convergence condition}, and \eqref{eqn: uniform convergence condition} are satisfied. In addition, assume that $\sG: \Lambda \times \sX \rightarrow \sY$ is a $\mathcal C^2$ operator, with $D^2\sG$ bounded on all subsets of $\Lambda \times \sX$. Then, there exists a neighborhood $\mathcal O$ of the origin in $\sX$, and for $n  \geq N$ large enough, a unique $\mathcal C^2$--function $\lambda \in \Lambda \rightarrow u_n(\lambda) \in \sX$ such that
$$
	\p{\lambda, u_n(\lambda); \lambda \in \Lambda}
		\textrm{ is a branch of nonsingular solutions of } \eqref{eqn: approximate FPP}
$$
and
$$
	u(\lambda) - u_n(\lambda) \in \mathcal O \quad \forall \lambda \in \Lambda.
$$
Furthermore, there exists a constant $C\in \mathbb R_+$ independent of $n$ and $\lambda$ with:
$$
	\norm{u(\lambda) - u_n(\lambda)}_{\sX} \leq C \norm{(T - T_n)\sG(\lambda, u(\lambda))}_{\sX} \quad \forall \lambda \in \Lambda.
$$
\end{theorem}

In the next two sections, we utilize this theorem to demonstrate that the error estimate for the approximate solution to the operator-valued Riccati equation, in principle, only depends on the the error induced by the approximate semigroup over the compact subspace characterized by the range of the coefficient operators in \eqref{eqn: Equation 1} and \eqref{eqn: a7}.

\section{An Error Estimate for Time-Dependent Approximations} \label{section: time-dependent approximation}
In this section, we will analyze the convergence of the approximate solutions of time-dependent operator-valued Riccati equations. Recall that on $t \in I$, where $I := [0,\tau]$, we have defined the Bochner integral form of the problem as
\begin{equation} \label{eqn: time-dependent problem}
    \Sigma(t) = S(t) \Sigma_0 S^*(t) + \int_0^t S(t-s)\p{F - \Sigma G \Sigma}(s) S^*(t-s)ds,
\end{equation}
where $F(\cdot),G(\cdot) \in \sC\p{I; \sJ_\infty^s}$ and $\Sigma_0 \in \sJ_\infty^s$ are self-adjoint and nonnegative operators on $H$. We wish to seek an approximation to \eqref{eqn: time-dependent problem} by solving the following approximating equation 
\begin{equation} \label{eqn: approximate time-dependent problem}
    \Sigma_n(t) = S_n(t) \Sigma_0 S^*_n(t) + \int_0^t S_n(t-s)\p{F - \Sigma_n G \Sigma_n}(s) S^*_n(t-s)ds,
\end{equation}
where $S(t)\in \sL\p{H}$ is approximated with $S_n(t) \in \sL\p{H}$.

In the following analysis, we will derive a fundamental estimate that will allow for the derivation of optimal convergence rates of $\Sigma_n(t) \rightarrow \Sigma(t)$ in the $\sJ_\infty$ topology. The convergence rate will depend principally on which compact subspace $K$ of $H$ the coefficient operators $\Sigma_0, F(t),$ and $G(t)$ map to. \emph{To make this point clear in the analysis, we will take $\Sigma_0, F(t), G(t) \in \sL\p{H,K} \cap \sJ_\infty^s$.} The choice of the compact space $K$ will remain abstract and will vary depending on the application problem and system design decisions. 

The analysis presented in the remainder of this section seeks to prove that the equations \eqref{eqn: time-dependent problem} and \eqref{eqn: approximate time-dependent problem} have the necessary nonlinear structure necessary to utilize the Brezzi-Rappaz-Raviart theorem (Theorem \ref{theorem: BRR}) so that we can bound the error of the approximate solution to the exact solution by the error associated with the way the semigroup in \eqref{eqn: approximate time-dependent problem} is approximated. To be more descriptive, we aim to prove that there exists a positive constant $C \in \mathbb R_+$ so that
$$
\norm{\Sigma(\cdot) - \Sigma_n(\cdot)}_{\sC\p{I,\sJ_\infty}} \leq C \sup_{t\in I} \norm{S(t) - S_n(t)}_{\sL\p{K,H}}.
$$
The following analysis will begin with an inspection of the linear operator from which \eqref{eqn: time-dependent problem} and \eqref{eqn: approximate time-dependent problem} are perturbations of. We then conclude our analysis by utilizing the derived properties of this operator to prove our desired result with the Brezzi-Rappaz-Raviart theorem (Theorem \ref{theorem: BRR}).

\subsection{Time-Dependent Lyapunov Operator}
Let us define $\Phi_t: \sJ_\infty^s \times \sC\p{I;\sJ_\infty^s} \rightarrow \sC\p{I; \sJ_\infty^s}$ as 
\begin{equation} \label{eqn: lyapunov operator definition}
    \Phi_t\stackvec{X}{Y(\cdot)} := S(t)XS^*(t) + \int_0^t S(t-s)Y(s)S^*(t-s)ds
\end{equation}
where $X\in \sJ_\infty^s$ and $Y(\cdot)\in \sC\p{ I;\sJ_\infty^s}$. 

\begin{remark}
We call $\Phi_t: \sJ_\infty^s \times \sC\p{I;\sJ_\infty^s} \rightarrow \sC\p{I;\sJ_\infty^s}$ the \emph{Lyapunov operator}, since it is the solution operator for the time-dependent variant of the Lyapunov equation shown below
$$
\left\{
\begin{aligned}
    \frac{dX}{dt} &= AX + XA^* \\
               X(0) &= X_0.
\end{aligned}
\right.
$$
This equation arises from the study of the stability of linear dynamical systems \cite{behr2019solution}. Though they are related, $\Phi_t(\cdot)$ should not be confused with the \emph{Lyapunov function} \cite[\S 4.4]{aastrom2021feedback}, often denoted $V(\cdot) : H \rightarrow \mathbb R_+$, used to determine the stability of control systems. 
\end{remark}

Let us now define $\sQ(\cdot): \sC\p{I; \sJ_\infty^s} \rightarrow \sL\p{H,K}\cap \sJ_\infty^s \times \sC\p{I; \sL\p{H,K} \cap \sJ_\infty^s}$ as 
\begin{equation} \label{eqn: Q def}
    \sQ(X(\cdot)) := \stackvec{\Sigma_0}{F(\cdot) - (X G X)(\cdot)},
\end{equation}
where we again take $\Sigma_0 \in \sL\p{H,K}\cap\sJ_\infty^s$ and $F(\cdot), G(\cdot) \in \sC\p{I; \sL\p{H,K}\cap\sJ_\infty^s}$. The operators $\Sigma_0, F(\cdot),$ and $G(\cdot)$ are fixed parameters and $X(\cdot)$ is the input variable in the definition of $\sQ\p{\cdot}$. It is then clear that \eqref{eqn: time-dependent problem} can be written succinctly as
\begin{equation} \label{eqn: time-dependent short form}
    \Sigma(t) = \Phi_t\br{\sQ\p{\Sigma(\cdot)}}.
\end{equation}
We will verify the range of $\sQ(\cdot)$ in the following.

\begin{proposition}
Let $H$ be a separable Hilbert space, and $K$ be a compact subspace of $H$ and $I:=[0,\tau]$ be a compact time interval, where $\tau \in \mathbb R_+$. Then, if $\Sigma_0 \in \sL\p{H,K}\cap \sJ_\infty^s$, and $F(\cdot),G(\cdot) \in \sC\p{I; \sL\p{H,K}\cap \sJ_\infty^s}$, we have that $\sQ(\cdot)$, as defined in \eqref{eqn: Q def}, is a mapping from $\sC\p{I; \sJ_\infty^s}$ into $\sL\p{H,K}\cap\sJ_\infty^s \times \sC\p{I; \sL\p{H,K} \cap \sJ_\infty^s}$. 
\end{proposition}
\begin{proof}
We only need to show that $(XGX)(\cdot)$ belongs to $\sC\p{I; \sL(H,K)\cap\sJ_\infty^s}$ since $\Sigma_0 \in \sL\p{H,K}\cap \sJ_\infty^s$ and $F(\cdot), G(\cdot) \in \p{I;\sL(H,K)\cap\sJ_\infty^s}$, by definition. The result of this proposition follows immediately from a direct application of Corollary \ref{corollary: c1}, and seeing that $(XGX)^*(t) = (X^*G^*X^*)(t) = (XGX)(t)$ on $H$ for any $X \in \sC(I; \sL(H,K) \cap \sJ_\infty^s)$.
\end{proof}

Similarly, we will define $\Phi_{n,t}: \sJ_\infty^s\times \sC\p{I; \sJ_\infty^s} \rightarrow \sC\p{I; \sJ_\infty^s}$, as the following approximation to the Lyapunov operator
\begin{equation} \label{eqn: lyapunov operator approximation definition}
    \Phi_{n,t}\stackvec{X}{Y(\cdot)} := S_n(t) X S^*_n(t) + \int_0^t S_n(t-s) Y(s) S^*_n(t-s) ds,
\end{equation}
for all $X \in \sJ_\infty^s$ and $Y \in \sC\p{I;\sJ_\infty^s}$. Likewise, \eqref{eqn: approximate time-dependent problem} can be written succinctly as
\begin{equation} \label{eqn: time-dependent short form approximation}
    \Sigma_n(t) = \Phi_{n,t}\br{\sQ\p{\Sigma_n(\cdot)}}
\end{equation}
after recalling, from the definition of $\sQ(\cdot)$, that $\sQ\p{\Sigma_n(\cdot)} = \left[ \begin{array}{c} \Sigma_0 \\ F(\cdot) - \p{\Sigma_n G \Sigma_n}\p{\cdot} \end{array} \right]$.

By now, it should be clear that the convergence of $\Sigma_n(\cdot)$ to $\Sigma(\cdot)$ will depend on the properties of the mappings $\Phi_t(\cdot), \Phi_{n,t}(\cdot)$, and $\sQ(\cdot)$. With this, we proceed to prove that $\Phi_{n,t}(\cdot) \rightarrow \Phi_t(\cdot)$ as $n \rightarrow \infty$. 

\subsection{Convergence of the Time-Dependent Lyapunov Operator Approximation}
We begin this section by establishing the linearity and continuity of $\Phi_t(\cdot)$ and $\Phi_{n,t}(\cdot)$ in the following. 
\begin{lemma} \label{lemma: l1}
    Let $H$ be a separable Hilbert space, and $S(t), S_n(t) \in \sL\p{H}$ be $C_0$-semigroups defined on a compact time interval $t \in I$, where $I:=[0,\tau]$ and $\tau \in \mathbb R_+$. Then for all $X\in \sJ_\infty^s$ and $Y\in\sC\p{I; \sJ_\infty^s}$, we have that $\Phi_t(\cdot)$ and $\Phi_{n,t}(\cdot)$ are bounded linear mappings from $\sJ_\infty^s \times \sC\p{I; \sJ_\infty^s}$ to $ \sC\p{I; \sJ_\infty^s}$. Furthermore, $\Phi_t(\cdot)$ and $\Phi_{n,t}(\cdot)$ are bounded linear mappings from $\sL\p{H,K}\cap \sJ_\infty^s \times \sC\p{I; \sL\p{H,K} \cap \sJ_\infty^s}$ to $\sC\p{I; \sJ_\infty^s}$.
\end{lemma}
\begin{proof}
For brevity, we only prove the result for $\Phi_t(\cdot)$. It will become clear that the same sequence of steps can be repeated with instead $\Phi_{n,t}(\cdot)$ to establish its linearity and continuity.

From the definition of $\Phi_t(\cdot)$ (see \eqref{eqn: lyapunov operator definition}), it is easy to see that 
\begin{equation*}
\begin{aligned}
    \Phi_t\stackvec{X + \alpha C}{Y(\cdot) + \alpha D(\cdot)} &= S(t)\p{X+\alpha C} S^*(t) + \int_0^t S(t-s)\p{Y + \alpha D}(s) S^*(t-s)ds  \\
        &= S(t)X S^*(t) + \int_0^t S(t-s)Y(s)S^*(t-s)ds \\
        &\quad + \alpha \p{S(t) C S^*(t) + \int_0^t S(t-s) D(s) S^*(t-s)ds} \\
        &= \Phi_t\stackvec{X}{Y(\cdot)} + \alpha \Phi_t\stackvec{C}{D(\cdot)}.
\end{aligned}
\end{equation*}
for all $X,C \in \sJ_\infty^s$, $Y(\cdot),D(\cdot) \in \sC\p{I; \sJ_\infty^s}$, and $\alpha \in \mathbb R$. Therefore, it is established that $\Phi_t: \sJ_\infty^s\times \sC\p{I;\sJ_\infty^s} \rightarrow \sC\p{I; \sJ_\infty^s}$ is a linear mapping, where the self-adjoint property of $\Phi_t\p{\cdot}$ is established by seeing that $S(t)XS^*(t)$ is self-adjoint on $H$ for all $X \in \sJ_\infty^s$. The continuity of $\Phi_t(\cdot)$ with respect to $t\in I$ then follows from the continuity of the input variables and fixed parameters in the definition of $\Phi_t(\cdot)$ with respect to $t \in I$.

The boundedness properties of the Lyapunov operator is then established by seeing that
\begin{equation*}
\begin{aligned}
\norm{\Phi_t\stackvec{X}{Y(\cdot)}}_{\sC\p{I; \sJ_\infty}} &\leq \sup_{t\in I} \Br{ \norm{S(t) X S^*(t)}_{\sJ_\infty} + \int_0^t \norm{S(t-s)Y(s)S^*(t-s)}_{\sJ_\infty}ds} \\
            &\leq M^2e^{2\omega \tau} \p{\norm{X}_{\sJ_\infty} 
                + \int_0^\tau e^{-2\omega s} \norm{Y(s)}_{\sJ_\infty} ds} \\
            &\leq M^2e^{2\omega \tau} \p{   
                \norm{X}_{\sJ_\infty} 
                + \frac{1}{2\omega}
                    \p{1 - e^{-2\omega\tau}}
                    \norm{Y}_{\sC\p{I; \sJ_\infty}}
                }\\
            &\leq C_{M,\omega,\tau} \p{\norm{X}_{\sJ_\infty} + \norm{Y}_{\sC\p{I;\sJ_\infty}}},
\end{aligned}
\end{equation*}
after applying \eqref{eqn: semigroup bound}. Hence, $\Phi_t: \sJ_\infty^s \times \sC\p{I; \sJ_\infty^s} \rightarrow \sC\p{I; \sJ_\infty^s}$ is a bounded operator on $\sJ_\infty^s \times \sC\p{I; \sJ_\infty^s}$. A continuation of the above series of inequalities yields
\begin{equation*}
\begin{aligned}
&\norm{\Phi_t\stackvec{X}{Y(\cdot)}}_{\sC\p{I; \sJ_\infty(H)}} \leq C_{M,\omega,\tau}  \p{\norm{X}_{\sL(H,K)} + \norm{Y(\cdot)}_{\sC\p{I;\sL(H,K)}}}
\end{aligned}
\end{equation*}
where we have used the fact that $\sL\p{K,H} \subset \sL\p{H}$. With this, we have also demonstrated that $\Phi_t$ is a bounded linear operator from $\sL(H,K)\cap \sJ_\infty^s \times \sC\p{I;\sL(H,K)\cap \sJ_\infty^s}$ into $\sC\p{I; \sJ_\infty^s}$.

Applying the same strategy, with instead the estimate \eqref{eqn: semigroup approximation bound}, establishes the linearity and boundedness of $\Phi_{n,t}$ as a mapping from $\sJ_\infty^s \times \sC\p{I;\sJ_\infty^s}$ and $\sL(H,K)\cap \sJ_\infty^s \times \sC\p{I; \sL(H,K)\cap\sJ_\infty^s}$ into $\sC\p{I; \sJ_\infty^s}$. This concludes this proof. 
\end{proof}

Now that we have established that $\Phi_t(\cdot)$ and $\Phi_{n,t}(\cdot)$ are continuous linear operators, we can then begin to analyze the convergence of $\Phi_{n,t}(\cdot)$ to $\Phi_t(\cdot)$ as $n\rightarrow \infty$ in the $\sJ_\infty$ norm. The analysis is presented in the proof of the following. 

\begin{lemma} \label{lemma: pointwise and uniform convergence}
    Let $H$ be a separable Hilbert space and $K$ be a compact subspace of $H$. Additionally, let $S(t) \in \sL(H)$ be a $C_0$-semigroup and $\left\{S_n(t) \right\}_{n=1}^\infty \in \sL(H)$ be an approximating sequence of $C_0$-semigroups on $H_n$ so that
$$
\lim_{n\rightarrow\infty}\norm{S(t)\phi - S_n(t)\phi}_H = 0 \quad \forall \phi \in H,
$$
for all values in a bounded time interval $t \in I$, where $I:= [0,\tau]$ with $\tau \in \mathbb R_+$. Furthermore, let $\Phi_t(\cdot), \Phi_{n,t}(\cdot) \in \sJ_\infty^s \times \sC\p{I; \sJ_\infty^s} \rightarrow \sC\p{I; \sJ_\infty^s}$ be defined as in \eqref{eqn: lyapunov operator definition} and \eqref{eqn: lyapunov operator approximation definition} respectively.
Then, we have that
\begin{equation} \label{eqn: a3}
\lim_{n\rightarrow\infty} \norm{\Phi_t\stackvec{X}{Y(\cdot)} - \Phi_{n,t}\stackvec{X}{Y(\cdot)}}_{\sC\p{I;\sJ_\infty}} = 0 
\end{equation}
for all non-negative self-adjoint operators $X \in \sJ_\infty^s$ and $Y \in \sC\p{I; \sJ_\infty^s}$, and that there exists a constant $C_{M, \omega, \tau} \in \mathbb R_+$ such that 
\begin{equation} \label{eqn: a4}
\begin{aligned}
&\norm{\Phi_t\stackvec{X}{Y(\cdot)} - \Phi_{n,t}\stackvec{X}{Y(\cdot)}}_{\sC\p{I; \sJ_\infty(H)}}\\
&\qquad\qquad\qquad\qquad\leq C_{M,\omega,\tau} \sup_{t \in I} \norm{S(t) - S_n(t)}_{\sL\p{K,H}} \p{\norm{X}_{\sL(H,K)} + \norm{Y(\cdot)}_{\sC\p{I; \sL(H,K)}} },
\end{aligned}
\end{equation}
for all $X \in \sL(H,K)\cap \sJ_\infty^s$ and $Y(\cdot) \in \sC\p{I; \sL(H,K)\cap\sJ_\infty^s}$. 
\end{lemma} 
\begin{proof} 
Throughout this proof, we will denote $M_* := \sup\Br{M_n}_{n=1}^\infty$ and $\omega_* := \sup\Br{\omega_n}_{n=1}^\infty$. From the definition of $\Phi_t(\cdot)$ and $\Phi_{n.t}(\cdot)$, we see that
$$
\begin{aligned}
&\norm{ \Phi_t \stackvec{X}{Y(\cdot)} - \Phi_{n,t}\stackvec{X}{Y(\cdot)}}_{\sC\p{I; \sJ_\infty}} \\
    &\quad \leq \sup_{t\in I} 
    \p{\norm{S(t) X \p{S(t) - S_n(t)}^*}_{\sJ_\infty} + \norm{\p{S(t)-S_n(t)}XS_n^*(t)}_{\sJ_\infty}}\\
    &\qquad + \int_0^\tau \norm{S(t-s)Y(s) \p{S(t-s) - S_n(t-s)}^*}_{\sJ_\infty} ds \\
    &\qquad + \int_0^\tau \norm{\p{S(t-s) - S_n(t-s)}Y(s)S_n^*(t-s)}_{\sJ_\infty}ds \\
    &\leq \sup_{t\in I} \p{\norm{\p{S(t) - S_n(t)}X}_{\sJ_\infty}\p{\norm{S(t)}_{\sL(H)} + \norm{S_n(t)}_{\sL(H)}}} \\
    &\qquad + \int_0^\tau \norm{\p{S(t-s) - S_n(t-s)}Y(s)}_{\sJ_\infty} \p{\norm{S(t-s)}_{\sL(H)} + \norm{S_n(t-s)}_{\sL(H)}} ds.
\end{aligned}
$$
The assertion \eqref{eqn: a3} then follows by an application of Lemma \ref{lemma: t2} and seeing that the integral is defined on a bounded interval. 

It then follows by continuing the previous sequence of inequalities that
$$
\begin{aligned}
&\norm{ \Phi_t \stackvec{X}{Y(\cdot)} - \Phi_{n,t}\stackvec{X}{Y(\cdot)}}_{\sC\p{I; \sJ_\infty}} \\
    &\quad \leq 2M_* e^{\omega_* \tau} \p{\sup_{t\in I}\norm{S(t) - S_n(t)}_{\sC\p{I; \sL(K,H)}} \norm{X}_{\sL\p{H,K}} + \int_0^\tau e^{-\omega_*s}\norm{\p{S(t-s) - S_n(t-s)}}_{\sL\p{K,H}} \norm{Y(s)}_{\sL\p{H,K}} ds} \\
    &\quad \leq 2M_*e^{\omega_* \tau} \sup_{t\in I}\norm{S(t) - S_n(t)}_{\sC\p{I; \sL(K,H)}}\p{ \norm{X}_{\sL(H,K)} + \frac{1}{\omega_*} \p{1 - e^{-\omega_*\tau}} \norm{Y(\cdot)}_{\sC\p{I; \sL(H,K)}}},
\end{aligned}
$$
thus proving \eqref{eqn: a4}.
\end{proof}

The well-posedness and the convergence of approximate Lyapunov operators, in pointwise and uniform sense, to the exact Lyapunov operator was established in the results presented in this subsection. In the final part of this section, we will apply Theorem \ref{theorem: BRR} to determine that the error of the approximate solution of the operator-valued Riccati equation is bounded by the semigroup approximation error.  

\subsection{Approximation Error}
Throughout this subsection, we will denote our function spaces of interest in the following manner
$$
    \sX:=\sC\p{I; \sJ_\infty^s} \quad \sY:=\sJ_\infty^s\times\sC\p{I;\sJ_\infty^s} \quad \sZ := \sL(H,K)\cap \sJ_\infty^s \times \sC\p{I; \sL(H,K)\cap \sJ_\infty^s}
$$
to make the application of Theorem \ref{theorem: BRR} easier to visualize. We note that $\sZ \subset \sY$ is a continuous embedding by recalling the definition of these spaces. 

We begin by seeing again that \eqref{eqn: time-dependent problem} and \eqref{eqn: approximate time-dependent problem} can be written in a fixed-point form as 
$$
    \Sigma(\cdot) = \Phi_t\br{\sQ(\Sigma(\cdot))} \quad \textrm{ and } \quad
    \Sigma_n(\cdot) = \Phi_{n,t}\br{\sQ(\Sigma_n(\cdot))}
$$
respectively. The linearity and boundedness of $\Phi_t(\cdot), \Phi_{n,t}(\cdot) \in \sL\p{\sY, \sX}$ was established in Lemma \ref{lemma: l1}. It is then easy to see from the definition of $\sQ(\cdot): \sX \rightarrow \sY$ that it is infinitely Fr\'echet with respect to $X \in \sX$. The first Gat\'eux derivative of $\sQ(\cdot)$ with respect to $X$ is given explicitly as
\begin{equation} \label{eqn: Derivative 1}
    \p{D_{X} \sQ(X)}T = \stackvec{0}{(TGX)(\cdot) + (XGT)(\cdot)},
\end{equation}
for any test operator $T \in \sX$. Since $G\in \sC\p{I; \sL\p{H,K}\cap\sJ_\infty^s}$ and that $0$ is an element in every space, we can see that $D_{X}\sQ \in \sL\p{\sX, \sZ}$ as a result of Lemma \ref{lemma: t1}.     Finally, we restate the pointwise and uniform convergence result that we have derived in Lemma \ref{lemma: pointwise and uniform convergence} in our new notation below, i.e.
$$
    \lim_{n\rightarrow \infty} \norm{\Phi_t(T) - \Phi_{n.t}(T)}_{\sX} = 0 \quad \forall T \in \sY
\quad\textrm{ and }\quad
    \lim_{n\rightarrow \infty} \norm{\Phi_t - \Phi_{n,t}}_{\sL\p{\sZ,\sX}} = 0.
$$
A direct comparison of the above observations against the conditions stated in \S\ref{subsection: BRR} allow us to apply the Brezzi-Rappaz-Raviart theorem (Theorem \ref{theorem: BRR}) to see that there exists a constant $C\in \mathbb R_+$ such that
\begin{equation*}
    \norm{\Sigma(\cdot) - \Sigma_n(\cdot)}_{\sC(I; \sJ_\infty)} \leq C \norm{\Phi_t\br{\sQ(\Sigma(\cdot))} - \Phi_{n,t}\br{\sQ(\Sigma(\cdot))}}_{\sC(I; \sJ_\infty)},
\end{equation*}
after recalling that $\sX := \sC(I; \sJ_\infty^s)$. Bounding this inequality with \eqref{eqn: a4} yields the first main result of this work.

\begin{theorem} \label{theorem: m1}
 Let $H$ be a separable Hilbert space and $K$ be a compact subset of $H$. Additionally, let $S(t) \in \sL(H)$ be a $C_0$-semigroup and $\left\{S_n(t) \right\}_{n=1}^\infty \in \sL(H)$ be an approximating sequence of $C_0$-semigroups on $H_n$ so that
$$
\lim_{n\rightarrow\infty}\norm{S(t)\phi - S_n(t)\phi}_H = 0 \quad \forall \phi \in H,
$$
for all values in a bounded time interval $t \in I$, where $I:= [0,\tau]$ with $\tau \in \mathbb R_+$. Furthermore, let $\Phi_t(\cdot), \Phi_{n,t}(\cdot) \in \sJ_\infty \times \sC\p{I; \sJ_\infty} \rightarrow \sC\p{I; \sJ_\infty}$ be defined as in \eqref{eqn: lyapunov operator definition} and \eqref{eqn: lyapunov operator approximation definition} respectively. And finally, let $\sQ(\cdot)$ be defined as in \eqref{eqn: Q def}. Then, we have that there exists a constant $\beta_{M,\omega,\tau} \in \mathbb R_+$ so that 
$$
    \norm{\Sigma(\cdot) - \Sigma_n(\cdot)}_{\sC\p{I; \sJ_\infty}} \leq \beta_{M, \omega, \tau} \sup_{t \in I} \norm{S(t) - S_n(t)}_{\sL(K,H)} \norm{\sQ(\Sigma(\cdot))}_{\sL(H,K)\times\sC\p{I; \sL(H,K)}},
$$
where $\Sigma(\cdot), \Sigma_n(\cdot) \in \sC\p{I; \sJ_\infty^s}$ are the solutions of the operator-valued Riccati equation \eqref{eqn: time-dependent problem} and its approximation \eqref{eqn: approximate time-dependent problem} respectively. 
\end{theorem}

\begin{remark}
In practice, the above result needs to be supplemented with additional analysis to resolve the time approximation error associated with time-stepping techniques. This error can be accounted for by simply utilizing the Minkowski inequality and applying known approximation results for time-stepping, i.e.,
\begin{equation}\label{eqn: time treatment}
    \norm{\Sigma(\cdot) - \Sigma_{n,k}(\cdot)}_{\sC(I; \sJ_\infty)} \leq \norm{\Sigma(\cdot) - \Sigma_n(\cdot)}_{\sC(I; \sJ_\infty)} + \norm{\Sigma_n(\cdot) - \Sigma_{n,k}(\cdot)}_{\sC(I; \sJ_\infty)},
\end{equation}
where we have introduced $\Sigma_{n,k}(\cdot) \in \sC(I; \sJ_\infty^s)$ as the solution approximation where both space and time-approximation has been applied. The error between $\Sigma_n(\cdot)$ and $\Sigma_{n,k}(\cdot)$ can then be bounded with existing results, under sufficient smoothness of $\Sigma$ with respect to $t \in I$ (see \cite[Chapter 11]{quarteroni2010numerical}).  
\end{remark}

With this, we have established an abstract error estimate for the approximate solution to the time-dependent operator-valued Riccati equation. This estimate shows that the error bound is a function of the approximation error of the semigroup approximation used in the approximate equation. The implication of this result is that many previous results in the approximation of time-dependent processes can be applied in deriving error bounds for the approximate solution of the operator-valued Riccati equation in a control system. We will qualify this claim in \S \ref{Section: Implications}. In the next section, we demonstrate in Theorem \ref{theorem: m2} that this observation can be applied to the time-independent variant of operator-valued Riccati equation.

\section{An Error Estimate for Time-Independent Approximations}
\label{section: time-independent approximation}
In this section, we will continue our analysis by deriving an abstract error estimate for approximate solutions to the time-independent operator Riccati equation. We begin by recalling again that the steady-state form of the operator Riccati equation is given as
\begin{equation} \label{eqn: steady-state riccati equation}
    \Sigma_\infty = \int_0^\infty S(s) \p{F - \Sigma_\infty G \Sigma_\infty} S^*(s) ds,
\end{equation}
where $S(t)\in \sL\p{H}$ is an exponentially stable semigroup satisfying $\norm{S(t)}_{\sL(H)} \leq Me^{-\alpha t}$, $F, G \in \sJ_\infty^s$ are compact, non-negative, and self-adjoint operators on $H$. We will again make the range of $F$ and $G$ explicit in this analysis by indicating that $F,G \in \sL\p{H,K} \cap \sJ_\infty^s$, where $K\subset H$ is a compact subset. For reference, it was demonstrated in Lemma \ref{lemma: steady-state well-posedness} that \eqref{eqn: steady-state riccati equation} is a well-posed problem. The corresponding approximation to \eqref{eqn: steady-state riccati equation} is then given by
\begin{equation} \label{eqn: steady-state riccati equation approximation}
    \Sigma_{\infty, n} = \int_0^\infty S_n(s) \p{F - \Sigma_{\infty,n} G \Sigma_{\infty,n}} S^*_n(s) ds,
\end{equation}
where $S_n(t) \in \sL\p{H}$ is an approximation to $S(t) \in \sL\p{H}$ satisfying $\norm{S_n(t)}_{\sL(H)} \leq M_n e^{-\alpha_n t}$. The analysis presented in this section will follow the same structure as the analysis presented in Section \ref{section: time-dependent approximation}. 

It is our goal to attempt to present the analysis in this section succinctly since its underlying logic is similar to (and simpler than) what we have presented in the previous section. To that end, we will first define the steady-state Lyapunov operator and its approximation and show that they are linear and bounded operators on $\sJ_\infty^s$. Then we move on to demonstrate that the approximate Lyapunov operator converges to the exact Lyapunov operator in the pointwise and uniform sense on compact subspaces of $\sJ_\infty^s$. And finally, we apply the Brezzi-Rappaz-Raviart theorem (Theorem \ref{theorem: BRR}) to derive the main result of this section. 

\subsection{Time-Independent Lyapunov Operator}
Let us define the \emph{time-independent Lyapunov operator} $\Phi_\infty(\cdot): \sJ_\infty^s \rightarrow \sJ_\infty^s$ as 
\begin{equation} \label{eqn: steady-state lyapunov operator}
    \Phi_\infty (X) := \int_0^\infty S(s)XS^*(s) ds
\end{equation}
for all $X \in \sJ_\infty^s$. Additionally, let us define $\sQ_\infty\p{\cdot}: \sJ_\infty^s \rightarrow \sL\p{H,K}\cap \sJ_\infty^s$ as 
\begin{equation}
    \sQ_\infty(X) := F - XGX.
\end{equation}
where $F,G \in \sL(H,K) \cap \sJ_\infty^s$. It is then clear that we can represent \eqref{eqn: steady-state riccati equation} by
\begin{equation} \label{eqn: short form 1}
    \Sigma_\infty = \Phi_\infty \br{\sQ_\infty(\Sigma_\infty)}.
\end{equation}
With this, we verify the range of $\sQ_\infty(\cdot): \sJ_\infty^s \rightarrow \sL(H,K) \cap \sJ_\infty^s$ in the following. 
\begin{proposition}
Let $H$ be a separable Hilbert space and $K$ be a compact subspace of $H$. Then if $F,G \in \sL(H,K)\cap \sJ_\infty^s$, we have that $\sQ_\infty(\cdot)$ is a mapping from $\sJ_\infty^s$ to $\sL(H,K) \cap \sJ_\infty^s$.
\end{proposition}
\begin{proof}
We need to only show that $XGX \in \sL(H,K)\cap \sJ_\infty^s$ for any $X \in \sJ_\infty^s$. This is easily seen by applying Lemma \ref{lemma: t1} and seeing that $\p{XGX}^* = X^*G^*X^* = XGX$ by definition. 
\end{proof}

The approximation of interest to the steady-state Lyapunov operator given in \eqref{eqn: steady-state lyapunov operator} is defined as 
\begin{equation} \label{eqn: steady-state lyapunov operator approximation}
    \Phi_{\infty,n} (X) := \int_0^\infty S_n(s) X S_n^*(s) ds.
\end{equation}
Similarly, the approximate steady-state operator Riccati equation \eqref{eqn: steady-state riccati equation approximation} can be written as
\begin{equation} \label{eqn: short form 2}
    \Sigma_{\infty,n} = \Phi_{\infty, n} \br{\sQ_\infty(\Sigma_{\infty,n})}.
\end{equation}
From inspecting \eqref{eqn: short form 1} and \eqref{eqn: short form 2}, that Theorem \ref{theorem: BRR} can be used to determine error estimates for the approximate solution to the time-independent operator-valued Riccati equation. 

We now verify the linearity and boundedness of $\Phi_\infty(\cdot)$ and $\Phi_{\infty, n}(\cdot)$ in the following.
\begin{lemma} \label{lemma: linearity and boundedness of steady-state lyapunov operator}
    Let $H$ be a separable Hilbert space and let $\Phi_\infty(\cdot)$ and $\Phi_{\infty,n}(\cdot)$ be defined as in \eqref{eqn: steady-state lyapunov operator} and \eqref{eqn: steady-state lyapunov operator approximation} respectively. Then $\Phi_\infty(\cdot), \Phi_{\infty, n}(\cdot): \sJ_\infty^s \rightarrow \sJ_\infty^s$ are bounded linear operators on $\sJ_\infty^s$.
\end{lemma}
\begin{proof}
    In this proof, we will only prove the linearity and boundedness of $\Phi_\infty(\cdot):\sJ_\infty^s \rightarrow \sJ_\infty^s$ for brevity. The linearity and boundedness of $\Phi_{\infty, n}(\cdot): \sJ_\infty^s \rightarrow \sJ_\infty^s$ then follows from the same argument.
    
    We begin this proof by observing that $\p{S(t)XS^*(t)}^* = S(t)X^*S^*(t) = S(t)XS^*(t)$ on $H$ for any $X \in \sJ_\infty^s$. It is then clear from the definition of $\Phi_\infty(\cdot): \sJ_\infty^s \rightarrow \sJ_\infty^s$ that
$$
\begin{aligned}
    \Phi_\infty (X + \alpha Y) &= \int_0^\infty S(s)(X + \alpha Y) S^*(s) ds \\
                               &= \int_0^\infty S(s)XS^*(s)ds + \alpha\int_0^\infty S(s)YS^*(s)ds \\
                               &= \Phi_\infty(X) + \alpha \Phi_\infty(Y)
\end{aligned}
$$
for all $X,Y \in \sJ_\infty^s$ and $\alpha \in \mathbb R$, and hence $\Phi(\cdot)$ is a linear operator. 

The boundedness of $\Phi(\cdot)$ on $\sJ_\infty^s$, is established by seeing that
$$
\begin{aligned}
    \norm{\Phi_\infty (X)}_{\sJ_\infty} &\leq \int_0^\infty \norm{S(s)}^2 \norm{X}_{\sJ_\infty} ds \\
                         &\leq M^2\norm{X}_{\sJ_\infty} \int_0^\infty e^{-2\alpha s} ds \\
                         &= \frac{M^2}{2\alpha} \norm{X}_{\sJ_\infty}
\end{aligned}
$$
for all $X \in \sJ_\infty^s$. This concludes the proof. 
\end{proof}

\subsection{Convergence of Time-Independent Lyapunov Operator Approximations}
We now move forward to determine the convergence of $\Phi_{\infty, n}\p{\cdot} \rightarrow \Phi_{\infty}\p{\cdot}$ in the following.

\begin{lemma} \label{lemma: pointwise and uniform convergence of steady-state Lyapunov operator}
Let $H$ be a separable Hilbert space and $K$ be a compact subspace of $H$. Additionally, let $S(t) \in \sL\p{H}$ be an exponentially stable $C_0$-semigroup and $\Br{S_n(t)}_{n=1}^\infty \in \sL(H)$ be an approximating sequence of $C_0$-semigroups so that 
\begin{equation*}
    \lim_{n \rightarrow \infty} \norm{S(t)\phi - S_n(t)\phi }_{\sL(H)} = 0 \quad \forall \phi \in H
\end{equation*}
for any $t \in \mathbb R_+$. Under these assumptions, we observe that
\begin{equation} \label{eqn: az1}
    \lim_{n\rightarrow \infty} \norm{\Phi_\infty (X) - \Phi_{\infty, n} (X)}_{\sJ_\infty} = 0 \quad \forall X \in \sJ_\infty^s
\end{equation} 
and
\begin{equation} \label{eqn: w2}
    \norm{\Phi_\infty (X) - \Phi_{\infty, n} (X)}_{\sJ_\infty} \leq 2M\norm{X}_{\sL(H,K)} \int_0^\infty e^{-\alpha s}\norm{S(s) - S_n(s)}_{\sL\p{H,K}} ds \quad \forall X \in \sL(H,K), 
\end{equation}
where $M,\alpha \in \mathbb R_+$, and the operators $\Phi_\infty(\cdot), \Phi_{\infty,n}(\cdot): \sJ_\infty^s \rightarrow \sJ_\infty^s$ are defined as in \eqref{eqn: steady-state lyapunov operator} and \eqref{eqn: steady-state lyapunov operator approximation} respectively.
\end{lemma}
\begin{proof}
We begin the proof by seeing that $\Br{\Phi_{\infty, n}}_{n=1}^\infty: \sJ_\infty^s \rightarrow \sJ_\infty^s$ is a bounded sequence as a consequence of Lemma \ref{lemma: linearity and boundedness of steady-state lyapunov operator} since $S_n(t) \in \sL(H)$ is assumed to be exponentially stable for any $n \in [1, \infty)$. This observation validates the notion of convergence of $\Phi_{\infty, n} \rightarrow \Phi_\infty$. 

From the definition of $\Phi_\infty(\cdot)$ and $\Phi_{\infty, n}(\cdot)$, we see that
$$
\Phi_\infty (X) - \Phi_{\infty, n} (X) = \int_0^\infty \p{S(s) - S_n(s)}X S^*(s) ds + \int_0^\infty S_n(s) X \p{S(s) - S_n(s)}^* ds.
$$
Taking the $\sJ_\infty$-norm on both sides of the above equation then gives us the following inequality
$$
    \norm{\Phi_\infty (X) - \Phi_{\infty, n} (X)}_{\sJ_\infty} = \int_0^\infty \norm{S(s)X - S_n(s)X}_{\sJ_\infty} \p{\norm{S(s)}_{\sL(H)} + \norm{S_n(s)}_{\sL(H)}} ds.
$$
We establish \eqref{eqn: az1} by seeing that $\Phi_{\infty, n}(\cdot)$ and $\Phi_{n}(\cdot)$ are bounded and that $\norm{S(s)X - S_n(s)X}_{\sJ_\infty}$ is pointwise convergent by Lemma \ref{lemma: t2}. 

We continue bounding the inequality above by seeing that 
$$
   \norm{\Phi_\infty (X) - \Phi_{\infty, n} (X)}_{\sJ_\infty} \leq 2M_* \norm{X}_{\sL(H,K)} \int_0^\infty \norm{S(s) - S_n(s)}_{\sL(K,H)} e^{-\alpha s} ds,
$$
by observing that
$$
    \norm{S(t)X - S_n(t) X}_{\sJ_\infty} \leq \norm{S(t) - S_n(t)}_{\sL(K,H)} \norm{X}_{\sL(H,K)}
$$
for any $t \in \mathbb R_+$, and that $\norm{S(t)}_{\sL(H)} \leq Me^{-\alpha t}$ and $\norm{S_n(t)}_{\sL(H)} \leq M_n e^{-\alpha_n t}$. Thus \eqref{eqn: w2} is verified. 
\end{proof}

\subsection{Approximation Error}
In this section, we will utilize the Brezzi-Rappaz-Raviart Theorem (Theorem \ref{theorem: BRR}) to derive an error estimate for the approximation of the operator-valued time-independent Riccati equation on Hilbert spaces. To that end, we begin by defining 
\begin{equation*}
\sX := \sJ_\infty^s \quad  \sY := \sJ_\infty^s \quad  \sZ := \sL(H,K) \cap \sJ_\infty^s.
\end{equation*}
From the definition of $\sQ_\infty(\cdot): \sX \rightarrow \sZ$, we see that $\sQ_\infty(X)$ is infinitely Fr\'echet differentiable. Going further, we have that 
$$
    \left[ D_X\sQ_\infty(X) \right] T = TGX + XGT \quad \forall T \in \sX,
$$
therefore $D_X\sQ_\infty(X) \in \sL(\sX, \sZ)$ for any $X \in \sX$. In Lemma \ref{lemma: pointwise and uniform convergence of steady-state Lyapunov operator}, we have shown that 
$$
    \lim_{n\rightarrow \infty}\norm{\Phi_\infty (X) - \Phi_{\infty,n}(X)}_{\sX} = 0
$$ 
and 
$$
    \lim_{n\rightarrow \infty}\norm{\Phi_\infty(\cdot) - \Phi_{\infty, n}(\cdot)}_{\sL\p{\sZ, \sX}} = 0. 
$$
With these observations, we see that there exists a constant $C \in \mathbb R_+$ that satisfies
$$
    \norm{\Sigma_\infty - \Sigma_{\infty,n}}_{\sJ_\infty} \leq C \norm{\Phi_\infty\br{\sQ_\infty(\Sigma_\infty)} - \Phi_{\infty,n}\br{\sQ_\infty(\Sigma_\infty)}}_{\sJ_\infty},
$$
after applying Theorem \ref{theorem: BRR} and recalling that $\sX := \sJ_\infty^s$. Applying \eqref{eqn: w2} then allows us to establish the main result of this section.  
\begin{theorem} \label{theorem: m2}
    Let $H$ be a separable Hilbert space and $K$ be a compact subspace of $H$. Additionally, let $S(t) \in \sL(H)$ be an exponentially stable semigroup and $\Br{S_n(t)}_{n=1}^\infty \in \sL(H)$ be an approximating sequence of exponentially stable semigroups so that 
$$
    \lim_{n\rightarrow \infty} \norm{S(t)\phi - S_n(t)\phi}_{\sL(H)} = 0.
$$
Additionally, let $\Sigma_\infty, \Sigma_{\infty,n} \in \sJ_\infty^s$ be the solution of \eqref{eqn: steady-state riccati equation} and \eqref{eqn: steady-state riccati equation approximation} respectively. Then there exists a constant $\gamma_M,\alpha \in \mathbb R_+$ so that
$$
    \norm{\Sigma_\infty - \Sigma_{\infty, n}}_{\sJ_\infty(H)} \leq \gamma_M \norm{\sQ_\infty(\Sigma_\infty)}_{\sL(H,K)} \int_0^\infty e^{-\alpha s} \norm{S(s) - S_n(s)}_{\sL(K,H)} ds. 
$$
\end{theorem}

With this result, we conclude the analysis for the approximation of the time-independent operator Riccati equation. Like in the approximation theorem presented in the previous section (Theorem \ref{theorem: m1}), we observe that the error bound is a function of the error of the semigroup approximation. Similarly, we can utilize the plethora of existing results on the approximation of time-dependent processes to derive error estimates for time-independent operator Riccati equations. We will apply this knowledge in the following sections.

\section{A Scalar Model LQR Control Problem}
In this section, we will demonstrate that the general theory for an asymptotic approximation of time-independent operator-valued Riccati equation applies to scalar  systems. This example serves as a first verification of the results presented in this work, in that it illustrates that the error estimate computed using the general theory provided in \S \ref{section: time-independent approximation} converges at the same rate as the error estimate derived from direct calculation. To this end, we consider the following exact system
\begin{equation}\label{eqn: simple exact system}
\left\{
\begin{aligned}
\frac{dz(t)}{dt} &= -az(t) + bu(t) \quad \textrm{ on } t\in(0, \infty) \\
            z(0) &= z_0,
\end{aligned}
\right.
\end{equation}
where $a \in \mathbb R_+$, $b \in \mathbb R_{ext}$, $z(\cdot) \in \sC\p{\mathbb R_+; \mathbb R_{ext} }$ is the solution variable, $z_0 \in \mathbb R$ is the initial condition, and $\mathbb R_{ext}:=\mathbb R \cup \Br{-\infty} \cup \Br{+\infty}$ denotes the extended real number system. 

Let us now define $a_\epsilon := (a + \epsilon)$, where $\epsilon \in \mathbb R$ is a small perturbation, then the approximate system we wish to analyze is given by 
\begin{equation}\label{eqn: simple approximate system}
\left\{
\begin{aligned}
\frac{dz_\epsilon(t)}{dt} &= -a_\epsilon z_\epsilon (t) + bu_\epsilon(t) \quad \textrm{ on } t\in(0, \infty) \\
            z_\epsilon(0) &= z_0,
\end{aligned}
\right.
\end{equation}
where $z_\epsilon(\cdot) \in \sC\p{\mathbb R_+; \mathbb R_{ext}}$ is the approximate solution with the initial condition coinciding with that of the exact problem. We will assume that $\epsilon := \epsilon(n)$ such that $\lim_{n\rightarrow \infty} \epsilon(n) = 0$ so that ``$\lim_{n\rightarrow\infty}$'' is equivalent to ``$\lim_{\epsilon \rightarrow 0}$''. In practice, the error between $a$ and $a_\epsilon$ is often induced by imprecise measurement of model parameters in real systems. 

The control inputs $u(\cdot), u_\epsilon(\cdot) \in L^2(\mathbb R)$ are chosen optimally as to minimize the following cost functionals
$$
    \mathcal J(u) := \int_0^\infty \br{ c^2 z^2(t;u) + u^2(t)} dt
    \quad
    \textrm{and}
    \quad
    \mathcal J_\epsilon (u_\epsilon) := \int_0^\infty \br{c^2 z_\epsilon^2(t;u) + u_\epsilon^2(t)}dt.
$$
for the exact and approximate systems respectively, where we have defined $c \in \mathbb R \cup \Br{-\infty} \cup \Br{+\infty}$ to be a weighting scalar. The optimal control for the exact and approximate systems are then given by $u_{opt}(t):= -bkz(t)$ and $u_{\epsilon, opt}(t) := -bk_\epsilon z_\epsilon(t)$, where the feedback gains $k, k_\epsilon \in \mathbb R$ are defined by $k:= b\sigma$ and $k_\epsilon := b\sigma_\epsilon$ respectively, and $\sigma, \sigma_\epsilon$ are the unique positive indefinite solutions of the following scalar Riccati equations
\begin{equation} \label{eqn: scalar Riccati equation}
    -2a\sigma -g\sigma^2 + f = 0 \quad \textrm{and} \quad
    -2a_\epsilon\sigma_\epsilon -g\sigma_\epsilon^2 + f = 0 ,
\end{equation}
where we have defined $g:= b^2$ and $f:= c^2$ for notational simplicity. 

We begin our analysis by deriving an error estimate through a direct computation. To that end, we have from the quadratic formula that the positive solutions to \eqref{eqn: scalar Riccati equation} are given by
$$
    \sigma = -\frac{a}{g} + \frac{\sqrt{a^2 + gf}}{g}
    \quad
    \textrm{ and }
    \quad 
    \sigma_\epsilon =  -\frac{a_\epsilon}{g} + \frac{\sqrt{a_\epsilon^2 + gf}}{g}.
$$
It then follows that
$$
\begin{aligned}
    |\sigma - \sigma_\epsilon| &\leq \frac{|a - a_\epsilon|}{g} + \frac{1}{g} \left|\sqrt{a^2 + gf} - \sqrt{a_\epsilon^2 + gf}\right| \\
    & \leq \frac{\epsilon}{g} + \frac{1}{g}\left| a + \sqrt{gf} - \sqrt{a_\epsilon^2 + gf} \right| \\
    &\leq \frac{\epsilon}{g} + \frac{|a - a_\epsilon|}{ g } \\
    &= \frac{2\epsilon}{g},
\end{aligned}
$$
after recalling the definition of $a_\epsilon$, and the inequality $\left|\sqrt{a} - \sqrt{b}\right| \leq \left|\sqrt{a-b}\right|$ for any $a,b \in \mathbb R_+$, where $|\cdot| : \mathbb C \rightarrow \mathbb R_+$ is interpreted here as the modulus operator for complex numbers. 

We now use Theorem \ref{theorem: m2} to derive an error estimate for the approximation error of the scalar control system. We begin by first relating this problem to the abstract problem given in Section \ref{section: time-independent approximation}. To that end, we first observe that we can set $H = \mathbb R_{ext}$ as the extended real line and $K = H$. Note that $K$ is a compact subset of $H$ because $H$ itself is compact \cite{tao2008compactness}. Then, letting $F = f$ and $G = g$, we have that $F$ and $G$ are both elements of $\sL\p{H,K}$. Denoting $A:=-a$ and $A_\epsilon := -a_\epsilon$, we can define $S(t) = e^{-at}$ and $S_\epsilon (t) = e^{-a_\epsilon t}$. It is then clear that $S(t), S_\epsilon (t) \in \sL\p{H}$ for all $t \in \mathbb R_+$. It then follows that the solution variable can be defined as $\Sigma_\infty = \sigma$. The absolute value $|\cdot|$ will then serve as the norm for $\sJ_\infty(H)$, $\sL\p{H,K}$, and $\sL\p{H}$. These observations then allow us to apply Theorem \ref{theorem: m2} to the scalar control system. We will now revert back to the scalar notation in the remainder of this section to be consistent with the analysis presented in the previous paragraph.

\begin{figure}
    \includegraphics[width = 0.65\textwidth]{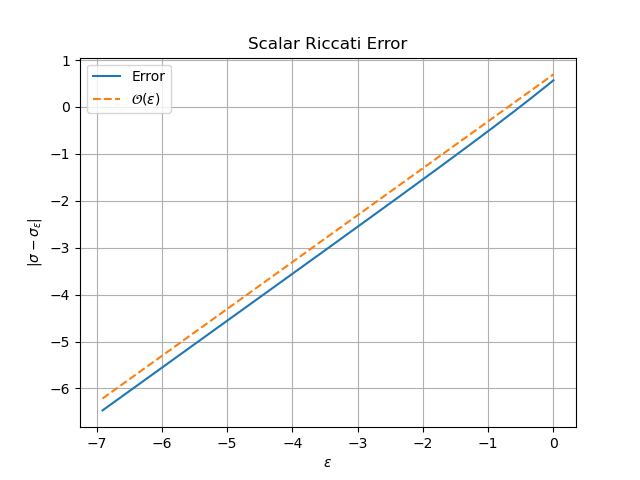}
    \caption{Numerical convergence plot with $\epsilon \in (0, 1]$ for the scalar-valued Riccati equation approximation \eqref{eqn: scalar Riccati equation} with $a,f,g = 1$.}
        \label{fig: scalar approximation convergence}
\end{figure}

We begin the analysis by determining an error estimate for the semigroups associated with the exact and approximate scalar systems. The uncontrolled state for the exact system is governed by
$$
\left\{
\begin{aligned}
    \frac{d\phi(t)}{dt} &= -a \phi(t) \quad \textrm{ on } t\in(0, \infty) \\ 
    \phi(0) &= \phi_0,
\end{aligned}
\right.
$$
and likewise, the approximate system is governed by
$$
\left\{
\begin{aligned}
    \frac{d\phi_\epsilon(t)}{dt} &= -a_\epsilon \phi_\epsilon(t) \quad \textrm{ on } t\in(0, \infty)\\ 
    \phi_\epsilon(0) &= \phi_0.
\end{aligned}
\right.
$$
It then follows that the difference between the solution of the exact and approximate uncontrolled system is given by the following
$$
\begin{aligned}
    \frac{d}{dt} \p{\phi(t) - \phi_\epsilon(t)} &= -a \phi(t) + a_\epsilon \phi_\epsilon(t) \\
    &= -a \p{\phi(t) - \phi_\epsilon(t)} + \epsilon \phi_\epsilon(t)
\end{aligned}
$$
subject to the intial condition $\phi(0) - \phi_\epsilon(0) = 0$. Solving for $\phi(t) - \phi_\epsilon(t)$ as the solution variable yields
$$
\begin{aligned}
    \phi(t) - \phi_\epsilon(t) &= \epsilon e^{-at} \phi_0 \int_0^t e^{-\epsilon s} ds \\
    &= e^{-a t} \phi_0 \p{1 - e^{-\epsilon t} } \\
    &= e^{-a t} \phi_0 \p{-\epsilon t + R^2(\xi;\epsilon)} \\
\end{aligned}
$$
and hence 
$$
    |\phi(t) - \phi_\epsilon(t)| \leq C\epsilon te^{-at}. 
$$
where $C \in \mathbb R_+$ is a positive constant and  $R^2(\xi; \epsilon)$ is the second order remainder term in the Taylor series \cite[Theorem 5.15]{rudin1964principles} evaluated at $\xi \in \mathbb R$ chosen so that $\epsilon t + R^2(\xi; \epsilon) = e^{\epsilon t} - 1$. It should be noted that $ \lim_{\epsilon \rightarrow 0} R^2(\xi; \epsilon) = 0$. Since $\phi(t) = e^{-at}$ and $\phi_\epsilon(t) = e^{-a_\epsilon t}$, we have that
$$
    |e^{-a t} - e^{-a_\epsilon t}| \leq C_1 \epsilon t e^{-a t}. 
$$
Using this result with Theorem \ref{theorem: m2} then implies that
$$
\begin{aligned}
    |\sigma - \sigma_\epsilon| 
    &\leq C\epsilon |f - g\sigma^2| \int_0^\infty s e^{-2as}ds \\
    &\leq \frac{C|f|}{2a^2} \epsilon.
\end{aligned}
$$
We have therefore demonstrated that Theorem \ref{theorem: m2} yields the same $\mathcal O(\epsilon)$ convergence rate predicted by direct computation, albeit with a different positive constant. The $\mathcal O(\epsilon)$ convergence rate is verified through a computational experiment as illustrated in Figure \ref{fig: scalar approximation convergence}.

\section{Infinite Dimensional Model LQR Control Problems} \label{Section: Implications}
In this section, we demonstrate the applicability of the theory derived in the earlier sections of this work to the finite element approximation (see, e.g., \cite{strang1974analysis}) of LQR control systems. Here, the finite element method is applied to an one-dimensional and two-dimensional thermal control system and an one-dimensional weakly damped wave control system, where the time-dependent operator-valued Riccati equation is approximated in the first model problem and its time-independent variant is approximated in the second and third model problems. Because the error in the operator norm is difficult to compute, we instead choose to compute the error of the functional gain associated with the optimal feedback law for the approximate control system. We demonstrate that the error of this approximate functional gain is dependent on the error of the approximate solution to the relevant Riccati equation. We conclude each subsection by presenting error convergence computed from the computational implementation of the model problems. This section itself is concluded with two examples that illustrate the numerical error computed on the model problems for when the theoretical conditions for optimality are violated. 

\subsection{Notation}
In the model problems considered in the following, the problem domain will be denoted $\Omega \subset \mathbb R^d$, where $d\in \mathbb N$ denotes the domain dimension, We will then denote $\Gamma := \partial \Omega$ as the boundary of $\Omega$ and $\mathbf x \in \Omega$ as the spatial coordinate $\textbf{x} := \p{x_1, \ldots, x_d}$, where each $x_j$ denotes the $j^{\textrm{th}}$ element of $\mathbf x \in \mathbb R^d$. Let now $\boldsymbol{\alpha} := \p{\alpha_1, \ldots, \alpha_d}$ be a multi-index where each element $\alpha_j$ is a positive integer. We will then denote, for any function $v$ defined on $\Omega$,
$$
    \partial_{\mathbf x}^{\boldsymbol{\alpha}} v := \frac{\partial^{\alpha_1}\cdots\partial^{\alpha_d} v}{\partial {x_1}^{\alpha_1} \cdots \partial {x_d}^{\alpha_d}}
$$
as the generalized multi-dimensional partial differentiation operator. We will then define $|\boldsymbol{\alpha}| := \sum_{j=1}^d \alpha_j$ as the total order of differentiation. 

Moving onto functional analytic notation, we define
$$
    \dual{v,w} = \int_{\Omega} vw d\mathbf x
$$
as the duality pairing between $v\in X$ and $w \in X'$, where $X$ is a topological vector space of vectors defined on $\Omega$, and $X'$ denotes the dual space associated with $X$. The norm on $X$ is then denoted $\norm{v}_X$. 

Let now $Y$ be a second topological vector space. The tensor product space $X\times Y$ can then be defined as 
$$
    X\times Y:= \Br{\mathbf v = \p{v,w} : v \in X \textrm{ and } w \in Y},
$$
where $X$ and $Y$ are topological vector spaces. The norm on $X\times Y$ can then defined by
$$
    \norm{\mathbf v}_{X\times Y} = \norm{v}_X + \norm{w}_Y,
$$
for all $\mathbf v:=\p{v,w} \in X\times Y$. 

The main topological vector spaces of interest for the discussion presented here are the Sobolev spaces $L^2(\Omega)$ and $H^{k}(\Omega)$. These spaces are defined as follows:
\begin{equation*}
    \begin{aligned}
        L^2(\Omega) := \Br{ v \textrm{ is Lebesgue measurable on } \Omega : \int_\Omega {v(x)}^2 d\xb < \infty } \\
        H^{k}(\Omega) := \Br{v \in L^2(\Omega) : \partial_\mathbf{x}^{\boldsymbol{\alpha}} v \in L^2(\Omega) \textrm{ for } |\boldsymbol{\alpha}|=1,\ldots,k}.
    \end{aligned}
\end{equation*}
The Sobolev norms are then defined as 
$$
    \norm{v}_{L^{2}(\Omega)} := \sqrt{\int_\Omega v^2 d\xb} \quad \textrm{and} \quad \norm{v}_{H^{k}(\Omega)} = \sqrt{\sum_{|\boldsymbol{\alpha}|=0}^k \int_{\Omega} \p{\partial_{\textbf x}^{\boldsymbol{\alpha}} v}^2 d\xb}.
$$
We then define
$$
    H^1_0(\Omega) := \Br{v \in H^1(\Omega) : v = 0 \textrm{ on } \Gamma}.
$$
This space consists of all functions belonging to $H^1(\Omega)$ that is vanishes on $\Gamma$. Further, we define the subspaces
$$
    N^k(\Omega) := \Br{v \in H^{k}(\Omega) : \partial_\xb^{\alphab} \partial_\nb v = 0 \textrm{ on } \Gamma, \quad |\alphab|=0,\ldots,k }
$$ 
and
$$
    D^k(\Omega) := \Br{v \in H^k(\Omega) : \partial_\xb^{\alphab} z = 0 \textrm{ on } \Gamma, \quad |\alphab|=0,\ldots, k},
$$
where $N^k(\Omega)$ and $D^k(\Omega)$ describe the space of functions that satisfy the homogeneous Neumann and Dirichlet boundary conditions in its derivatives up to order $k \in \mathbb N$.

\subsection{Finite-Horizon LQR Control of a 1D Thermal Process} \label{sec: ex2}
\subsubsection{Problem Definition}
In this example problem, we consider an one-dimensional thermal control system whose state is controlled by adjusting the heat flux into and out of the interior domain. The domain for this problem is $\Omega:= \p{-1,1}$ with the boundary $\Gamma:= \Br{-1} \cup \Br{1}$. Additionally, we will let $I:= \br{0,\tau}$ be the finite time domain on which we want to exercise the control law. The mathematical model for this system is then given by 
\begin{equation}\label{eqn: example problem 2}
    \left\{
    \begin{aligned}
        \partial_t z(x,t) &= \alpha \partial_x^2z(x,t) - \beta z(x,t) + b(x)u(t)  & \quad \textrm{ on } \Omega \times I \\
        \partial_x z(x,t) &= 0  &\quad \textrm{ on } \Gamma \\
        z(x,0) &= z_0(x)    &\quad \textrm{ on } \Omega,
    \end{aligned}
    \right.
\end{equation}
where $u(\cdot) \in L^2(I, \mathbb R)$ is the control input and $b(x) \in N^{k+3}(\Omega)$ is the actuator function. The $\alpha, \beta \in \mathbb R_+$ coefficients determines the rate of thermal dissipation within the domain and out of the domain respectively. Here, we apply the LQR methodology and choose an optimal control $u_{opt}(\cdot) \in L^2(I; \mathbb R)$ that minimizes the following cost functional
\begin{equation} \label{eqn: thermal cost functional}
    \mathcal J(u) := \int_{-1}^1 z( x, \tau;u) \cdot P_\tau z(x, \tau;u) d x + \int_0^\tau \p{\int_{-1}^1 q(x)z(x, t;u)dx}^2 + |u(t)|^2 dt,
\end{equation}
where $q \in N^{k+3}(\Omega)$ and $P_\tau \in \sL\p{H,K} \cap \sJ_\infty^s$ is the terminal state weighting operator. 

We will now map the model equations into the infinite-dimensional state space setting. To that end, let us now define $H:= L^2(\Omega)$ as the state space for the thermal control problem and $A: \mathcal D(A) \rightarrow H$ as the generator $\alpha \partial_{xx}(\cdot) - \beta(\cdot)$, where $\mathcal D(A) := H^2(\Omega) \cap N^0(\Omega)$. The theory behind the generation of the dissipative semigroup $S(t)\in \sL(H)$ from $A: \mathcal D(A) \rightarrow H$ can be found in \cite[Chapter 7.4]{evans2022partial}. For this problem, we will define $B: \mathbb R \rightarrow H$ and $C: H \rightarrow \mathbb R$ as the following 
$$
    B(\cdot) := b(x)(\cdot) \quad \textrm{ and } \quad  C(\cdot) := \int_{-1}^1 q(x) (\cdot) dx.
$$
We will assume that $b(x), q(x) \in K$ for the remainder of this subsection, where $K = N^{k+3}(\Omega)$. The compactness of $H^{k+3}(\Omega)$ relative to $H$ is verified by the Rellich-Kondrachov theorem \cite[Chapter 6]{adams2003sobolev} and it then follows that $K$ is compact in $H$ since it is a subspace of $H^{k+3}(\Omega)$. This choice of $K$ is necessary to obtain $\mathcal O\p{h_n^{k+1}}$ convergence because higher order boundary compatibility conditions need to be satisfied to obtain optimal estimates in the semigroup approximation. Refer to Appendix \ref{section: error estimates} for further details.

We now write the state space representation of the model problem as
$$
\left\{
\begin{aligned}
    \frac{dz(\cdot, t)}{dt} &= Az(\cdot, t) + Bu(t) \\
             z(\cdot, 0) &= z_0,
\end{aligned}
\right.
$$
which is well-defined for all $z(\cdot, t)\in \mathcal D(A)$ and $u(t) \in \mathbb R$ since we have shown that the range of $A$ and $B$ are in $H$. The correspondent optimal feedback law is then given by
$$
    u_{opt}(t) = -\mathcal K_t z(\cdot, t),
$$
where $\mathcal K_t:=B^*P(t)$ and $P(t) \in \sJ_\infty^s$ satisfies the following weak time-dependent operator-valued Riccati equation
$$
\left\{
\begin{aligned}
    \frac{d}{dt}\p{\psi, P(t) \phi}_H &= 
    -\p{\psi, A^* \phi}_H
    - \p{\psi, A\phi}_H 
    + \p{\psi, \p{PBB^*P}(t)}_H 
    - \p{\psi, C^*C \phi}_H \\
    \p{\psi, P(\tau)\phi}_H &= \p{\psi, P_\tau \phi}_H
\end{aligned}
\right.
$$
for all $\phi, \psi \in \mathcal D(A)$. The theory discussed in \S\ref{section: time-dependent approximation} applies to this example problem, since $P_\tau, BB^*, C^*C \in \sJ_\infty^s \cap \sL(H,K)$, by definition.

\subsubsection{Finite Element Approximation with Piecewise Lagrange Elements} \label{finite element space}

One dimensional piecewise Lagrange polynomials \cite[Chapter 3, Section 1]{burden2015numerical} of arbitrary order $k \in \mathbb N$ are chosen as the approximation space for this model problem. The domain $\Omega$ is segmented into $n$ disjoint line segments $\mathscr T_n:= \Br{\ell_j}_{j=1}^n$ of equal length $h_n := \frac{1}{n}$. On $\mathscr T_n$, the piecewise Lagrange polynomial approximation space is then defined as
\begin{equation} \label{eqn: V def}
    V^k_n := \Br{ v \in C^0(\overline \Omega) : \left. v \right|_{\ell_j} \in \mathcal P^k(\ell_j) \quad \forall j=1,\ldots, n},
\end{equation}
where $C^0(\overline \Omega)$ denotes the space of continuous functions defined on $\overline \Omega$ and $\mathcal P^k(\ell_j)$ denotes the space of polynomials of order $k$ defined over $\ell_j$. The dimension of $V^k_n$ is then denoted as $n$ and we will denote $\Br{\phi_j}_{j=1}^{n}$ as an admissible basis of $V^k_n$. 

The finite element approximation of \eqref{eqn: example problem 2} is given by seeking a $z_n(\cdot, t) \in V^k_n$ that satisfies
\begin{equation} \label{eqn: finite element 2}
\left\{
\begin{aligned}
    \frac{d}{dt} \dual{z_n(\cdot,t), q_n} &= a\p{z_n, q_n} + \p{b(x)u_n(t), q_n} \\ 
    \dual{z_n(\cdot, 0) , q_n} &= \dual{z_0, q_n}
\end{aligned}
\right.
\end{equation}
for all $q_n \in V^k_n$ and $t \in I$, where  
$$
    a(q_n, r_n) := - \int_{-1}^1 \p{\alpha q_n'(x)r_n'(x) + \beta q_n(x) r_n(x)} dx 
$$
is the bilinear form associated with the diffusion operator $A$. Here, the homogeneous Neumann boundary condition is enforced naturally through integration-by-parts and setting the boundary  terms to zero. The solution $z_n(\cdot,t) \in V^k_n$ can then be represented as
$$
    z_n(x,t) := \sum_{j=1}^n \hat z_j(t) \phi_j(x),
$$
where $\Br{\phi_j}_{j=1}^n \in V^k_n$ is a set of functions representing the basis of the approximation space $V^k_n$. The coefficients $\hat{\mathbf z}(t) = \br{\hat z_1(t), \cdots \hat z_n(t)}$ is then the solution of the following matrix system
$$
\left\{
\begin{aligned}
    \frac{d}{dt} \hat{\mathbf z}(t) &= \br{\mathbf A_n} \hat{\mathbf z}(t) + \br{\mathbf B_n}u_n(t) \\
    \hat{\mathbf z}(0) &= \hat{\mathbf z}_0
\end{aligned}
\right.,
$$
where 
$$
    \br{\mathbf A_n} := \br{M_n}^{-1} \br{K_n} \textrm{ and }
    \br{\mathbf B_n} := \br{M_n}^{-1} \br{B_n},
$$
with 
$$
    \br{M_n}_{ij} := \int_{-1}^1 \phi_i(x)\phi_j(x) dx, \quad 
    \br{K_n}_{ij} := a(\phi_i, \phi_j),\textrm{ and } 
    \br{B_n}_i := \int_{-1}^1 b(x) \phi_i(x) dx.
$$
The resultant cost functional to be minimized is then given in matrix form as 
\begin{equation} \label{eqn: matrix thermal cost functional}
    \mathcal J(u_n) = 
    \hat{\mathbf{z}}^T(\tau) \br{\mathbf P_{\tau. n}} \hat{\mathbf{z}}(\tau)
    + \int_0^\tau \hat{\mathbf z}(t) \br{\mathbf C}^T\br{\mathbf C} \hat{\mathbf z}(t) + u_n^2(t) dt,
\end{equation}
where 
$$
\br{\mathbf P_{\tau, n}}_{ij} := \int_{-1}^1 \phi_i(x) \mathbf P_{\tau} \phi_j(x) dx \quad \textrm{ and } \quad 
\br{\mathbf C}_{i} = \br{\int_{-1}^1 q(x)\phi_i(x) dx, \ldots \int_{-1}^1 q(x)\phi_n(x) dx}.  
$$
This cost functional is derived by substituting $z(\cdot)$ with $z_n(\cdot)$ and $u$ with $u_n(\cdot)$ in \eqref{eqn: thermal cost functional} and recalling that $z_n(t) := \sum_{j=1}^n \hat z_j(t) \phi_j(x)$. 

The optimal feedback law for the finite element matrix control system is then given by 
$$
    u_{n,opt}(t) :=  \br{\mathcal K_{n,t}} \hat{\mathbf z}(t),
$$
where $\br{\mathcal K_{n,t}} := [\hat{\kappa}_1(t), \ldots, \hat{\kappa}_n(t)]$ is the vector representation of the functional gain given by
$$
    \br{\mathcal K_{n,t}} = - \br{\mathbf B}^T \br{\mathbf P_n(t)},
$$
and $\br{\mathbf P_n(t)} \in \mathbb R^{n\times n}$ is the solution of the following matrix Ricccati equation
$$
\left\{
\begin{aligned}
    \frac{d}{dt}\br{\mathbf P_n(t)}
    &= 
    -\br{\mathbf A_n}^T \br{\mathbf P_n(t)}
    -\br{\mathbf P_n(t)} \br{\mathbf A_n} 
    + \p{\br{\mathbf P_n(t)} \br{\mathbf B}^T \br{\mathbf B} \br{\mathbf P_n(t)}}
    -\br{\mathbf C}^T \br{\mathbf C} \\
    \br{\mathbf P_n}(\tau) &= \br{\mathbf P_{n, \tau}}
\end{aligned}
\right. 
$$
for all $t \in I$. The functional gain of the matrix system can then be defined as a function in $V^k_n$ by setting 
$\boldsymbol\kappa_n(t) := -\br{\mathbf B}^T \br{\mathbf P_n(t)} \br{\mathbf M_n}^{-1}$ and then defining $\kappa_n(t) := \sum_{j=1}^n \hat\kappa_j(t) \phi_j(x)$, where $\hat\kappa_j(t)$ are the elements of the vector $\boldsymbol \kappa_n(t)$. In functional form, we may write the functional gain as 
$$
\mathcal K_{n,t}(\cdot) := \dual{\kappa_n(t), \cdot}.
$$
With this, we move on to present the main result of this section.

\begin{lemma}
Let $H = L^2(\Omega)$ and $K = N^{k+3}(\Omega)$, then there exists a constant $M \in \mathbb R_+$ such that
$$
    \norm{\mathcal K_t - \mathcal K_{n,t}}_{\sC(I;\sL\p{H;\mathbb R)}} \leq M h_n^{k+1},
$$
where $\mathcal K_t$ and $\mathcal K_{n,t}$ are the functional gains associated with the optimal feedback laws for the infinite dimensional control system and its finite element approximation respectively. 
\end{lemma}
\begin{proof}
    For the purpose of analysis only, we may equivalently represent the finite element problem \eqref{eqn: finite element 2} as
    $$
    \left\{
    \begin{aligned}
        \frac{dz_n(\cdot,t)}{dt} & = A_n z_n(\cdot,t) + B_n u_n(t) \\
                    z_n(\cdot, 0) &= \pi_n z_0
    \end{aligned}
    \right.
    $$
    for all $t \in I$, where the operator $A_n: H_n \rightarrow H_n$ is defined by
    $$
    \begin{aligned}
        A_n v_n &= \sum_{j=1}^n \dual{\phi_j, A v_n}\phi_j \\
                    &= -\sum_{j=1}^n a\p{v_n, \phi_j} \phi_j
    \end{aligned}
    $$
    for all $v_n,\phi_i \in V^k_n$. Here, $\mathcal D(A_n) := \Br{v_n \in V^k_n : \partial_x v_n = 0 \textrm{ on } \Gamma}$ since homogeneous Neumann boundary condition is enforced implicitly as a natural condition. The actuator operator $B_n : \mathbb R \rightarrow V^k_n$ can be simply taken as the $L^2(\Omega)$ projection of $b$ into $V^k_n$, i.e.
    $$
    \begin{aligned}
        B_n\p{\cdot} &:= \pi_n B\p{\cdot} \\
                    &= \sum_{j=1}^n \dual{\phi_j, b(x)\p{\cdot}} \phi_j
    \end{aligned}
    $$
    where $\pi_n : H\rightarrow V^k_n$ is the orthogonal projection operator associated with the approximation space $V^k_n$. Likewise 
    $$
    \begin{aligned}
        C_n^*{\cdot} &:= \pi_n C^*\p{\cdot} \\
            &= \sum_{j=1}^n \dual{\phi_j, c(x) \p{\cdot} \phi_j}
    \end{aligned}
    $$    
    Then, the optimal feedback law for the state space representation of the finite element approximation can then be given as
    $$
        u_{n,opt}(t) = -\mathcal K_{n,t} z_n(t),
    $$
    where $\mathcal K_{n,t} := -B_n^*P_n(t)$ and $P_n(t) \in sC\p{I;\sJ_\infty^s}$ is the solution of the following operator-valued Riccati equation
    \begin{equation*}
    \left\{
        \begin{aligned}
        \frac{d}{dt} \p{\psi_n, P_n(t) \phi_n}_H &= 
        -\p{\psi_n, A_n^*P_n(t) \phi_n}_H -
        \p{\psi_n, A_n P_n(t) \phi_n}_H +
        \p{\psi_n, (P_n B_n B_n^* P_n)(t)\phi_n}_H -
        \p{\psi_n, C_n^*C_n \phi_n}_H \\
        \p{\psi_n, P_n(\tau) \phi_n}_H &= \p{\psi_n, P_\tau \phi_n}_H
        \end{aligned}
    \right.
    \end{equation*}
    for all $\phi_n, \psi_n \in V^k_n$ and $t \in I$. In the remainder of the proof, we will take $S(t)\in \sL\p{H}$ to be the semigroup generated by $A: \mathcal D(A) \rightarrow H$ and $S_n(t) \in \sL\p{H} := \iota_H \widehat S_n(t) \pi_n$, where $\widehat S_n(t) \in \sL\p{H_n}$ is the semigroup generated by $A_n: \mathcal D(A_n) \rightarrow V^k_n$.  
    
    We first utilize the work of \cite{fujita1976finite} to establish the pointwise convergence of $S_n(t)$ to $S(t)$, i.e.
    $$
        \lim_{n\rightarrow \infty} \norm{S(t)\phi - S_n(t) \phi}_H = 0 \quad \forall \phi \in H.
    $$
    Following this, we then apply Theorem \ref{theorem: m1} to observe that there exists a constant $C_1 \in \mathbb R_+$ so that
    $$
        \norm{P(\cdot) - P_n(\cdot)}_{\sC(I; \sJ_\infty)} \leq C_1 \sup_{t\in I} \norm{S(t) - S_n(t)}_{\sL(K,H)} \norm{\sQ(P)}_{\sC(I; \sL(H,K))}.
    $$
    Applying Lemma \ref{lemma: heat semigroup error} allows us to see that there exists a constant $C_1 \in \mathbb R_+$ so that  
    $$
        \norm{P(\cdot) - P_n(\cdot)}_{\sC(I; \sJ_\infty)} \leq C_1 h_n^{k+1} \norm{\sQ(P)}_{\sC(I; \sL(H,K))},
    $$
    where $\sQ(P) \in \sC(I; \sL(H,K))$ is verified by inspecting the definition of $B: \mathbb R \rightarrow N^{k+3}(\Omega)$ and the smoothness assumption made on $q(x) \in N^{k+3}(\Omega)$. Following this, we have that there exists a constant $C_2 \in \mathbb R_+$
    $$
        \norm{B - B_n}_{\sL(\mathbb R, H)} \leq C_2 h_n^{k+1} |b|_{H^{k+1}(\Omega)}
    $$
    after seeing that $\norm{v - \pi_n v}_{H} \leq Ch_n^{k+1}|v|_{H^{k+1}(\Omega)}$ for all $v \in K$ as a consequence of the one dimensional Lagrange interpolation error estimate given in \cite[Proposition 1.12]{ern2004theory}. Putting these error bounds together allows us to see that
    $$
    \begin{aligned}
        \norm{\mathcal K_t - \mathcal K_{n,t}}_{\sC\p{I; \sL(H,\mathbb R)}} 
        &=\norm{B^*P(\cdot) - B^*_nP_n(\cdot)}_{\sC\p{I; \sL(H,\mathbb R)}} \\ 
        &\leq\norm{P(\cdot) - P_n(\cdot)}_{\sC(I; \sJ_\infty)} \norm{B}_{\sL(\mathbb R; H)}
            + \norm{B - B_n}_{\sL(\mathbb R; H)}\norm{P_n(\cdot)}_{\sC(I; \sJ_\infty)} \\
        &\leq C_1 h_n^{k+1} \norm{\sQ(P)}_{\sC(I; \sL(H,K))} \norm{B}_{\sL(\mathbb R; H)}
            + C_2 h_n^{k+1} \norm{P_n(\cdot)}_{\sC(I; \sJ_\infty)}.
    \end{aligned}
    $$
    Thus proving the error bound presented in this Lemma.  
\end{proof}

This concludes the analytical discussion regarding the theoretical error estimate for the finite element approximation of the model thermal control system. We corroborate the error estimate for the functional gain in the computational results presented in below.

\subsubsection{Computational Results}
In this computational example, we set $I:= [0, 0.1]$, $\alpha = 1$, $\beta = 1$, and $b(x), q(x) = \exp\p{\frac{-1}{1-x^2}}$. With these parameters, the time-dependent matrix Riccati equation was solved using the trapezoidal time-stepping method \cite[\S 5.4]{burden2015numerical} with a time-step of $\Delta t = 10^{-4}$. The {\tt linalg.solve\_continuous\_are()} function in the SciPy numerical library was used to solve the matrix Riccati equation at every implicit time step. The terminal condition was taken to be $P_\tau = 0$. 

\begin{figure}[h!]
    \label{fig: Ex2 Results}
    \includegraphics[width = 0.65\textwidth]{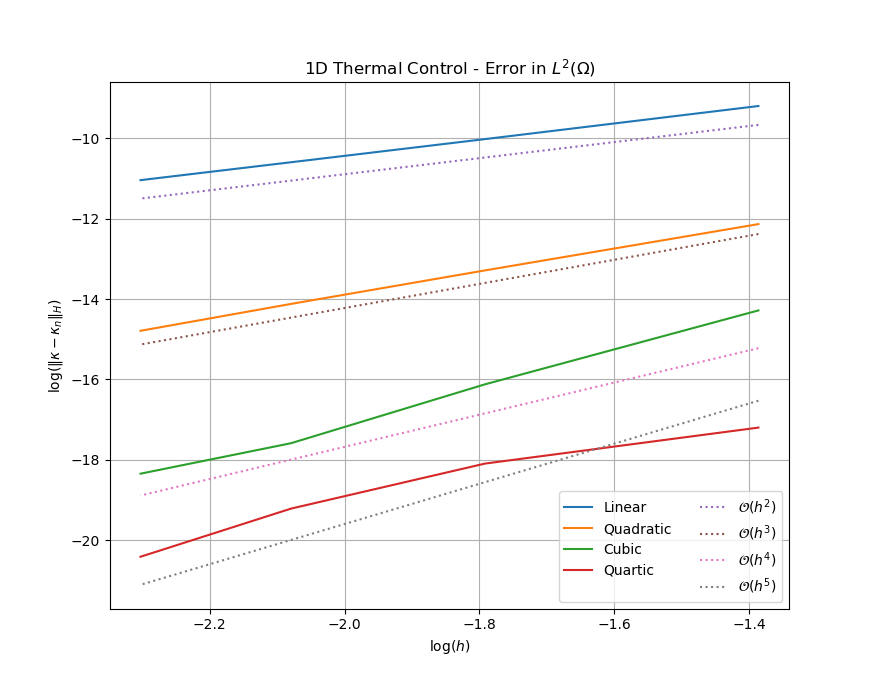}
    \caption{Approximation error in $\norm{\kappa - \kappa_n}_H$ for the finite element approximation of the thermal control system, where $H=L^2(\Omega)$.}
\end{figure}

The functional gain error is plotted in Figure \ref{fig: Ex2 Results} with respect to the $L^2(\Omega)$ norm, where we have compared the finite element solution to a spectral element solution on a single element with a 128-th order Lagrange polynomial basis defined on the Lobatto quadrature points. We observe good agreement between the numerically computed errors and our theoretical results as $h_n \rightarrow 0$.

We now move forward to discuss the error estimate for the finite element approximation of the two-dimensional generalization of this model problem.

\subsection{Infinite-Horizon LQR Control of a 2D Thermal Control Problem}
\label{subsection: 2D thermal problem}
In this numerical example, we consider the two-dimensional infinite-horizon analogue to the model problem presented in \S\ref{sec: ex2}. To that end, we define $\Omega:=[0, 1] \times [0,1]$ to be the problem domain. The model problem of interest is then given by
\begin{equation*}
     \min_{u(\cdot) \in L^2\p{(0, \infty); \mathbb R}} \int_0^\tau \br{\p{\int_\Omega q(\xb)z(\mathbf x,t)d\xb}^2 + r|u(t)|^2} dt,
\end{equation*}
where $z\p{\cdot} \in \sC(\mathbb R_+; H)$ is governed by the following set partial differential equation
\begin{equation}
    \left\{
    \begin{aligned}
        \partial_t z(\xb,t) &= \alpha \Delta z(\xb,t) - \beta z(\xb,t) + b(\xb)u(t)  & \quad \textrm{ on } \Omega \times \mathbb R_+ \\
        \partial_\nb z(\xb,t) &= 0  &\quad \textrm{ on } \Gamma \\
        z(\xb,0) &= z_0(\xb)    &\quad \textrm{ on } \Omega,
    \end{aligned}
    \right.
\end{equation}
In this problem, we choose $H := L^2(\Omega)$ and $K:= N^{k+3}(\Omega)$ and restrict $b,q$ to be elements of $K$. $\mathcal D(A)$ is then $H^2(\Omega) \cap N^0(\Omega)$ and $S(t) \in \sL\p{H}$ is the $C_0$-semigroup generated from $A(\cdot):= \alpha \Delta(\cdot) - \beta(\cdot)$. Defining $B(\cdot) := b(\xb) \cdot$ and $C(\cdot) := q(\xb)$, we have that the optimal feedback control law is given by the following 
$$
    u_{opt} := \mathcal Kz(t),
$$
where
$$
    \mathcal K(\cdot) := -\frac{1}{r}B^* P (\cdot)
$$
and $P \in \sJ_\infty$ is the solution to the following operator-valued Riccati equation. 
$$
   \p{\psi, A^* \phi}_H
    + \p{\psi, A\phi}_H 
    - \frac1r \p{\psi, PBB^*P\phi}_H 
    + \p{\psi, C^*C \phi}_H = 0
$$
for all $\phi,\psi \in \mathcal D(A)$. For brevity, the technical details on the finite element approximation of this model problem will be omitted since it is a simple extension to what has already been discussed previously in \S\ref{sec: ex1}. In short, this means that we also expect $\mathcal O\p{h^{k+1}}$ convergence for the $k$-order Lagrange element approximation of this model problem.

\begin{figure}
    \includegraphics[width = 0.65\textwidth]{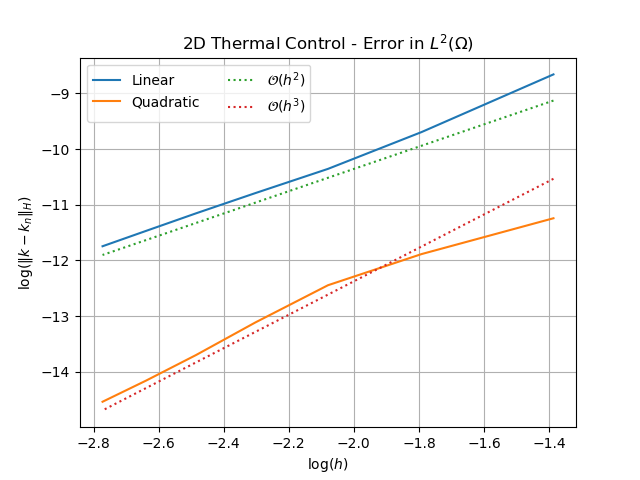}
    \caption{Approximation error in $\norm{\kappa - \kappa_n}_H$ for the finite element approximation of the two-dimensional thermal control system, where $H = L^2(\Omega)$. }
    \label{figure: 2D-thermal}
\end{figure}

In the following computational example, we have set $\alpha = 10^{-2}$ and $\beta = 1$. Additionally, we have set $b(\xb) = q(\xb) = \exp\p{\frac{-2}{1-(2x - 1)^2)}}\exp\p{\frac{-2}{1-(2y - 1)^2)}}$ and $r=10^{-4}$. The tensor product linear and quadratic elements were used on a square mesh consisting of squares with length and width $h_n := \frac1n$. The approximation error associated with the functional gain for this 2D problem is presented in Figure \ref{figure: 2D-thermal}, where optimal asymptotic convergence rates are observed. The {\tt FENICSx} \cite{BarattaEtal2023, ScroggsEtal2022, BasixJoss, AlnaesEtal2014} Python library was used to assemble the matrices used in the finite element approximation of the control system and the {\tt linalg.solve\_continuous\_are()} Python function was used to compute the solution to the matrix Algebraic Riccati Equation. 

\subsection{Infinite-Horizon LQR Control of a 1D Weakly Damped Wave System} \label{sec: ex1}
\subsubsection{Problem Definition}
    The domain for this problem is defined by $\Omega:= \p{-1,1}$ with the boundary defined as $\Gamma:= \Br{-1} \cup \Br{1}$. The model problem we approximate in this section is then given by
    \begin{equation} \label{eqn: damped wave cost functional}
		\min_{u(\cdot) \in L^2\p{(0,\infty); \mathbb R}} \int_0^\infty \left[\p{\int_{-1}^1 p_1(x) v(x,t) dx}^2 + \p{\int_{-1}^1 p_2(x) w(x,t) dx}^2 + \beta |u(t)|^2\right] dt
    \end{equation}
	subject to the 1D weakly damped wave system
	$$
	\left\{
	\begin{aligned}
		\partial_t v(x,t) &= w(x,t) + b_1(x)u(t)  &\textrm{ for } (x,t)\in(-1,1) \times (0,\infty)\\
		\partial_t w(x,t) &= c^2\partial_{x}^2 v(x,t) - \gamma w(x,t) + b_2(x) u(t) &\textrm{ for } (x,t) \in (-1,1) \times (0,\infty)\\
		v(x,0) &= v_0(x) &\textrm{ for } x \in (-1,1) \\
		w(x,0) &= w_0(x) &\textrm{ for } x \in(-1,1) \\
		v(-1,t) &= v(1,t) = 0 &\textrm{ for } t\in(0,\infty),
	\end{aligned}
	\right.
	$$
	for some arbitrary initial conditions $v_0\in H^1_0(\Omega)$, and $w_0 \in L^2(\Omega)$. The functions $p(x)$ and $b(x)$ are chosen to be $H^{k+3}(\Omega)$ smooth for this example and we require that their partial derivatives up to order $k+3$ are zero on $\Gamma$, i.e. $p(x),b(x) \in D^{k+3}(\Omega)$. This choice of $K$ is required to satisfy the higher order boundary compatibility conditions to derive an $\mathcal O\p{h_n^{k+1}}$ error estimate. See again Appendix \ref{section: error estimates} for details. Here $k \in \mathbb N$ will again denote the polynomial order of the approximation space used in the finite element approximation. 
 
    We can map this control problem into the state space form required by our theory by setting
	$$
		\mathbf z(x,t):= \left(\begin{array}{c} v(x,t)\\w(x,t) \end{array} \right), \qquad
        \mathcal A \mathbf z(\cdot, t) := \left[ \begin{array}{cc} 0 & v(x,t)  \\ c^2 \partial_{xx}w(x,t) & -\gamma v(x,t) \end{array}\right],
	$$
 	$$
	Bu(t):= \left[ \begin{array}{c} b_1(x) \\ b_2(x) \end{array} \right] u(t),\quad C \mathbf z(\cdot,t):= \br{\int_{-1}^1p_1(x) v(x,t)dx, \int_{-1}^1 p_2(x)w(x,t)dx }, \quad \textrm{ and } R:=\beta,
	$$
    We can write the above control system in the following state space form
    $$
        \min_{u(\cdot) \in L^2((0,\infty); \mathbb R)} \int_0^\infty \norm{C\mathbf z(\cdot, t)}_Y^2 + \p{u(t), R u(t)}_U dt  
    $$
    subject to 
    $$
    \left\{
    \begin{aligned}
        \frac{d\mathbf z(\cdot,t)}{dt} &= \mathcal A \mathbf z(\cdot,t) + Bu(t)\\
                \mathbf z(\cdot,0) &= \mathbf z_0,
    \end{aligned}
    \right.
    $$
    for all $t\in \mathbb R_+$, where we have chosen $U:= \mathbb R^2$, $Y:= \mathbb R$, $H:=H^1_0(\Omega) \times L^2(\Omega)$ and $\mathcal D(\mathcal A):= \p{H^2(\Omega)\cap H^1_0(\Omega)}\times H^1_0(\Omega)$. This choice of $H$ and $\mathcal D(\mathcal A)$ is required in order for $\mathcal A: \mathcal D(\mathcal A) \rightarrow H$ to generate an exponentially stable $C_0$-semigroup (e.g., see \cite{cox1994rate}). With this, it is clear that this example problem falls into the set of problems presented in \S\ref{sec: LQR}. The associated Riccati equation for this model problem is then given in the weak form by
    $$
        \p{\psi, \mathcal A^* P\phi}_H 
        + \p{\psi,P \mathcal A \phi}_H - 
         \p{\psi,PBR^{-1}B^*P\phi}_H 
        + \p{\psi, C^*C\phi}_H
        = 0
    $$
    for all $\phi,\psi \in H$. It is then clear that $BR^{-1}B^*$ and $C^*C$ are operators belonging to $\sJ_\infty^s \cap\sL(H,K)$, where $K := D^{k+3}(\Omega) \times H^{k+2}(\Omega)$. The choice of $K$ is derived from our earlier definition of the functions $q(x)$ and $b(x)$ as well as the analysis given in \S\ref{subsection: finite elment approximation of the weakly damped wave semigroup}. The optimal feedback law for this model control system is then given by
    \begin{equation} \label{eqn: functional gain 1}
        u_{opt}(t) = -\mathcal K \mathbf z(\cdot, t)
    \end{equation}
    for all $t \in \mathbb R_+$, where $\mathcal K: H \rightarrow \mathbb R$ is the functional gain defined as $\mathcal K:= R^{-1}B^*P$. 
    
\subsubsection{Finite Element Approximation with Piecewise Lagrange Elements} \label{section: finite element approximation}

The primary finite element approximation space used in this problem approximation is $V^k_n \subset H^k(\Omega)$ described in \S\ref{finite element space}. We additionally define the space
\begin{equation}
    V^k_{n,0} := V^k_n \cap H^1_0(\Omega),
\end{equation}
which consists of all functions belonging to $V^k_n$ that are zero on $\Gamma$. The approximation space associated with this model problem is given by 
$$
    H_n:= V^k_{n,0} \times V^k_{n,0}.
$$
An inspection indicates that $H_n \subset H$ and hence the orthogonal projection operator $\pi_n: H \rightarrow H_n$ is well defined. The choice of the approximation space $V^k_0$ as the approximation space for the $L^2(\Omega)$ variable $w(x,t)$ stems from the observation that $w = \partial_t v(x,t)$, and hence $w(\cdot, t)$ must also vanish on the boundary $\Gamma$ for all $t \in \mathbb R_+$.

In the finite element approximation of the weakly damped wave system, the approximate system is defined in the weak form by 
$$
\left\{
\begin{aligned}
    \frac{d}{dt} \dual{v_n, q_n} &= \dual{w_n, q_n} + \dual{b_1 u_n(t), q_n}, \quad &\dual{v_n(0), q_n} = \dual{v_0, q_n} \quad &\forall q_n \in V^k_{n,0} \\
    \frac{d}{dt} \dual{w_n, r_n} &= a(v_n, r_n) - \gamma \dual{w_n, r_n} + \dual{b_2u_n(t), r_n}, \quad &\dual{w_n(0), r_n} = \dual{w_0, r_n} \quad &\forall r_n \in V^k_{n,0},
\end{aligned}
\right.
$$
where 
$$
    v_n(x,t) = \sum_{j=1}^{n} \hat v_j(t) \phi_j(x) \qquad \textrm{ and } \qquad 
    w_n(x,t) = \sum_{j=1}^{n} \hat w_j(t) \psi_j(x)
$$
are the finite element approximation of $v(t)$ and $w(t)$ respectively with $\Br{\phi_j, \psi_j}_{j=1}^n \in V^k_{n,0}\times V^k_{n,0}$ are the basis of $H_n$ and $\hat v_i(t), \hat w_i(t) \in \mathbb R$ as time-dependent coefficients, and 
$$
    a(\phi,\psi) := -\int_{-1}^1 c^2 \phi'(x) \psi'(x)dx
$$
is a bilinear form mapping any $\phi, \psi \in H^1_0(\Omega)$ to $\mathbb R$. It then follows by the finite dimensional nature of the space $V^k_{n,0}$, and since $\phi_j$ and $\psi_j$ are both elements of $V^k_{n,0}$, that the finite element description of the approximate problem can be written in the following equivalent matrix form
$$ 
\left\{
\begin{aligned}
    \relax \br{M_n} \frac{d}{dt} \hat{\mathbf v} &= \br{M_n} \hat{\mathbf w}(t) + \br{b_{1,n}}, \quad &\hat{\mathbf v}(0) = \hat{\mathbf v}_0 \\
    \relax \br{M_n} \frac{d}{dt} \hat{\mathbf w} &= \br{K_n} \hat{\mathbf v}(t) - \gamma \br{M_n} \hat{\mathbf w}(t) + \br{b_{2,n}}u_n(t) \quad  &\hat{\mathbf w}(0) = \hat{\mathbf w}_0,
\end{aligned}
\right.
$$
where for $i,j=1,\ldots, n$,
$$
\begin{aligned}
\br{M_n}_{ij} = \int_{-1}^1 \phi_i(x) \phi_j(x) dx, \quad
\br{K_n}_{ij} = \int_{-1}^1 \phi_i'(x) \phi'_j(x) dx, \\
[b_{1,n}]_i = \int_{-1}^1 b_1(x) \phi_i(x) dx, \quad [b_{2,n}]_i = \int_{-1}^1 b_2(x) \phi_i(x) dx,
\end{aligned}
$$
with $\br{b_{1,n}}, \br{b_{2,n}}: \mathbb R^{n} \rightarrow \mathbb R$ being column vectors. Denoting
$$
    \br{\mathbf M_n} := \br{\begin{array}{cc} \br{M_n} & 0 \\ 0 & \br{M_n}  \end{array}}, \quad 
    \br{\mathbf K_n} := \br{
    \begin{array}{cc}
        0 & \br{M_n} \\
        \br{K_n} & -\gamma\br{M_n} 
    \end{array}
    }, 
    \quad
    \textrm{ and } 
    \quad 
    \br{\mathbf b_n} := \br{\begin{array}{c} \br{b_{1,n}} \\ \br{b_{2,n}}\end{array}},
$$
we may write the finite dimensional problem in the following block form
\begin{equation} \label{eqn: block matrix form}
\left\{
\begin{aligned}
    \frac{d}{dt}\hat{\mathbf z} &= \br{\mathbf A_n}\hat{\mathbf z} + \br{\mathbf B_n} u_n(t) \\
    \hat{\mathbf z}(0) &= \hat{\mathbf z}_0,
\end{aligned}
\right.
\end{equation}
where $\hat{\mathbf z}(t) = \br{\hat{\mathbf v}(t), \hat{\mathbf w}(t)}^T$, $\br{\mathbf A_n} = \br{\mathbf M_n}^{-1} \br{\mathbf K_n}$, and $\br{\mathbf B_n} = \br{\mathbf M_n}^{-1} \br{\mathbf b_n}$. 

With this, we wish to choose $u_n(\cdot)$ that minimizes following cost functional
\begin{equation} \label{eqn: matrix damped wave cost functional}
    \int_0^\infty \hat{\mathbf z}^T(t)\br{\mathbf C_n}^T\br{\mathbf C_n}
    \hat{\mathbf z}(t) + Ru_n^2(t) dt,
\end{equation}
where $\br{\mathbf C_n} := \br{\br{q_{1,n}}, \br{q_{2,n}}}$ is a row vector with $\br{q_{1,n}}_{i} = \int_{-1}^1 q_1(x) \phi_i(x) dx$ and $\br{q_{2,n}}_{i} = \int_{-1}^1 q_2(x) \phi_i(x) dx$  for $i=1,\ldots,n$. The matrix form of the cost functional given in \eqref{eqn: matrix damped wave cost functional} arises by replacing $v$ with $v_n$ and $u$ with $u_n$ in \eqref{eqn: damped wave cost functional}. The optimal feedback law for the finite element approximation of the control system is then given in matrix form by
$$
    u_{n,opt}(t) = -\br{\mathcal K_n} \hat{\mathbf z}(t),
$$
where $\br{\mathcal K_n}:= R^{-1} \br{\mathbf B_n}^T \br{\mathbf P_n}$ and $\br{\mathbf P_n}$ is the solution of the following matrix Riccati equation
\begin{equation} \label{eqn: matrix riccati equation}
    \br{\mathbf A_n}^T \br{\mathbf P_n} + \br{\mathbf P_n} \br{\mathbf A_n} 
    - \br{\mathbf P_n}\br{\mathbf B_n} R^{-1} \br{\mathbf B_n}^T \br{\mathbf P_n}
    + \br{\mathbf C_n}^T \br{\mathbf C_n} = 0.
\end{equation}

We now wish to represent the optimal feedback law in the functional sense, where 
$$
    u_{n, opt}(t) := -\mathcal K_n \mathbf z_n(\cdot, t)
$$
is the approximate functional gain represented as an element of $H$, by first defining $\kappa_n := \br{\sum_{j=1}^{n} \hat{\kappa}_j^1 \phi_j(x), \sum_{j=1}^{n} \hat{\kappa}_j^2 \phi_j(x)}^T$, where the coefficients $\hat \kappa_j^1, \hat \kappa_j^2  \in \mathbb R$ are computed using
$$
    \br{{\boldsymbol\kappa}} = R^{-1}[\mathbf B_n]^T[\mathbf P_n][\mathbf M_n]^{-1},
$$
where we have denoted $\br{\boldsymbol\kappa}:= \br{\hat\kappa_1^1,\ldots,\hat\kappa_{N}^1, \hat\kappa_1^2,\ldots, \hat\kappa_{N}^2}^T$. The functional form is then given by setting
$$
    \mathcal K_n\p{\cdot}:= \dual{\kappa_n, \cdot}.
$$
With this, we can now prove that $\mathcal K_n$ converges to $\mathcal K$ with an optimal order of convergence. 
\begin{lemma}
    Let $H=H^1_0(\Omega) \times L^2(\Omega)$ and $K=D^{k+3} \times H^{k+2}(\Omega)$ for all $k\in \mathbb N$. If $b(\cdot),q(\cdot) \in K$, then there exists a constant $M \in \mathbb R_+$ so that
    $$ 
        \norm{\mathcal K - \mathcal K_n}_{\sL\p{H, \mathbb R}} \leq Mh_n^{k} \norm{\mathcal Q(\Sigma)}_{\sL\p{H, K}},
    $$
    where $\mathcal K, \mathcal K_n \in \sL(H; \mathbb R)$ are the functional gains associated with the optimal feedback control law for the infinite dimensional control system \eqref{eqn: functional gain 1} and its finite element approximation respectively. 
\end{lemma}
\begin{proof}  
    For the purpose of analysis, we may represent the finite element approximation of the weakly damped wave system in the following form,
    $$
        \left\{
        \begin{aligned}
        \frac{d\mathbf z_n(\cdot,t)}{dt} &= \mathcal A_n\mathbf z_n(\cdot,t) + B_n u(t) \\
        \mathbf z_n(\cdot,0) &= \mathbf z_{n,0},
        \end{aligned}
        \right.
    $$
    for all $t\in \mathbb R_+$, where we have defined
    $$ 
        \mathcal A_n \mathbf z_n(\cdot, t) := \left[ \begin{array}{cc} 0 & v_n(x,t)  \\ c^2 \sum_{j=1}^n a(w_n(\cdot,t),\phi_j)\phi_j\p{x} & -\gamma v_n(x,t) \end{array}\right],
    $$
    $$
        B_nu_n(t):= \left[ \begin{array}{c} \sum_{j=1}^n \dual{\phi_j, b_1(\cdot)}\phi_j(x) \\ \sum_{j=1}^n \dual{\phi_j, b_2(\cdot)}\phi_j(x)
        \end{array} \right] u_n(t),
    $$
    and
    $$
        C_n \mathbf z_n(\cdot,t):= \br{\int_{-1}^1 \sum_{j=1}^n \dual{\phi_j, p_1(\cdot)}\phi_j\p{x} v_n(x,t) dx, \int_{-1}^1 \sum_{j=1}^n \dual{\phi_j, p_2(\cdot)}\phi_j\p{x} v_n(x,t) dx}.
    $$
    Note that we have defined $\mathbf z_n(x,t) := \br{v_n(x,t), w_n(x,t)}^T$. Then, the optimal feedback law for the state space representation of the finite element approximation can then be given as
    $$
        u_{n,opt}(t) = -\mathcal K_n \mathbf z_n(\cdot, t),
    $$
    where $\mathcal K_{n} := -B_n^*P_n$ and $P_n \in \sJ_\infty^s$ is the solution of the following operator-valued Riccati equation
    \begin{equation*}
        \begin{aligned}
        \p{\psi_n, A_n^*P_n(t) \phi_n}_H +
        \p{\psi_n, A_n P_n(t) \phi_n}_H -
        \p{\psi_n, (P_n B_n B_n^* P_n)(t)\phi_n}_H +
        \p{\psi_n, C_n^*C_n \phi_n}_H = 0\\
        \end{aligned}
    \end{equation*}
    for all $\phi_n, \psi_n \in H_n$.

    From the definition of $\mathcal K, \mathcal K_n: H\rightarrow \mathbb R$, we see that
    $$
    \begin{aligned}
        \norm{\mathcal K - \mathcal K_n}_H 
        &\leq |R^{-1}| \norm{B^*P - B^*_n P_n}_{\sL(H,\mathbb R)} \\
        &\leq |R^{-1}| \p{\norm{B^*}_{\sL(H, \mathbb R)} \norm{P - P_n}_{\sJ_\infty} + \norm{B^* - B_n}_{\sL(H, \mathbb R)}\norm{P_n}_{\sJ_\infty}}.
    \end{aligned}
    $$
    We first observe that 
    $$
        \lim_{n\rightarrow \infty}\norm{S(t)\mathbf q - S_n(t)\mathbf q}_H = 0
    $$ 
    for all $\mathbf q \in H$ as a consequence of Lemma \ref{eqn: damped wave pointwise convergence}. This then allows us to apply Theorem \ref{theorem: m2} to see that 
    $$
    \begin{aligned}
        \norm{P-P_n}_{\sL(H)} 
        &\leq \gamma_n \norm{\sQ(P)}_{\sL(H,K)} \int_0^\infty e^{-\alpha s} \norm{S(s) - S_n(s)}_{\sL(K,H)} ds \\
        &\leq C_\gamma h^k \norm{\sQ(P)}_{\sL(H,K)} \int_0^\infty e^{-\alpha s} ds \\
        &=C_{\gamma, \alpha} h^k \norm{\sQ(P)}_{\sL(H,K)}
    \end{aligned}
    $$
    as a consequence of Lemma \ref{lemma: wave semigroup error}. We observe that $\sQ(P)$ is bounded in $\sL(H, K)$ since we have assumed that $b(s)$ and $q(x)$ are elements of $K$. 

    We continue by seeing that
    $$
        \begin{aligned}
            \norm{B - B_n}_{\sL(H,\mathbb R)} &= \norm{\br{b_1(x) - \pi_n b_1(x), b_2(x) - \pi_n b_2(x)}^T}_H \\
            &=\norm{b_1(x) - \pi_n b_1(x)}_{L^2(\Omega)} + \norm{b_2(x) - \pi_n b_2(x)}_{L^2(\Omega)} \\
            &\leq Ch_n^{k+1}\p{|b_1|_{H^{k+1}(\Omega)} + |b_2|_{H^{k+1}(\Omega)}},
        \end{aligned}
    $$
    after bounding the error associated with the $L^2(\Omega)$ projection onto $V^k_n$ by the error associated with one-dimensional Lagrange interpolation. (\cite[Proposition 1.12]{ern2004theory}). With this, it becomes clear that there exists a constant $M \in \mathbb R_+$ so that
    $$
        \norm{\mathcal K - \mathcal K_n}_{H} \leq Mh_n^{k}.
    $$
    This concludes the proof.
\end{proof}

\subsubsection{Computational Results}
\begin{figure}[h!] \label{fig: Ex1 Results}
    \includegraphics[width=0.47\textwidth]{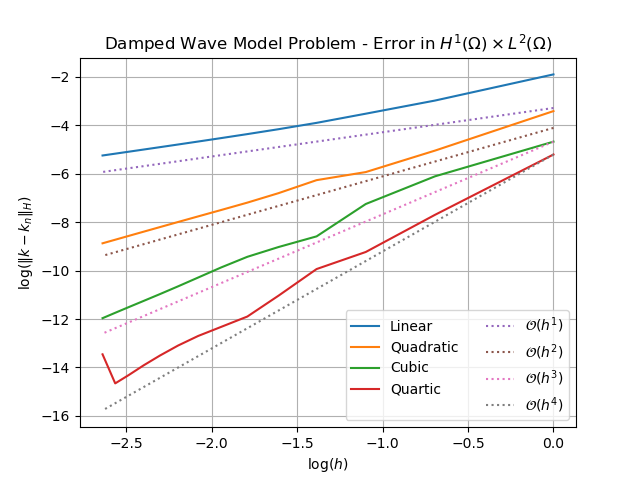}
    \includegraphics[width=0.47\textwidth]{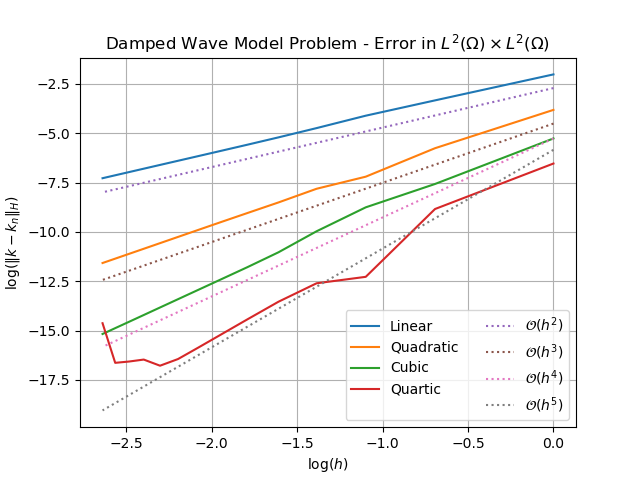}
        \caption{Approximation error in $\norm{\kappa - \kappa_n}_H$ for the finite element approximation of the weakly damped wave control system, where $H=H^1_0(\Omega) \times L^2(\Omega)$ in the left plot and $H=L^2(\Omega)\times L^2(\Omega)$ in the right figure.}
\end{figure}

In the computational illustration, the author's personal code was utilized to generate the operator matrices defined above. The {\tt linalg.solve\_continuous\_are()} function in SciPy was utilized to solve \eqref{eqn: matrix riccati equation}. We set $c = 1$, $\gamma = 10^{-4}$, and $b_2(x)= q_1(x) := \exp\p{\frac{-1}{1-x^2}}$ and $b_1(x) = q_2(x) = 0$. This choice of $b_1,b_2, q_1, q_2$ are chosen to model the system with physically possible actuators while minimizing the easily observable vertical displacement. With these parameters, the matrix Riccati equation was solved and the functional gain was computed. The functional gain error is plotted in Figure \ref{fig: Ex1 Results} with respect to the $H^1_0(\Omega)\times L^2(\Omega)$ norm and $L^2(\Omega) \times L^2(\Omega)$ norm, where we have compared the finite element solution to a spectral element solution on a single element with a 128$^{th}$ order Lagrange polynomial basis defined on the Lobatto quadrature points. We observe good agreement between the numerically computed errors and our theoretical results, albeit with numerical precision error past $10^{-8}$.

\subsection{Violation of Theoretical Conditions}
We have observed in our model examples that the finite element approximation of the functional gain associated with the linear quadratic regulator converges optimally with respect to the classical finite element theory. In each example, we have chosen the actuator and weighting functions to be smooth functions that vanish on the boundary. We will now demonstrate in the following two examples what happens computationally when this condition is violated.

\begin{figure}[h!]
    \includegraphics[width = 0.47\textwidth]{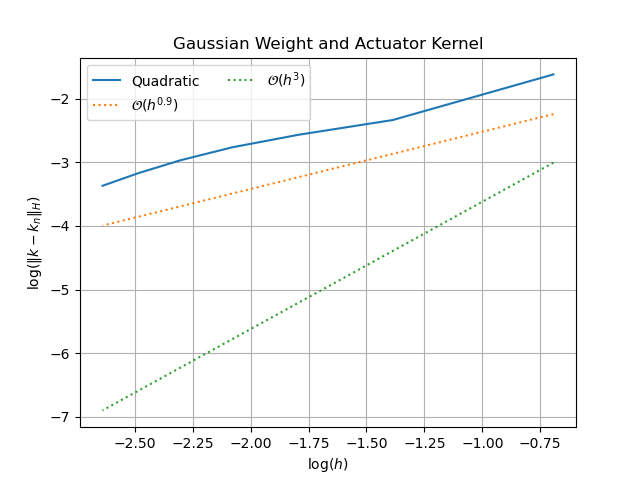}
    \includegraphics[width = 0.47\textwidth]{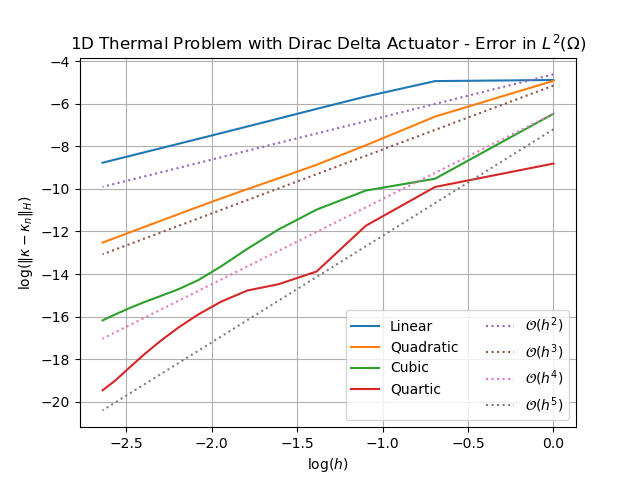}
    \caption{\textbf{Left:} Error of the functional gain approximation for the two-dimensional thermal model problem for when $b,q$ are Gaussian functions. \textbf{Right:} Error of the functional gain approximation for the 1D thermal model problem for when $b(x) = \delta(x)$.}
    \label{fig: counter_ex}
\end{figure}

First, we consider the model problem of \S \ref{subsection: 2D thermal problem}, where we take instead $b(\xb) = q(\xb) = \exp\p{-x^2 - y^2}$. From inspection, it becomes clear that $b$ and $q$ are not elements of $K=D^{k+3}(\Omega)$ since they are are nonzero everywhere on $\mathbb R^2$. A computational convergence experiment with quadratic elements indicates that the convergence rate of the functional gain approximation in $L^2(\Omega)$ converges suboptimally on the order of $\mathcal O(h^{0.9})$. The plot of the error is presented in the left of  Figure \ref{fig: counter_ex}. 

Next, we demonstrate that there is one particular violation that allows for optimal convergence rates. Here, we consider the one-dimensional thermal problem presented in \S\ref{sec: ex2} with $q(x) := \exp\p{\frac{-1}{1-x^2}}$ and $b(x) = \delta(x)$, where $\delta(x)$ denotes the Dirac delta function. This function clearly violates the conditions specified by the theory presented in this work since it is not a bounded function in any subspace of $L^2(\Omega)$. However, we observe optimal convergence for linear, quadratic, cubic, and quartic elements. The plot of the error is presented in the right of Figure \ref{fig: counter_ex}. The observed optimal convergence rates indicates to us that the functional gain associated with unbounded control systems (e.g. boundary and point control) can be optimally approximated. The reason behind this lies in the fact that the solution to the operator-valued Riccati equation maps elements in $K'$ to $K$ when the coefficient operator $F\in \sL\p{K', K}$. The full details behind the analysis of the approximation of operator-valued Riccati equations for when the control operator $G$ is unbounded is the topic of an upcoming work. 

\section{Conclusion}

In this work, we have presented an abstract theory for the approximation of operator-valued Riccati equations. This theory indicates that the error of the approximate solution, in its time-dependent and time-independent forms, is bounded by the semigroup approximation error. This is significant because it allows for the use of the existing approximation theory to derive error estimates. Indeed, we have derived error estimates for the solution approximation to the operator-valued Riccati equation using the classical estimates in our model problems. Additionally, our theory only requires the condition of (1) pointwise semigroup convergence, and (2) compactness in the coefficient operators to guarantee convergence of the approximate solution. This is two fewer conditions than what is considered necessary and sufficient in the classical set of convergence conditions. These observations imply that fewer steps are necessary establish convergence of the approximation and to derive error estimates using our general theory. Thus, our general theory simplifies the once arduous approximation theory for operator-valued Riccati equations. 

We have then demonstrated that the numerical error obtained from our computational experiments match the error estimates derived from our theory. To be specific, we have derived and verified $\mathcal O(\epsilon)$ convergence for the asymptotically perturbed scalar system, $\mathcal O(h^{k+1})$ convergence for the thermal system in the $\sL\p{L^2(\Omega)}$ norm, and $\mathcal O(h^{k})$ convergence for the damped wave system in the $\sL\p{H^1(\Omega) \times L^2(\Omega)}$ norm. This computational evidence strongly backs up our derived error estimates and the general theory presented in this work. 

While we believe that the theory presented in this work has wide-reaching implications, we are also aware of several theoretical shortcomings. The most notable shortcoming is the treatment of unbounded operators. These operators typically occur in e.g. boundary and point control and observations systems. While this is a clear violation to the compactness condition we have imposed on the coefficient operators in this work, we have observed optimal convergence rates for when the operator associated with the control actuator is unbounded in a thermal system. This indicates to us that our theory can be to be extended to the case of unbounded control operators. This is the topic on an ongoing study and will be released in a future work. 

\section*{Acknowledgements}
We would like to acknowledge the anonymous reviewers of this manuscript for their helpful suggestions. Their feedback pushed the author to significantly improve the quality of this work. Additionally, the author would like to thank Prof. John A. Burns for his mentorship in control theory. 

\bibliographystyle{plain}
\bibliography{bibliography}

\appendix
\section{Technical Lemmas} \label{sec: technical appendix}
In this section, we present the technical lemmas that we will use throughout the analysis presented \S\ref{section: time-dependent approximation} and \S\ref{section: time-independent approximation}. We begin by demonstrating that $\sL(H,K)$ forms a two-sided *-ideal (stated in a different manner).

\begin{lemma} \label{lemma: t1}
Let $H$ be a separable Hilbert space and $K$ be a compact subspace of $H$. In addition, let $U,V \in \sL\p{H}$ and $X\in \sL\p{H,K}$, Then $UXV\in \sL\p{H,K}$. 
\end{lemma}
\begin{proof}
We begin by seeing that
$$
    \norm{UXV}_{\sL\p{H,K}} \leq \norm{UX}_{\sL\p{H,K}} \norm{V}_{\sL\p{H}}.
$$
It then follows by Proposition \ref{proposition: properties of compact operators} that
$$
    \norm{UX}_{\sL\p{H,K}} = \norm{X^*U^*}_{\sL\p{K,H}} \leq \norm{X^*}_{\sL\p{H}} \norm{U^*}_{\sL\p{K,H}} \leq \norm{X}_{\sL\p{H}} \norm{U^*}_{\sL\p{H}}.
$$
Therefore 
$$
    \norm{UXV}_{\sL\p{H,K}} \leq \norm{U}_{\sL\p{H}} \norm{X}_{\sL\p{H}} \norm{V}_{\sL\p{H}}.
$$
The right hand side of the above inequality is bounded because $\sL\p{H,K}\subset\sL\p{H}$. The boundedness of $UXV$ in the $\sL\p{H,K}$ norm implies that $UXV \in \sL\p{H,K}$.  
\end{proof}

The following result is then easily proven.

\begin{corollary} \label{corollary: c1}
    Let $H$ be a separable Hilbert space and $K$ be a compact subspace of $H$, and $I$ be a compact subset of $\mathbb R_+$. Additionally, let $U(\cdot), V(\cdot) \in \sC\p{I; \sL(H)}$ and $X(\cdot) \in \sC\p{I; \sL(H,K)}$. Then $(UXV)(\cdot) \in \sC\p{I; \sL(H,K)}$. 
\end{corollary}
\begin{proof}
    For any $t \in I$, we have that $U(t), V(t) \in \sL(H)$ and $X \in \sL(H,K)$. A direct application of Lemma \ref{lemma: t1} allows us to see that $(UXV)(t) \in \sL(H,K)$. The result of this corollary is proven after recalling the definition of $\sC(I; \sL(H))$ and $\sC(I; \sL(H,K))$ and recalling that the continuous functions form a Banach algebra \cite[Theorem 4.12]{adams2003sobolev}.
\end{proof}

We conclude this section with the following uniform convergence result.

\begin{lemma} \label{lemma: t2}
Let $H$ be a separable Hilbert space and $S(t) \in \sL\p{H}$ be a $C_0$-semigroup. In addition let $\left\{S_n(t)\right\}_{n=1}^\infty \in \sL\p{H}$ be an approximating sequence of $C_0$-semigroups so that 
$$
    \lim_{n \rightarrow \infty} \norm{S(t)\phi - S_n(t)\phi}_H = 0 \quad \forall \phi \in H.
$$
Then we have that
\begin{equation} \label{eqn: a1}
    \lim_{n \rightarrow \infty} \norm{\p{S(t) - S_n(t)}X}_{\sJ_\infty} = 0
\end{equation}
and
\begin{equation} \label{eqn: a2}
    \lim_{n \rightarrow \infty} \norm{X\p{S(t) - S_n(t)}^*}_{\sJ_\infty} = 0
\end{equation}
for all $X\in \sJ_\infty$ and $t \in \mathbb R_+$.
\end{lemma}
\begin{proof}
 From the definition of the operator norm, we have that
 $$
    \norm{\p{S(t) - S_n(t)}X}_{\sJ_\infty} = \sup_{\phi \in B_1(0)} \norm{\p{S(t) - S_n(t)}X\phi}_H,
 $$
 where we have denoted $B_1(0)$ as the unit ball in $H$. It then follows that
 $$
    \norm{\p{S(t) - S_n(t) }X\phi}_H = \sup_{\psi \in X\p{B_1(0)}} \norm{\p{S(t) - S_n(t)} \psi}_H.
 $$
 Since $X \in \sJ_\infty$ is a compact operator, we have that $X\p{B_1(0)}$ is a compact set in $H$. The compactness of $X\p{B_1(0)}$ then implies that 
 $$ 
    \lim_{n\rightarrow \infty} \sup_{\psi \in X\p{B_1(0)}} \norm{\p{S(t) - S_n(t)} \psi}_H = 0
 $$
 because bounded operators converge uniformly on compact sets as a consequence of the Uniform Boundedness Principle \cite[Theorem 5.3-2]{ciarlet2013linear}. This convergence property is satisfied for all $t \in \mathbb R_+$ because of the pointwise convergence assumption made in the statement of the lemma and hence we have proven \eqref{eqn: a1}.
 
 Assertion \eqref{eqn: a2} then follows by seeing that $\norm{X\p{S(t) - S_n(t)}^*}_{\sJ_\infty} = \norm{\p{S(t) - S_n(t)}X^*}_{\sJ_\infty}$ and then applying \eqref{eqn: a1} after seeing that $X^*$ is also an element of $\sJ_\infty$.

\end{proof}

\section{Error Estimates for Finite Element Approximation of Semigroups} \label{section: error estimates}
In this section, we will present results on the error estimates for the finite element approximation of the heat and the weakly damped wave equation for when the solution operator of these equations are represented as a $C_0$-semigroup. We will begin by analyzing the heat semigroup approximation error and then move on to analyze the approximation error for the weakly damped wave semigroup. For simplicity, \emph{we will consider only the case where the problem domain is exactly represented} in the finite element meshing process to avoid the theoretical technicalities regarding geometric approximation. See e.g. \cite{strang1974analysis, ciarlet2002finite} for a discussion on how curved domains may be theoretically addressed. 

\subsection{Finite Element Approximation of the Thermal Semigroup}
In this analysis, we will carry over the notation used in \S\ref{sec: ex2}. Recall that we have defined $H:= L^2(\Omega)$ for the thermal model problem, where $\Omega$ is an $d$-dimensional domain. It is known that $A: \mathcal D(A) \rightarrow H$ generates an exponentially stable semigroup $S(t) \in \sL(H)$ \cite[\S 7.4]{evans2022partial}. Letting $z(\cdot, t) = S(t) g$, for any $g \in H$, we have that $z\p{\cdot,\cdot}\in \sC\p{I;H}$ is the solution of the following equation
\begin{equation} \label{eqn: heat equation}
\left\{
\begin{aligned}
    \partial_t z(\xb,t) &= \alpha \Delta z(\xb,t) - \beta z(\xb,t) &\textrm{ in } \Omega\times I \\
            \partial_\nb z(\xb, t) &= 0           &\textrm{ on } \Gamma \\ 
            z(\xb,0) &= g(\xb) &\textrm{ in } \Omega,
\end{aligned}   
\right. 
\end{equation}
where $\Delta \p{\cdot} := \sum_{j=1}^d \frac{\partial^2 (\cdot)}{\partial x_i^2}$ is the $d$-dimensional Laplacian operator, $\nb$ denotes the unit outer normal vector with respect to $\Gamma$, and $\partial_{\nb}(\cdot)$ denotes the normal derivative operator. 

Additionally, we have that $A_n: \mathcal D(A_n) \rightarrow V^k_n$ also generates an exponentially stable semigroup $\widehat S_n(t)\in \sL\p{V^k_n}$ since $V^k_n$ is a conforming subspace embedded in $H^1_0(\Omega)$. We extend the domain and range of $\widehat S_n(t)$ through an injection and projection process, i.e., $S_n(t) := \iota_H \widehat S_n(t) \pi_n$. Letting $z_n(\cdot, t) = S_n(t) g$, for any $g \in H$, we have that $z_n\p{\cdot,\cdot}\in \sC\p{I; V^k_n}$ is the solution of the following weak equation
$$
\left\{
\begin{aligned}
    \frac{d}{dt}\dual{z_n(\cdot, t), r_n} &= -a(z_n(\cdot, t), r_n) \\
                            \dual{z_n(\cdot, 0), r_n} &= \dual{g, r_n}
\end{aligned}
\right.
$$
for all $r_n \in V^k_n$ and $t \in I$. 

Recall that we have defined
$$
    N^k(\Omega):= \Br{v \in H^{k}(\Omega) : \partial_\xb^{\alphab} \partial_\nb v = 0 \textrm{ on } \Gamma, \quad |\alphab|=0,\ldots,k }
$$
as the space of functions whose partial derivatives up to order $k$ satisfies the homogeneous Neumann condition. With these notions defined, we are now ready to begin deriving an error bound on $S_n(t) \in \sL(H)$. We start with the following.
\begin{lemma} \label{lemma: ll1}
    Let $S(t)\in \sL(H)$ be the parabolic semigroup generated by $A: \mathcal D(A) \rightarrow H$, then we have that
    $$
        \partial_{\xb}^{\alphab} S(t) g = S(t) \partial_{\xb}^{\alphab} g
    $$
    for all $g \in N^k(\Omega)$. Furthermore, we have that $z(\cdot, t) \in \sC\p{I; N^{k}(\Omega)}$.
\end{lemma}
\begin{proof}
    Taking the $k^{th}$ order partial derivative with respect to $\xb \in \Omega$ on both sides of \eqref{eqn: heat equation} results in 
    \begin{equation*} 
\left\{
\begin{aligned}
    \partial^{\alphab}_\xb \partial_t z(\xb,t) &= \alpha \partial^{\alphab}_\xb \Delta z(\xb,t) - \beta z(\xb,t) &\textrm{ in } \Omega\times I \\
            \partial^{\alphab}_\xb \partial_\nb z(\xb, t) &= 0           &\textrm{ on } \Gamma \\ 
            \partial^{\alphab}_\xb z(\xb,0) &=  \partial^{\alphab}_\xb g(\xb) &\textrm{ in } \Omega.
\end{aligned}   
\right. 
\end{equation*}
Switching the order of differentiation for all the differentiation operations in the equation then allows us to see that 
\begin{equation}
\left\{
\begin{aligned}
    \partial_t \xi(\xb,t) &= \alpha \Delta \xi(\xb,t) - \beta \xi(\xb,t) &\textrm{ in } \Omega\times I \\
            \partial_\nb\xi(\xb,t) &= 0           &\textrm{ on } \Gamma \\ 
            \xi(\xb,0) &= \partial^{\alphab}_\xb g(\xb) &\textrm{ in } \Omega,
\end{aligned}   
\right. 
\end{equation}
where we have denoted $\xi(\xb,t) = \partial^{\alphab}_\xb z(\xb,t)$. Since $S(t) \in \sL\p{H}$ is the solution to the homogeneous Neumann initial boundary value problem, it then follows that 
\begin{equation} \label{eqn: thermal semigroup commutation}
    \partial^{\alphab}_\xb S(t)g = \partial^{\alphab}_\xb z(\xb,t) = \xi(\xb,t) = S(t) \partial^{\alphab}_\xb g.
\end{equation}
This proves the first statement provided in the lemma. 

Since $S(t) \in \sL\p{H}$ is the solution operator to the homogeneous Neumann initial boundary value problem, we have that \eqref{eqn: thermal semigroup commutation} is uniformly bounded for every $\xb \in \Omega$ only $g \in N^k\p{\Omega}$, e.g. if all partial derivatives of $g$ up to order $k$ satisfy the higher order Neumann boundary compatibility condition. The second result of the Lemma follows from \eqref{eqn: thermal semigroup commutation} after observing that $S(t) \in \sL\p{H}$ and that $S(t)$ is $\sL\p{H}$-norm continuous with respect to $t\in\mathbb R_+$. 
\end{proof}

This result demonstrates that the smoothness of the solution to the homogeneous parabolic problem in the spatial dimension is dependent on whether or not its partial derivatives satisfy the Neumann boundary condition, irrespective of the smoothness of the domain geometry. This is, in essence, what has allowed us to obtain optimal error convergence in our two-dimensional thermal model problem posed on the unit square (see \S\ref{subsection: 2D thermal problem}). This is in contrast to what one would initially assume from the general regularity theory for parabolic partial differential equations, where lower solution regularity is expected if $\Gamma$ is not sufficiently smooth \cite[Section 7.1.3]{evans2022partial}.

We now move on to derive the main result of this subsection in the following. 

\begin{lemma} \label{lemma: heat semigroup error}
    Let $S(t) \in \sL(H)$ be the parabolic semigroup generated by $A: \mathcal D(A) \rightarrow H$, and $S_n(t) \in \sL(H)$ be its finite element approximation, then we have that there exists a constant $C \in \mathbb R_+$ so that
    $$
        \norm{S(t) - S_n(t)}_{\sL(N^{k+3}(\Omega), L^2(\Omega))} \leq Ch^{k+1}_n \p{1 + t}e^{-\alpha t},
    $$
    where $\alpha \in \mathbb R_+$ is the fundamental rate of dissipation of the parabolic semigroup and $k \in \mathbb N_+$ is the polynomial order of the finite element approximation space. 
    \begin{proof}
        Following the analysis presented in \cite[Theorem 7.1]{strang1974analysis}, we have the following basic inequality for the error estimate for the finite element approximation of the solution to linear parabolic equations
        $$
            \norm{z(\cdot, t) - z_n(\cdot, t)}_{L^2(\Omega)} \leq Ch_n^{k+1}\p{\norm{z(\cdot, t)}_{H^{k+1}(\Omega)} + e^{-\alpha t} \norm{g}_{H^{k+1}(\Omega)} + \int_0^t e^{\alpha(s - t)} \norm{\partial_s z(\cdot,s)}_{H^{k+1}(\Omega)} ds }.
        $$
        We begin by seeing that 
        $$
            \norm{z(\cdot, t)}_{H^{k+1}(\Omega)} = \norm{S(t) g}_{H^{k+1}(\Omega)} \leq e^{-\alpha t} \norm{g}_{H^{k+1}(\Omega)},
        $$
        as a consequence of applying Lemma \ref{lemma: ll1} to the definition of the $H^{k+1}(\Omega)$ norm and recalling that $g \in N^{k+3}(\Omega)$. Since $g \in N^{k+3}\p{\Omega}$, it then follows that
        $$
        \begin{aligned}
            \norm{\partial_s z(\cdot,s)}_{H^{k+1}(\Omega)} &= 
            \norm{\alpha \Delta z(\cdot, s) - \beta z(\cdot,s)}_{H^{k+1}(\Omega)} \\
            &=\norm{\alpha \Delta S(s)g - \beta S(s) g}_{H^{k+1}(\Omega)} \\
            & \leq Ce^{-\alpha s} \norm{g}_{H^{k+3}(\Omega)},
        \end{aligned}
        $$
        which follows again from applying Lemma \ref{lemma: ll1}, and equation \eqref{eqn: heat equation}. It then follows that 
        $$
        \norm{z(\cdot, t) - z_n(\cdot, t)}_{L^2(\Omega)} \leq Ch_n^{k+1}\p{1 + t}e^{-\alpha t} \norm{g}_{H^{k+3}(\Omega)}.
        $$
        We then arrive at the result of this Lemma by recalling that $z(\cdot, t) = S(t) g$ and that $z_n(\cdot, t) = S_n(t) g$, and then finally applying the definition of the operator norm. 
    \end{proof}
\end{lemma}

We have shown above that the classical finite element error estimate for parabolic equations can be used to derive uniform error estimates for the approximate semigroup in the $\sL\p{K,H}$ operator norm. We will derive a similar estimate for the approximation of hyperbolic semigroups in the following subsection.

\subsection{Finite Element Approximation of the Weakly Damped Wave Semigroup} \label{subsection: finite elment approximation of the weakly damped wave semigroup}
Here, we present results pertaining to the finite element approximation of the weakly damped wave semigroup. Recall that we have set $H:= H^1_0(\Omega) \times L^2(\Omega)$ with $\Omega:(-1,1)$. It is known that $\mathcal A: \mathcal D(\mathcal A) \rightarrow H$ generates an exponentially stable semigroup \cite{conti2022optimal}, which we denote $S(t) \in \sL(H)$. Letting $\mathbf z(\cdot, t) = S(t) \mathbf q$ for any $\mathbf q \in H$, we have that $\mathbf 
 z(x,t) = \br{v(x,t), w(x,t)}$ is the solution of the following weakly damped wave equation
\begin{equation} \label{eqn: wave}
\left\{
\begin{aligned}
    \partial_t v(x,t) &= w &\textrm{ in } \Omega \times \mathbb R_+\\
    \partial_t w(x,t) &= c^2\partial_x^2 v(x,t) - \gamma w(x,t) &\textrm{ in } \Omega \times \mathbb R_+ \\
    v(x,t) &= 0 &\textrm{ on } \Gamma \times \mathbb R_+ \\
    v(x,0) &= q(x) &\textrm{ in } \Omega \\
    w(x,0) &= r(x) &\textrm{ in } \Omega,
\end{aligned}
\right. 
\end{equation}
where we have denoted $\mathbf q := \br{q, r}$. Additionally, it is known that $\mathcal A_n: \mathcal D(\mathcal A_n) \rightarrow H_n$ also generates an exponentially stable semigroup $\widehat S_n(t) \in \sL(H_n)$. The domain and range of $\widehat S_n(t)$ can then be extended to $\sL\p{H}$ through an injection and projection process, i.e.,
$S_n(t) := \iota_H \widehat S_n(t) \pi_n$. Setting $\mathbf z_n(\cdot, t) = S_n(t) \mathbf q$, for any $\mathbf q \in H$, we have that $\mathbf  z_n(\cdot, \cdot) \in \sC\p{\mathbb R_+; H_n}$ satisfies the following set of weak equations
$$
\left\{
\begin{aligned}
    \frac{d}{dt}\dual{v_n(t), \phi_n} &= \dual{w_n, \phi_n} \quad &\forall \phi_n \in V^k_{n,0} \\
    \frac{d}{dt}\dual{w_n(t), \psi_n} &= -a(v_n, \psi_n) - \gamma \dual{w_n, \psi_n} \quad &\forall \psi_n \in V^k_{n,0} \\
    \dual{v_n(0), \phi_n} &= \dual{q, \phi_n} &\forall \phi_n \in V^k_{n,0}  \\
    \dual{w_n(0), \psi_n} &= \dual{r, \psi_n} &\forall \psi_n \in V^k_{n,0} 
\end{aligned}
\right.
$$
for all $t \in \mathbb R_+$. With the approximate semigroup defined, we now determine its pointwise convergence to $S(t) \in \sL\p{H}$ for all $t \in I$ in the following. 

\begin{lemma} \label{eqn: damped wave pointwise convergence}
    Let $S(t), S_n(t) \in \sL\p{H}$ be the weakly damped wave semigroup and its finite element approximation respectively, then we have that 
    $$
        \lim_{n\rightarrow \infty} \norm{S(t)\mathbf q - S_n(t)\mathbf q}_H = 0
    $$
    for all $t\in I$ and $\mathbf q \in H$.
\end{lemma}
\begin{proof}
    This proof will rely on Theorem \ref{theorem: Trotter-Kato}, where the convergence of the approximate resolvent operator implies convergence of the approximate semigroup.
    The resolvent operator of $S(t)\in \sL\p{H}$ is given as the $\sL\p{H}$ operator $R(\mathcal A, s) := \p{\mathcal A - s\mathcal I_H}^{-1}$, where $\mathcal I_H$ denotes the identity operator in $H$. Let now $\mathbf z \in H$ be defined by $\mathbf z:= \p{v, w}$ and $\mathbf q \in H$ be defined by $\mathbf q:= \p{q, r}$. It then follows that $\mathbf z = R(\mathcal A, s)\mathbf q$ is the solution to the following set of equations
    \begin{equation} \label{eqn: AAA}
    \left\{
    \begin{aligned}
        sv(x) - w(x) &= q(x) \quad &\textrm{ on } \Omega \\
        sw(x) - c^2 \partial_x^2 v(x) + \gamma w(x) &= r(x) \quad  &\textrm{ on } \Omega \\
        v(x) &= 0 \quad &\textrm{ on } \Gamma.
    \end{aligned}
    \right.
    \end{equation}
    By seeing that $w(x) = \frac{c^2}{s + \gamma} \partial_x^2 v(x) + \frac{1}{s+\gamma} r(x)$, we have that we can reduce \eqref{eqn: AAA} into the following elliptic boundary value problem
    \begin{equation*}
    \left\{
    \begin{aligned}
        -\frac{c^2}{s+\gamma} \partial_x^2 v(x) + sv(x) &= q(x) + \frac{1}{s+\gamma}r(x) \quad &\textrm{ on } \Omega \\
        v(x) &= 0 \quad &\textrm{ on } \Gamma.
    \end{aligned}
    \right.
    \end{equation*}
    We see that ellipticity of this equation is preserved as long as $s \in (-\alpha, \infty)$ for some $\alpha \in \mathbb R_+$ chosen sufficiently small. Since the domain $\Omega := [-1, 1]$ does not have any corners, we have from the elliptic regularity theory \cite[Chapter 6, Theorem 4]{evans2022partial} that
    \begin{equation} \label{eqn: reg estimate}  
        \norm{v}_{H^2(\Omega)} \leq C_s \p{\norm{q}_{L^2(\Omega)}  + \norm{r}_{L^2(\Omega)}},
    \end{equation}
    where $C_s \in \mathbb R_+$ is a positive constant that depends on $s \in (-\alpha, \infty)$. 
    
    We now prove convergence of the resolvent operator associated with $S_n(t) \in S_n(t)$. We begin by first defining the weak form of \eqref{eqn: AAA} by seeking an $\p{v,w} \in H^1_0(\Omega) \times L^2(\Omega)$ that satisfies
    \begin{equation} \label{eqn: weak AAA}
    \left\{
    \begin{aligned}
        s\dual{v, \phi} - \dual{w, \phi} &= \dual{q, \phi} \quad &\forall \phi \in L^2(\Omega) \\
        \p{s-\gamma}\dual{w, \psi} - c^2\dual{\partial_x v, \partial_x \psi} &= \dual{r, \psi} \quad &\forall \psi \in H^1_0(\Omega).
    \end{aligned}
    \right.
    \end{equation}
    It then follows that $R(\mathcal A_n, s) := \iota_H \p{\mathcal A_n - s\mathcal I_{H_n}}^{-1} \pi_n$ is the solution operator to the following weak problem, where we seek $(v_n, w_n)\in H_n$ that satisfies
    \begin{equation} \label{eqn: AAA FEM}
    \left\{
    \begin{aligned}
        s\dual{v_n, \phi_n} - \dual{w_n, \phi_n} &= \dual{q, \phi_n} \quad &\forall \phi_n \in V^k_{n,0} \\
        \p{s-\gamma}\dual{w_n, \psi_n} - c^2\dual{\partial_x v_n, \partial_x \psi_n} &= \dual{r, \psi_n} \quad &\forall \psi_n \in V^k_{n,0}.
    \end{aligned}
    \right.
    \end{equation}
    Choosing $\phi \in H^1_0(\Omega)$ and $\phi_n \in V^k_{n,0}$, and subsequently eliminating $\dual{w, \psi}$ from \eqref{eqn: weak AAA} and $\dual{w_n, \psi_n}$ from \eqref{eqn: AAA FEM} respectively yields the following variational equations
    \begin{subequations}
    \begin{equation} \label{eqn: CCC}
        Q_s\p{v, \phi} = \dual{q + \frac{1}{s+\gamma}r, \phi} \quad \forall \phi \in H^1_0(\Omega)
    \end{equation}
    \begin{equation} \label{eqn: DDD}
        Q_s\p{v_n, \phi_n} = \dual{q + \frac{1}{s+\gamma}r, \phi_n} \quad \forall \phi_n \in V^k_{n,0},
    \end{equation}
    \end{subequations}
    where we have defined the bilinear form
    $$
        Q_s(\phi,\psi) := \frac{c^2}{s+\gamma}\dual{\partial_x \phi, \partial_x \psi} + s\dual{\phi, \psi}
    $$
    for any $\phi, \psi \in H^1_0(\Omega)$. Choosing $\phi = \phi_n$ in \eqref{eqn: CCC} and then subsequently subtracting \eqref{eqn: DDD} allows us to see that the Galerkin orthogonality condition
    $$
        Q_s\p{v - v_n,\phi_n} = 0
    $$
    is satisfied for all $\phi_n \in V^k_{n,0}$. It then follows by the definition of the $H^1(\Omega)$ norm, the boundedness, and coercivity of $Q_s(\cdot,\cdot)$, that there exists two positive constants $C_1, C_2 \in \mathbb R_+$ so that 
    $$
    C_1\norm{\phi}_{H^1(\Omega)}^2 \leq Q_s(\phi,\phi) \leq C_2\norm{\phi}_{H^1(\Omega)}^2
    $$
    for any $\phi \in H^1_0(\Omega)$. With this, we finally have that
    $$
    \begin{aligned}
        \norm{v - v_n}^2_{H^1(\Omega)} &\leq Q_s(v - v_n, v - v_n) \\
            &= Q_s(v-v_n, v-\phi_n) \\
            &\leq C\norm{v-v_n}_{H^1(\Omega)} \inf_{\phi_n \in V^k_{n,0}} \norm{v-\phi_n}_{H_1(\Omega)} \\
            &\leq Ch_n \norm{v - v_n}_{H^1(\Omega)} \norm{v}_{H^2(\Omega)},
    \end{aligned}
    $$
    where have applied the inequality $\inf_{\phi_n \in V^k_{n,0}} \norm{v - \phi_n}_{H_1(\Omega)} \leq Ch_n\norm{v}_{H^2(\Omega)}$ \cite[Proposition 1.12]{ern2004theory}. And hence we arrive at
    \begin{equation} \label{eqn: v error}
        \norm{v-v_n}_{H^1(\Omega)} \leq C_s'h_n\p{\norm{q}_{L^2(\Omega)} + \norm{r}_{L^2(\Omega)}}, 
    \end{equation}
    after applying \eqref{eqn: reg estimate}. 
    
    Returning to \eqref{eqn: weak AAA} and \eqref{eqn: AAA FEM}, we observe that
    \begin{equation} \label{eqn: nonorthogonality}
        \p{w - w_n , \phi_n} = s\p{v - v_n, \phi_n} 
    \end{equation}
    for all $\phi_n \in V^k_{n}$. Using this, it then follows that
    $$
    \begin{aligned}
        \norm{w-w_n}_{L^2(\Omega)} &= \dual{w-w_n, w - w_n} \\
            &= \dual{w - w_n, w - \phi_n} + \dual{w- w_n, \phi_n - w_n} \\
            &= \dual{w - w-n, w - \phi_n} + s\dual{v - v_n, \phi_n - w_n} \\
            &\leq \norm{w - w_n}_{L^2(\Omega)} \inf_{\phi_n \in V^k_{n,0}} \norm{w - \phi_n}_{L^2(\Omega)} + s\norm{v - v_n}_{L^2(\Omega)}\p{ \norm{w - \phi_n}_{L^2(\Omega)} + \norm{w - w_n}_{L^2(\Omega)}  } \\
            &\leq Ch_n \norm{w}_{H^1(\Omega)} \norm{w-w_n}_{L^2(\Omega)} + 2s\norm{v - v_n}_{L^2(\Omega)}\norm{w-w_n}_{L^2(\Omega)} \\
            &\leq C_s'h_n\br{\p{\norm{v}_{H^1(\Omega)} + \norm{q}_{H^1(\Omega)}} + \p{\norm{q}_{L^2(\Omega)} + \norm{r}_{L^2(\Omega)}}} \norm{w - w_n}_{L^2(\Omega)},
    \end{aligned}
    $$
    where we have applied the inequality $\lim_{\phi_n \in V^k_{n,0}} \norm{w-\phi_n}_{L^2(\Omega)} \leq Ch_n \norm{w}_{H^1(\Omega)}$ and \eqref{eqn: v error}. And hence 
    \begin{equation}\label{eqn: w error}
        \norm{w-w_n}_{L^2(\Omega)} \leq C_s''h_n\p{\norm{q}_{H^1(\Omega)} + \norm{r}_{L^2(\Omega)}}, 
    \end{equation}
    where we have applied the elliptic regularity estimate \eqref{eqn: reg estimate} on $\norm{v}_{H^1(\Omega)}$ in the final inequality. Combining \eqref{eqn: v error} and $\eqref{eqn: w error}$ then allows us to see that
    \begin{equation} \label{eqn: overall error}
    \begin{aligned}
        \norm{R\p{\mathcal A, s}\mathbf q - R\p{\mathcal A_n, s} \mathbf q}_{H} &= \norm{\mathbf z - \mathbf z_n}_{H} \\
        &=\norm{v - v_n}_{H^1(\Omega)} + \norm{w - w_n}_{L^2(\Omega)}\\
        &\leq C_s'''h_n \p{\norm{q}_{H^1(\Omega)} + \norm{r}_{L^2(\Omega)}},
    \end{aligned}
    \end{equation}
    and hence $\lim_{n\rightarrow \infty} \iota_H \p{\mathcal A_n - s \mathcal I_{H_n}}^{-1} \pi_n \mathbf q = \p{\mathcal A - s \mathcal I_H}^{-1} \mathbf q$ for all $s\in (-\alpha, \infty)$. The result of the lemma then follows from Theorem \ref{theorem: Trotter-Kato}.
\end{proof}

We now introduce the following differential operator for notational convenience
$$
    \boldsymbol \partial_{j}\mathbf z := 
    \br{
    \begin{array}{cc}
    \partial_x^j    &   0 \\
    0   &   \partial_x^j
    \end{array}
    }
    \br{\begin{array}{c} v(x) \\ w(x) \end{array}}
$$
for any $\mathbf z := [v(x), w(x)]^T$. We then move forward to prove that $\boldsymbol \partial_j$ and $S(t)$ commute under necessary smoothness and boundary conditions. Recall that we have defined
$$
    D^k(\Omega) := \Br{v \in H^k(\Omega) : \partial_x^j z = 0 \textrm{ on } \Gamma, \quad j=0,\ldots, k}.
$$
\begin{lemma} \label{lemma: la2}
Let $S(t) \in \sL(H)$ be the weakly damped wave semigroup, then we have that
\begin{equation} \label{eqn: la1}
    \boldsymbol \partial_j S(t) \mathbf q := S(t) \boldsymbol \partial_j \mathbf q
\end{equation}
for all $\mathbf q \in D^{k+1}(\Omega) \times H^{k}(\Omega)$. Furthermore, we have that $\mathbf z(t) \in \sC\p{I; D^{k+1}(\Omega) \times H^{k}(\Omega)}$.
\end{lemma}
\begin{proof}
    First, we see that $S(t) \boldsymbol \partial_j \mathbf q$ is well defined since $\boldsymbol \partial_j \mathbf q \in H$ for any $\mathbf q \in H^{k+1}(\Omega) \times H^{k}(\Omega)$ from the definition of Sobolev spaces. This is necessary since $S(t) \in \sL(H)$. The remainder of the proof follows using the same argument as in the proof of Lemma \ref{lemma: ll1}.

    Next, taking the $j$-th order partial derivative with respect to $x$ on both sides of \eqref{eqn: wave} results in
    $$
    \left\{
    \begin{aligned}
        \partial_t \partial_x^j v(x,t) &= \partial_x^j w &\textrm{ in } \Omega \times \mathbb R_+\\
        \partial_t \partial_x^j w(x,t) &= c^2\partial_x^2 \partial_x^j v(x,t) - \gamma \partial_x^j w(x,t) &\textrm{ in } \Omega \times \mathbb R_+ \\
        \partial_x^j v(x,t) &= 0 &\textrm{ on } \Gamma \times \mathbb R_+ \\
        \partial_x^j v(x,0) &= \partial_x^j q(x) &\textrm{ in } \Omega \\
        \partial_x^k w(x,0) &= \partial_x^k r(x) &\textrm{ in } \Omega,
    \end{aligned}
    \right. 
    $$
    after interchanging the order of partial differentiation. This then allows to see that 
    \begin{equation} \label{eqn: wave semigroup commutation}
    S(t)\boldsymbol \partial_j\mathbf q = \boldsymbol \partial_j \mathbf z(t) = \boldsymbol \partial_j S(t) \mathbf q
    \end{equation}
    for any $\mathbf q \in H^{k+1}(\Omega) \times H^{k}(\Omega)$. Since $S(t) \in \sL\p{H}$ is the solution operator to the weakly damped wave equation with homogeneous boundary conditions on $v$, we must have that $q \in D^{k+1}(\Omega)$ in order for \eqref{eqn: wave semigroup commutation} make sense, e.g. so that the partial derivatives of $\mathbf q$ up to order $k+1$ are compatible with the imposed boundary condition. This proves the first assertion made by the Lemma. The second result follows from the fact that $S(t) \in \sL\p{H}$ for all $t \in \mathbb R_+$ and that $S(t)$ is $\sL\p{H} $norm continuous with respect to $t\in\mathbb R_+$.
\end{proof}

We now move on to derive an error estimate on the finite element approximation of the weakly damped wave semigroup in the following.
\begin{lemma} \label{lemma: wave semigroup error}
    Let $S(t) \in \sL(H)$ be the weakly damped wave semigroup generated by $\mathcal A: \mathcal D(\mathcal A) \rightarrow H$, and $S_n(t) \in \sL(H)$ be its finite element approximation, then we have that there exists a constant $C \in \mathbb R_+$ so that
    $$
        \norm{S(t) - S_n(t)}_{\sL(V, H^1(\Omega) \times L^2(\Omega))} \leq Ch_n^{k},
    $$
    where $\alpha \in \mathbb R_+$ is the fundamental rate of dissipation of the parabolic semigroup and $V:= D^{k+3}(\Omega) \times H^{k+2}(\Omega)$. 
    \begin{proof}
    In this analysis, we take $\mathbf q \in D^{k+3}(\Omega) \times H^{k+2}(\Omega)$ as any arbitrary vector and $\mathbf z(t) = S(t)\mathbf q$ with $\mathbf z(t) = [v(t), w(t)]$. Lemma \ref{lemma: la2} then implies that $\mathbf z(t) \in D^{k+3}(\Omega)\times H^{k+2}(\Omega)$. 
    
    From the analysis presented in \cite{doi:10.1137/0710073}, we have the following basic error estimate for the weakly damped wave equation
    $$
    \begin{aligned}
        \norm{v - v_n}_{H^1(\Omega)} + \norm{w - w_n}_{L^2(\Omega)} &\leq Ch_n^{k}\p{\norm{v}_{H^{k+1}(\Omega)} + \norm{w}_{H^{k+1}(\Omega)}} \\
        &+ C_\gamma h_n^{k+1}\int_0^t \norm{\partial_s w(s)}_{H^{k+1}(\Omega)} ds.
    \end{aligned}
    $$
    Seeing that $\partial_t w = c^2\partial_x^2 v - \gamma w$, we observe that
    $$
        \norm{\partial_s w(s)}_{H^{k+1}(\Omega)} \leq C\p{\norm{v}_{H^{k+3}(\Omega)} + \norm{w}_{H^{k+2}(\Omega)}},
    $$
    after applying the embedding inequality $\norm{v}_{H^{k+1}(\Omega)} \leq \norm{v}_{H^{k+2}(\Omega)}$. Applying Lemma \ref{lemma: la2} and the exponential stability of $S(t) \in \sL(H)$ allows us to see that
    $$
        \norm{v}_{H^{k+3}(\Omega)} + \norm{w}_{H^{k+2}(\Omega)} \leq Ce^{-\alpha t}\p{\norm{p}_{H^{k+3}(\Omega)} + \norm{q}_{H^{k+2}(\Omega)}}.
    $$
    Putting this all together, we observe that
    $$
    \norm{v - v_n}_{H^1(\Omega)} + \norm{w-w_n}_{L^2(\Omega)} \leq Ch_n^k\p{\norm{p}_{H^{k+3}(\Omega)} + \norm{q}_{H^{k+2}(\Omega)}}.
    $$
    We conclude the proof by seeing that 
    $$
        \norm{S(t)\mathbf q - S_n(t)\mathbf q}_{H^1(\Omega) \times L^2(\Omega)} = \norm{v-v_n}_{H^1(\Omega)} + \norm{w - w_n}_{L^2(\Omega)}
    $$
    and applying the definition of the $\sL\p{H}$ operator norm. 
    \end{proof}
\end{lemma}
\end{document}